\documentclass[12pt,centertags,oneside]{amsart}

\usepackage{amsmath,amstext,amsthm,amscd,typearea,hyperref,stmaryrd,mathabx}
\usepackage{amssymb}
\usepackage{a4wide}
\usepackage[mathscr]{eucal}
\usepackage{mathrsfs}
\usepackage{typearea}
\usepackage{charter}
\usepackage{pdfsync}
\usepackage[a4paper,width=16.2cm,top=3cm,bottom=3cm,lines=50]{geometry}

\numberwithin{equation}{section}

\newtheorem{theorem}{Theorem}[section]
\newtheorem{definition}[theorem]{Definition}
\newtheorem{proposition}[theorem]{Proposition}
\newtheorem{corollary}[theorem]{Corollary}
\newtheorem{lemma}[theorem]{Lemma}

\newtheorem{example}[theorem]{Example}

\newcommand{\cali}[1]{\mathscr{#1}}

\newcommand{\Tan}{\mathop{\mathrm{Tan}}\nolimits}
\newcommand{\Fin}{\mathop{\mathrm{Fin}}\nolimits}
\newcommand{\Neg}{\mathop{\mathrm{Neg}}\nolimits}

\newcommand{\supp}{{\rm supp}}

\newcommand{\ddc}{{dd^c}}
\newcommand{\ddcx}{{dd^c_x}}
\newcommand{\ddcy}{{dd^c_y}}
\newcommand{\dc}{{d^c}}
\newcommand{\dbar}{{\overline\partial}}
\newcommand{\ddbar}{{\partial\overline\partial}}

\newcommand{\ind}{{\bf 1}}

\renewcommand{\Re}{{\rm Re}}
\renewcommand{\Im}{{\rm Im}}

\newcommand{\Cc}{\cali{C}}

\newcommand{\Fc}{\cali{F}}
\newcommand{\Gc}{\cali{G}}

\newcommand{\Oc}{\cali{O}}
\newcommand{\Pc}{\cali{P}}

\newcommand{\Uc}{\cali{U}}

\newcommand{\A}{\mathbb{A}}
\newcommand{\C}{\mathbb{C}}
\newcommand{\D}{\mathbb{D}}
\newcommand{\E}{\mathbb{E}}
\renewcommand{\H}{\mathbb{H}}
\newcommand{\N}{\mathbb{N}}
\newcommand{\Z}{\mathbb{Z}}

\newcommand{\R}{\mathbb{R}}
\newcommand{\T}{\mathbb{T}}
\newcommand{\B}{\mathbb{B}}

\newcommand{\U}{\mathbb{U}}

\renewcommand{\S}{\mathbb{S}}
\renewcommand{\P}{\mathbb{P}}

\newcommand{\Ebf}{\mathbf{E}}

\renewcommand{\L}{{\mathcal L}}

\newcommand{\Sing}{\mathop{\mathrm{Sing}}\nolimits}

\newcommand{\CommaBin}{\mathbin{\raisebox{0.5ex}{,}}}

\title[]{Unique Ergodicity for foliations on compact K\"ahler surfaces}

\author{Tien-Cuong Dinh}
\address{Department of Mathematics, National University 
of Singapore, 10 Lower Kent Ridge Road, Singapore 119076. 
{\tt  http://www.math.nus.edu.sg/$\sim$matdtc} }
\email{matdtc@nus.edu.sg}

\author{Vi{\^e}t-Anh Nguy{\^e}n}
\address{Universit\'e de Lille, 
Laboratoire de math\'ematiques Paul Painlev\'e, 
CNRS U.M.R. 8524,  
59655 Villeneuve d'Ascq Cedex, 
France. {\tt http://www.math.univ-lille1.fr/$\sim$vnguyen}}
\email{Viet-Anh.Nguyen@univ-lille.fr}

\author{Nessim Sibony}
\address{Laboratoire de Math\'ematiques d'Orsay, Univ. Paris-Sud, CNRS, Universit\'e
Paris-Saclay, 91405 Orsay, France}

\address{and Korea  Institute for Advanced Study (KIAS), 85 Hoegiro, Dongdaemun-gu,
Seoul 02455, Republic of Korea}
\email{Nessim.Sibony@math.u-psud.fr}

\date{April 22, 2019}

\begin{document}


\begin{abstract}
Let $\Fc$ be a holomorphic  foliation  by  Riemann surfaces  on a compact  K\"ahler surface $X.$ 
Assume it is generic in the sense that all the singularities  are hyperbolic and that the foliation  admits no directed positive closed $(1,1)$-current.
Then there exists a  unique (up to  a  multiplicative  constant)  positive  $\ddc$-closed $(1,1)$-current directed  by $\Fc.$
This  is  a very strong ergodic property  of  $\Fc.$  
Our proof uses an extension of  the theory of densities to  a class of non-$\ddc$-closed  currents. A complete description of the cone of directed positive  $\ddc$-closed $(1,1)$-currents is also given when $\Fc$ admits directed positive closed currents.
\end{abstract}

\maketitle
\tableofcontents

\medskip\medskip

\noindent
{\bf MSC 2010:} Primary  37F75 - 37A

 \medskip

\noindent
{\bf Keywords:}  Singular holomorphic  foliation, $\ddc$-closed  currents, hyperbolic  singularity,  density of currents, tangent current.


\section{Introduction} \label{S:Intro}

Let $X$ be a compact K\"ahler surface endowed with a K\"ahler form $\omega$. 
Let $\Fc$ be  a  (possibly  singular) holomorphic  foliation on   $X$.
Recall that the foliation $\Fc$ is given by  an open covering $\{\U_j\}$ of $X$ and  holomorphic vector fields $v_j\in H^0(\U_j,\Tan(X))$ with isolated singularities (i.e. isolated zeros) such that
$$v_j=g_{jk}v_k\quad\text{on}\quad \U_j\cap \U_k$$
for some non-vanishing holomorphic functions $g_{jk}\in H^0(\U_j\cap \U_k, \Oc^*_X).$
Its leaves are locally integral curves of these vector fields.  The set of singularities  of $\Fc$ is  precisely the union of the zero  sets of these vector fields. This  set is  finite. 

Using rational vector fields, we see that projective complex surfaces admit large families of foliations.  
Foliations can be also given locally by a non-zero holomorphic $1$-form and the leaves are Riemann surfaces on which these
forms vanish. In the case of complex
dimension $2$ that we consider, these leaves always exist without any integrability condition, i.e. the Frobenius condition is always satisfied for  bi-degree reasons.

If a holomorphic vector field has an isolated zero at some point $p,$ we
say that the singularity $p$ is hyperbolic if the two eigenvalues of the linear
part of the vector field at $p$ have non-real quotient. According to Poincar\'e,
 if $p$
is such a singular point, then there are local holomorphic coordinates centered
at $ p$ such that the vector field has the form
$$\eta x_1{\partial\over  \partial  x_1}+ x_2{\partial\over  \partial  x_2}  \, \CommaBin$$
where $(x_1 , x_2 )\in\C^2$, $ \eta = a + ib$ with $a, b \in\R$  and $b \not= 0.$

In order to  develop an ergodic theory for foliations, in the Riemannian case,  L. Garnett \cite{Garnett} introduced the notion  of  {\it harmonic measures}
for nonsingular foliations which  are generalizations of the {\it foliation cycles} of Sullivan \cite{Sullivan}. According to Sullivan \cite{Sullivan}, the existence of
a positive closed current, directed by the foliation,  corresponds to the existence of measures on transversals, invariant by the holonomy maps.

In the complex case, it is more fruitful to consider rather  the formalism of {\it directed $\ddc$-closed  currents.} This permits to  use  the interplay between cohomological intersection and geometric
intersection. In the present article, we use the cohomological properties of tangent currents.

Recall that   $d$ and $\dc$ denote the  real  differential operators on $X$  defined by
$d:=\partial+\overline\partial,$  $\dc:= {1\over 2\pi i}(\partial -\overline\partial) $ so that 
$\ddc={i\over \pi} \partial\overline\partial.$ 
A positive $\ddc$-closed current $T$ of bi-dimension $(1, 1)$ is directed by the foliation $\Fc$
if $T \wedge \Omega = 0$ for every local holomorphic $1$-form $\Omega$ defining $\Fc .$ Let $\U$ be any
flow box of $\Fc$ outside the singularities and denote by $V_\alpha$ the plaques of $\Fc$ in
$\U$ parametrized by $\alpha$ in some transversal $\Sigma$ of $\U.$ 
On the flow box $\U,$ such a current has the form
\begin{eqnarray} \label{e:T-box}
T|_\U=\int_{\alpha\in\Sigma} h_\alpha [V_\alpha] d\mu(\alpha),
\end{eqnarray}
where $h_\alpha$ is  a positive harmonic function on  $V_\alpha,$ and $[V_\alpha]$ denotes the current of integration on the plaque $V_\alpha$  (see e.g. \cite[Prop.\,2.3]{DinhNguyenSibony12}).
In  \cite{BerndtssonSibony} it is  shown that  for a foliation $\Fc$ by Riemann surfaces  
 with   finitely  many  singular points as above, there exists a non-zero directed  positive $\ddc$-closed  current.  
If $T$ is a positive $\ddc$-closed current of bi-dimension 
$(1, 1)$ directed by $\Fc ,$ then it has no mass on the singularities of $\Fc$
because this set is finite, see e.g. \cite{BerndtssonSibony,Skoda}. 

One of our main results gives the unique ergodicity for foliations
$\Fc$ which do not admit a positive directed closed current. This hypothesis implies that there are no invariant closed curve, 
and that $\Fc$ is {\it hyperbolic}, i.e. the leaves are hyperbolic or equivalently uniformized by the unit disc,  see  \cite{BurnsSibony}.  Unique ergodicity for the case where there is an invariant closed curve was studied in \cite{DinhSibony18}.

Now  we  briefly discuss  the family of holomorphic  foliations  on $\P^2$ with a given degree $d>1$. Foliations on $\P^2$ are always singular.  Recall that {\it the (geometric) degree} $d$ here is the number of tangencies of the foliation with a generic line.
 This  family can be identified  with  a Zariski dense open set $\Uc_d$ of some projective space.
 We will  say that   a  property is {\it  typical} for this family  if it is valid for  $\Fc$ in  a  set of full Lebesgue measure
 of $\Uc_d.$
 Here  are some typical properties of a  foliation in $\Uc_d$, see also Ilyashenko--Yakovenko \cite{IY}, Shcherbakov \cite{Sh} and \cite{SW}. 
 \begin{enumerate}
  \item {\rm  (Jouanolou \cite{Jouanolou} and  Lins Neto-Soares \cite{NetoSoares})}     
  all  the singularities of $\Fc$ are    hyperbolic and  $\Fc$ does not possess any invariant algebraic  curve.
  \item {\rm   (Glutsyuk  \cite{Glutsyuk} and  Lins Neto \cite{Neto})}
 $\Fc$ is hyperbolic.
 \item  {\rm  (Brunella \cite{Brunella})}    $\Fc$ admits no directed positive closed current.
  \end{enumerate}

Let  $\Fc$ be a hyperbolic foliation in a compact complex  manifold. Denote by $L_x$ the leaf of $\Fc$ through a point $x$.
Forn{\ae}ss and the
third author in \cite{FornaessSibony05} introduced
an average on each leaf $L_x$   which allows us to get another construction of directed
positive $\ddc$-closed currents.

More precisely, let $\D$ and $r\D$ denote the unit disc and the disc of center $0$
and radius $r$ in $\C.$ Let $\phi^x :\ \D \to L_x$ be a universal covering map for $L_x$ with
$\phi^x (0) = x.$ Define the Ahlfors-Shimizu characteristic function for $\phi^x$ by
$$
T^x(r):=\int_0^r {dt\over t}\int_{t\D} (\phi^x)^*(\omega),
$$
where we  recall that $\omega$ is a fixed K\"ahler  form on $X.$ Define the Nevanlinna current of
index $r$, $0 < r < 1$, associated with $L_x$ by
\begin{equation}\label{e:tau^x_r}
\tau^x_r:={1\over  T^x(r)} (\phi^x)_*\left\lbrack \log^+{r\over |\zeta|}  \right\rbrack={1\over  T^x(r)}\int_0^r {dt\over t}(\phi^x)_*[t\D].
\end{equation}
Here, $\log^+ := \max(\log, 0)$ and $\zeta$ is the standard coordinate of $\C$ so that the
unit disc $\D$ is equal to $\{|\zeta | < 1\}.$ Note that for each $x,$ the map $\phi^x$ is uniquely
defined up to a rotation in $\D.$ So the above definitions do not depend on the
choice of $\phi^x .$

 When  the singularities  of $\Fc$ are  all isolated (not necessarily hyperbolic), it was shown in \cite{FornaessSibony05} (see also \cite{DinhSibony18}) that $T^x (r ) \to\infty$ 
 as $r \to 1$ (this result still holds on manifolds of higher dimension).  Consequently, the cluster points
of $\tau^x_r $ are all $\ddc$-closed currents directed by $\Fc .$ It turns out that a Birkhoff
type theorem implies that for a generic foliation all extremal directed positive $\ddc$-closed currents of mass $1$ can be obtained in this way \cite{DinhNguyenSibony12}. General directed positive $\ddc$-closed currents are averages of
the extremal ones.

Here are the main results of the present paper which also hold for bi-Lipschitz laminations by Riemann surfaces (without singularities) in $X$. Recall that such a  lamination is a compact subset of $X$ which is locally a union of disjoint graphs of holomorphic functions depending in a bi-Lipschitz way on parameters, see Subsections \ref{SS:positive_ddc_closed_currents} and  \ref{SS:outside} for a precise local description.

\begin{theorem} \label{T:main_1}
 Let $\Fc$ be a holomorphic foliation by Riemann surfaces  with   only hyperbolic singularities or
 a bi-Lipschitz lamination by Riemann surfaces in a compact K\"ahler surface $(X,\omega)$. 
Assume that $\Fc$ admits no directed positive closed current.
Then there exists a  unique  positive  $\ddc$-closed  current $T$ of mass  $1$  directed  by $\Fc.$  
In particular, if $\phi^x :\ \D \to L_x$ is a universal covering map of an arbitrary   leaf $L_x$
as above, then $\tau_r^x \to T,$ in the sense of currents, as $r \to 1$. Moreover, the cohomology class $\{T\}$ of $T$ is nef and big, i.e. it belongs to the closure of the K\"ahler cone of $X$ and can be represented by a strictly positive closed $(1,1)$-current.
\end{theorem}

Note that the current $T$ is necessarily extremal in the cone of all positive $\ddc$-closed currents on $X$. Indeed, if $T'$ is such a current and $T'\leq T$, then $T'$ is necessarily directed by the foliation and according to the theorem, $T'$ is proportional to $T$. Note also that the nef property of $\{T\}$ is a consequence of a general result of independent interest, see Corollary \ref{c:nef} below. That corollary is a byproduct of our theory of densities of currents.

When  $X=\P^2$  the theorem  was  proved by Forn{\ae}ss and  the third author in  \cite{FornaessSibony10}. In that case  according  to \cite{Brunella}, if  all the singularities of $\Fc\in\Uc_d$ are hyperbolic
and $\Fc$ does not possess any invariant algebraic curve, then $\Fc$ admits no
directed positive closed current. So  the conclusion of Theorem \ref{T:main_1} is a typical property of the family $\Uc_d.$
The proof in  \cite{FornaessSibony10} is  based on two   ingredients.  The first one is an energy theory for positive $\ddc$-closed currents which was previously developed in \cite{FornaessSibony05}.
The second one is  a geometric  intersection  calculus for  these  currents. For the  second ingredient,
 the  transitivity of the automorphism group of $\P^2$ is heavily  used. Moreover, the proof is quite technical. The computations needed to estimate the geometric intersections are quite involved.
Using these techniques, P\'erez-Garrand\'es \cite{Perez-Garrandes}  has studied the case where $X$ is a homogeneous compact K\"ahler surface.

The  new  idea in the proof of  Theorem  \ref{T:main_1} is  to introduce 
a more flexible  tool 
which is   a  density theory for tensor products of  positive $\ddc$-closed currents.  The method allows us to bypass   the assumption of   homogeneity of $X.$    
The proof is  more conceptual and also  far less technical. The strategy is as follows. Given a positive $\ddc$-closed current $T$ on a surface $X$, we consider 
the positive current $T\otimes T$ near the diagonal  $\Delta$ of $X\times X$ which, in general, is not $\ddc$-closed.
We study the tangent currents to $T\otimes T$ along the diagonal $\Delta$. 
As one can expect this is related to the self-intersection properties of the current $T$. It turns out that the geometry of the tangent currents is quite simple. 
They are positive closed currents and are the pull-back of positive measures $\vartheta$ on $\Delta$ to the normal bundle of $\Delta$ in $X\times X$.
We relate the mass of $\vartheta$ to a cohomology class of the current $T$ and its energy.

The foliation or lamination enters in the picture to prove that $\vartheta$ is zero when $T$ is directed by a foliation or lamination as above. This is done using the local properties of the foliation or lamination,
the local description of $T$ and in particular, that the singularities are hyperbolic.
The vanishing of $\vartheta$ gives easily the uniqueness using a kind of Hodge-Riemann relations.

We expect that our results   could  have  numerous applications.
Using Theorem  \ref{T:main_1},
the  second  author has  very recently  shown in \cite{Nguyen18c}  that under the  assumption of this  theorem with the extra assumption that $X$ is projective, the Lyapunov exponent of $\Fc$ defined  in \cite{Nguyen17,Nguyen18b}
is  strictly negative. Moreover, when $X=\P^2$ the Lyapunov exponent of a  typical  foliation $\Fc\in\Uc_d$  is equal to  $-{d+2\over  d-1}\cdot$ 
The following result gives us a more complete picture of the strong ergodicity obtained in the present study.

\begin{theorem}\label{T:main_1bis} Let $\Fc$    be  
a    holomorphic foliation by Riemann surfaces  with   only hyperbolic singularities or a bi-Lipschitz lamination  by Riemann surfaces in a 
   compact K\"ahler surface $(X,\omega)$. Then one and only one of the following three possibilities occurs.
\begin{itemize}
\item[(a)] $\Fc$ admits invariant closed analytic curves and all positive directed $\ddc$-closed $(1,1)$-currents are linear combinations, with non-negative coefficients, of the currents of integration on those curves. In particular, these currents are all closed.
\item[(b)] $\Fc$ admits a directed positive closed $(1,1)$-current $T$ of mass $1$ having no mass on invariant closed analytic curves (this property holds when there is no such a curve). Every directed positive $\ddc$-closed $(1,1)$-current is closed, and if it has no mass on invariant closed analytic curves, then it has no mass on each single leaf and its cohomology class is proportional to $\{T\}$. Moreover, $\{T\}$ is nef (i.e. it belongs to the closure of the K\"ahler cone of $X$) and $\{T\}^2=0$.
\item[(c)] $\Fc$ admits a unique 
directed positive $\ddc$-closed and non-closed $(1,1)$-current $T$ of mass $1$ having no mass on each single leaf. 
Every directed positive $\ddc$-closed $(1,1)$-current is a combination, with non-negative coefficients, of $T$ and the currents of integration on invariant closed analytic curves. 
Moreover, $\{T\}$ is nef and big.
\end{itemize}   
\end{theorem}

A polynomial vector field in $\C^2$ induces a holomorphic foliation in $\P^2$.
When we fix the maximum of the degrees of its coefficients, if the vector field is generic,
the line at infinity $L_\infty:=\P^2\setminus\C^2$ is an invariant curve, see Ilyashenko-Yakovenko \cite{IY}. 
The current $[L_\infty]$ is the only directed positive $\ddc$-closed $(1,1)$-current of mass 1. So Property (a) holds in that case, see \cite{DinhSibony18} for details and also Rebelo \cite{Rebelo} 
for a related result. Note also that when Property (a) holds, a general theorem by Jouanolou says that there are only finitely many invariant closed analytic curves \cite{Jouanolou78}.

If $\Fc$ is a smooth fibration on $X$, then the directed positive $\ddc$-closed currents are all closed and are generated by the fibers of $\Fc$. They belong to the same cohomology class which is nef with zero self-intersection. So Property (b) holds in that case. Using a suspension one can also construct examples satisfying Property (b) which are not fibrations, see \cite[Ex.\,1]{FSW} and replace the circle there by $\P^1$. In such examples, there are two invariant closed curves and infinitely many directed positive closed $(1,1)$-currents of mass 1 having no mass on those curves.

Property (c) holds for foliations which are, in some sense, generic. There are many examples of such foliations in $\P^2$ without invariant closed analytic curves. The cohomology class of the unique directed $\ddc$-closed $(1,1)$-current here is K\"ahler because $H^2(\P^2,\R)$ is of dimension 1. 
If we blow up the singularities of the foliation, we get examples satisfying the same property and having invariant closed analytic curves. Then, the cohomology class of the unique directed $\ddc$-closed $(1,1)$-current is no more K\"ahler but it is big. In fact, we have the following general result which is a direct consequence of Theorem \ref{T:main_1bis}.

\begin{corollary}
Let $\Fc$ be a holomorphic foliation by Riemann surfaces  with   only hyperbolic singularities or a bi-Lipschitz lamination  by Riemann surfaces in a 
   compact K\"ahler surface $X$. Let $T$ be a positive $\ddc$-closed current directed by $\Fc$ having no mass on invariant closed analytic curves. Then the following properties are equivalent : 
   
\ \     (1) \ $T$ is not closed; \qquad (2) \ $\{T\}$ is big;  \qquad (3) \ $\{T\}^2>0$; \quad and  \quad (4) \ $\{T\}^2\not=0$.
\end{corollary}
 
Note that the hyperbolicity of the singularities is necessary in this result. The foliation on $\P^2$, given on an affine chart by the holomorphic 1-form 
$x_2dx_1-ax_1dx_2$ with $a\in\R$,  admits a non-hyperbolic singularity at 0 as well as diffuse invariant positive closed $(1,1)$-currents whose cohomology classes are K\"ahler. See also Corollary \ref{c:nef} and Theorem \ref{T:auxiliary} below which apply for foliations with arbitrary singularities.  
 
The  paper is  organized as follows. In  Section   \ref{S:Densities_and_Main_Theorem}, we introduce the densities for the tensor product of positive $\ddc$-closed currents using a notion of tangent current which is  described in  Theorem  \ref{T:main_2}.
Then we state Theorem  \ref{T:main_3} dealing with  the tensor square  power of a   
positive $\ddc$-closed current  directed by a  foliation or a lamination.
These are the key ingredients in the proofs of  the main theorems 
which will be presented at the end of this  section.   The proof  of Theorem  \ref{T:main_2}  occupies   Section 
\ref{S:tangent_currents}. Section \ref{S:Vanishing} is devoted to the proof of Theorem  \ref{T:main_3}. In Appendices \ref{a:Young}, \ref{a:current}, \ref{a:Harnack}, we present some basic facts on Young's inequality, $\ddc$-closed currents, directed $\ddc$-closed currents, Harnack's  inequality and their consequences that we use in the previous sections.

Note that after  we had finished  the article, Deroin informed us that with Kleptsyn, they had independently obtained a result similar to our first main theorem under  stronger hypotheses on the foliation and on the surface.

\smallskip

\noindent
{\bf Main Notation.} For the reader's convenience, we list here the main notations which are used through the paper. We consider a compact K\"ahler surface $(X,\omega)$ and denote by $\Fc$ a foliation by Riemann surfaces or a bi-Lipschitz lamination without singularities on $X$. Denote by $\Pi:\widehat{X\times X}\to X\times X$ the blow-up along the diagonal $\Delta$ of $X\times X$ and $\widehat\Delta:=\Pi^{-1}(\Delta)$ the exceptional hypersurface. The K\"ahler form $\widehat\omega$ on $\widehat{X\times X}$, the negative quasi-potential $\phi$ of $\Pi_*(\widehat\omega)$  will be chosen in Subsection \ref{SS:Mass}.
Denote by $\pi_j:X\times X\to X$ the projection onto the $j$-th factor and we use the K\"ahler form $\widetilde\omega:=\pi_1^*(\omega)+\pi_2^*(\omega)$ on $X\times X$. The constants $c$ and $c_j$ that we will use depend only on the above choices of $\omega, \widehat\omega, \phi$ and some other auxiliary parameters. 

Let $\D$ and $r\D$ denote respectively the unit disc and the disc of center 0 and radius $r$ in $\C$, $\B$ and $r\B$ the unit ball and the ball of center 0 and radius $r$ in $\C^2$. When we use local coordinates $x=(x_1,x_2)$ (or $y=(y_1,y_2)$) on $X$, we often identify a chart of $X$ with $10\B=\{\|x\|<10\}$ and we work with a fixed finite covering of $X$ by open subsets of the form ${1\over 4}\B$. The diagonal $\Delta$ is then covered by a finite number of charts which are identified with ${1\over 4}\B \times {1\over 4}\B$; they are contained in the chart $10\B\times 10\B$. With the above local coordinates $x$ on $X$, denote also by $\B(x,r)$ the ball of center $x$ and of radius $r$. 

On the chart $10\B\times 10\B$, we use two local coordinate systems: the first system is the standard one $(x,y)=(x_1,x_2,y_1,y_2)$ and the second system is $(z,w):=(x-y,y)$ on which 
$\Delta$ is given by the equation $z=0$. The tangent bundles of $X\times X$ and $\Delta$ are denoted by $\Tan(X\times X)$ and $\Tan(\Delta)$. 
The normal vector bundle of $\Delta$ in $X\times X$ is denoted by $\E:=\Tan(X\times X)|_\Delta/\Tan(\Delta)$, where $\Delta$ is also identified to the zero 
section of $\E$. Denote by $\pi:\E\to\Delta$ the canonical projection. The fiberwise multiplication by $\lambda\in\C^*$ on $\E$ is denoted by $A_\lambda$. 
Over $\Delta\cap (5\B\times 5\B)$, with the coordinates $(z,w)$, $\E$ is identified to $\C^2\times 5\B$, $\pi$ is the projection $(z,w)\mapsto w$ and $A_\lambda$ 
is equal to the map $a_\lambda(z,w):=(\lambda z,w)$. 

The notations $\lesssim$ and $\gtrsim$ stand for inequalities up to a positive multiplicative constant. The pairing $\langle \cdot,\cdot\rangle$ often denotes the value of a current on a test form. This is often an integral on the manifold where the current is defined. We will also use some test forms which are smooth outside a point in $X$ or outside the diagonal $\Delta$ in $X\times X$. The paring $\langle \cdot,\cdot\rangle_0$ 
denotes an integral taken outside these singularities. 

Finally, several notations introduced in Appendix \ref{a:Harnack} are heavily used in Subsection \ref{SS:near_singularities}.
 
\smallskip

\noindent
{\bf Acknowledgments. }  
The first author is supported by the grants 
C-146-000-047-001 and R-146-000-248-114 from 
the National University of Singapore (NUS). 
The second author is supported by the Labex CEMPI (ANR-11-LABX-0007-01).
The paper was partially prepared 
during the visit of the second author at the Vietnam  Institute for Advanced Study in Mathematics (VIASM) and 
at the NUS. He would like to express his gratitude to these organizations for hospitality and  for  financial support.

\section{Theory of densities and strategy for the  proofs of the main theorems}  \label{S:Densities_and_Main_Theorem}

In this section, we will present the main tool used in this article: the theory of densities for a class of non $\ddc$-closed currents.
We refer the reader to \cite{DinhSibony18, DinhSibony18b} for the case of $\ddc$-closed currents. The proofs of the main theorems stated in the Introduction will be 
provided in this section modulo Theorems \ref{T:main_2} and \ref{T:main_3} whose proofs will be given respectively in Section \ref{S:tangent_currents} and Section \ref{S:Vanishing}.

\subsection{Tangent  currents of  tensor products of positive $\ddc$-closed currents}   \label{SS:positive_ddc_closed_currents}

Consider two positive $\ddc$-closed $(1,1)$-currents $T_1$ and $T_2$  on $X$.
We will study the density of $T_1\otimes T_2$ near the diagonal $\Delta$ of $X\times X$ via a notion of ``tangent cone'' to $T_1\otimes T_2$ along $\Delta$ that we introduce now.

\begin{definition}[see also \eqref{e:local_tau}, \eqref{e:local_dtau}, \eqref{e:local_dtau-1}]   \label{D:admissible-maps} 
\rm A {\it smooth admissible map} is a smooth bijective map $\tau$ from a
neighbourhood of $\Delta$ in $X\times X$ to a neighbourhood of $\Delta$ in $\E$ such that
\begin{enumerate}
\item The restriction of $\tau$ to $\Delta$ is the identity map on $\Delta$; in particular,
the restriction of the differential $d\tau$ to $\Delta$ induces a map from
$\Tan(X\times X ) |_\Delta$ to $\Tan(\E)|_\Delta$; since $\Delta$ is pointwise fixed by $\tau$, the differential $d\tau$ also induces two endomorphisms of $\Tan(\Delta )$ and $\E$  respectively;
\item The differential $d\tau (x,x),$ at each point $(x,x) \in \Delta,$ is a $\C$-linear map from the tangent space to $X\times X$ at $(x,x)$ to 
the tangent space to $\E$ at $(x,x)$;
\item  The endomorphism of $\E$, induced by $d\tau$ (restricted to $\Delta$), is the identity map.
\end{enumerate}
\end{definition}

Note that the dependence of $d\tau (x,x)$ in $(x,x) \in  \Delta$ is  in general not holomorphic.
Consider the exponential map from $\E$ to $X\times X$ with respect to any Hermitian metric on $X\times X$. It defines a smooth bijective map from a neighbourhood of $\Delta$ in $\E$ to  a neighbourhood of $\Delta$ in $X\times X$. The inverse map is smooth and admissible, see also \cite[Lem.\,4.2]{DinhSibony18b}.
 
Let $\tau$ be any smooth admissible map as
above. Define
\begin{equation*}
(T_1\otimes T_2)_\lambda  := (A_\lambda )_* \tau_*  (T_1\otimes T_2 ).
 \end{equation*}
This is a current of degree $4.$ Its domain of definition is some open subset of
$\E$ containing $\Delta$ which increases to $\E$ when $|\lambda|$ increases to infinity. 

Observe that $(T_1\otimes T_2)_\lambda$ is not a $(2, 2)$-current and we cannot speak of its positivity. Moreover,
it is not $\ddc$-closed in general and we cannot speak of its cohomology class.
The present situation is more involved than the case where $T_1$ and $T_2$ are closed because in this case the current
$(T_1\otimes T_2)_\lambda$ is also closed.

By \eqref{e:decompo_posi_har} from Appendix \ref{a:current}, we can write for $j\in\{1,2\},$
\begin{equation} \label{e:decompo_posi_har_T12}
T_j=\Omega_j+\partial S_j+\overline{\partial S_j} + i\ddbar u_j,
\end{equation}
where $\Omega_j$ is a closed real smooth $(1,1)$-form, $S_j$ is a current of bi-degree $(0,1)$ and $u_j$ is a real current of bi-degree $(0,0)$. Note that $\partial \overline S_j$ and $\dbar S_j$ are forms of class $L^2$ which are independent of the choice of $\Omega_j,S_j,u_j$.
It turns out that a crucial argument in the proof of Theorem \ref{T:main_2} below is a result on the regularity of the potentials $u_j$ and their gradients, see Proposition \ref{p:T-rep} in Appendix \ref{a:current}.

The following theorem will be proved in  Section \ref{S:tangent_currents}. We refer to Appendix \ref{a:current} for the notion of Lelong number $\nu(T_j,\cdot)$ and the energy $E(T)$. 

\begin{theorem}\label{T:main_2} 
Let $T_1$ and $T_2$ be  two positive $\ddc$-closed $(1,1)$-currents on a compact K\"ahler surface $(X,\omega)$.   Assume that $T_1$ has no mass on the set $\{\nu(T_2,\cdot)>0\}$ and $T_2$ has no mass on the set $\{\nu(T_1,\cdot)>0\}$. Then, with the above notations, we have the following properties.
\begin{enumerate}
\item The mass of $(T_1\otimes T_2)_\lambda$ on any given compact subset of $\E$ is bounded uniformly on $\lambda$ for $|\lambda|$ large enough. If $\T$ is a cluster value of $(T_1\otimes T_2)_\lambda $ when $\lambda\to\infty ,$ then it
is a positive closed $(2, 2)$-current on $\E$ given by $\T = \pi^* (\vartheta)$ for some
positive measure $\vartheta$ on $\Delta.$ Moreover, if $(\lambda_n )$ is a
sequence tending to infinity such that $(T_1\otimes T_2)_{\lambda_n}  \to \T,$ then $\T$ may depend on $(\lambda_n )$
but it does not depend on the choice of the map $\tau .$
\item
The mass of $\vartheta$ does not depend on the choice of $\T$ and it is given by
\begin{equation}\label{e:mass}
\|\vartheta\|= \int_X\Omega_1\wedge \Omega_2-  \int_X \dbar S_1\wedge \partial \overline{S}_2 - \int_X \dbar S_2\wedge \partial\overline{S}_1  .
\end{equation}
In particular, if $T_1=T_2=T$ with $T=\Omega+\partial S+\dbar\overline S+i\ddbar u$ as in \eqref{e:decompo_posi_har}, then  
\begin{equation}\label{e:mass_bis}
\|\vartheta\|= \int_X\Omega^2-  2\int_X \dbar S\wedge \partial \overline{S} =\int_X\Omega^2-2E(T)  .
\end{equation}
\end{enumerate}
\end{theorem}

Note that in general $\T$ is not unique as this is already the case for positive
closed currents, see \cite{DinhSibony18b} for details. However, the mass formula shows that if one of such currents is zero then all of them are zero.
We can now introduce the following notion.

\begin{definition} \rm  
Any current $\T$ obtained as in Theorem \ref{T:main_2} is called a {\it tangent
current} to $T_1\otimes T_2$ along the diagonal $\Delta .$
\end{definition}

We have the following result and refer to McQuillan \cite{McQuillan} and  Burns--Sibony \cite{BurnsSibony} for some related results in the foliation setting.

\begin{corollary} \label{c:nef}
Let $T$ be a positive $\ddc$-closed $(1,1)$-current of a compact K\"ahler surface $X$. Assume that the set $\{\nu(T,\cdot)>0\}$ is of Hausdorff $2$-dimensional measure $0$.
Then the cohomology class $\{T\}$ of $T$ is nef, and when $T$ is not closed, $\{T\}$ is also big. In particular, if $T$ is a positive closed $(1,1)$-current having no mass on proper analytic subsets of $X$, then $\{T\}$ is nef.
\end{corollary}
\proof
We consider the first assertion on the nefness of $\{T\}$.
Let $Z$ be any irreducible analytic subset of dimension 1 of $X$. Denote by $[Z]$ the positive closed $(1,1)$-current of integration on $Z$ and $\{Z\}$ its cohomology class. To prove the nefness, we only need  to check that $\{T\}^2\geq 0$ and $\{T\}\smallsmile\{Z\}\geq 0$, see Demailly-Paun \cite[Cor.\,0.3]{DP}. 

We first show that $T$ has no mass on $Z$. Let $T'$ denote the restriction of $T$ to $X\setminus Z$. Since $T'$ is positive $\ddc$-closed with finite mass, we can extend it by zero through $Z$ and we still denote by $T'$ the extended current. 
This current $T'$ is positive and we have $\ddc T'\leq 0$, see \cite{AB, DinhSibony07}. On the other hand, by Stoke's theorem, we have
$$\|\ddc T'\|=\langle -\ddc T', 1\rangle =\langle -T',\ddc 1\rangle =0.$$
It follows that $\ddc T'=0$. Therefore, $T-T'$ is a positive $\ddc$-closed current supported by $Z$. So it is equal to $h[Z]$ for some non-negative harmonic function $h$ on $Z$. By maximum principle, $h$ should be constant. If $h\not=0$, we see that $T$ has a positive Lelong number at each point of $Z$. This contradicts the hypothesis on $T$. So $h=0$ and we deduce that $T=T'$ or equivalently $T$ has no mass on $Z$.

Since  $\{\nu(T,\cdot)>0\}$ is of Hausdorff 2-dimensional measure $0$, we also deduce that $[Z]$ has no mass on $\{\nu(T,\cdot)>0\}$.
Therefore, we can apply Theorem \ref{T:main_2} to $T_1:=T$ and $T_2:=[Z]$. From \eqref{e:decompo_posi_har_T12}, since $[Z]$ is closed, we get that $\dbar S_2=0$, see the discussion after \eqref{e:energy}.  Hence
$$\{T\}\smallsmile \{Z\} =\int_X\Omega_1\wedge\Omega_2=\|\vartheta\|\geq 0.$$

Since  $\{\nu(T,\cdot)>0\}$ is of Hausdorff 2-dimensional measure $0$, $T$ has no mass on this set, see \cite{BerndtssonSibony}. By the last assertion in Theorem \ref{T:main_2}, since $E(T)\geq 0$, we also have 
$$\{T\}^2=\int_X\Omega^2=\|\vartheta\|+2E(T)\geq 0.$$
So $\{T\}$ is nef. Moreover, if $T$ is not closed, then $E(T)>0$, see the discussion after \eqref{e:energy}. Therefore, $\{T\}^2>0$ and hence $\{T\}$ is big, i.e. it can be represented by a strictly positive closed $(1,1)$-current, see Demailly-Paun \cite[Th.\,0.5]{DP}. This ends the proof of the first assertion.

For the second assertion, since $T$ is closed and has no mass on proper analytic subsets of $X$, by Siu's theorem, the set $\{\nu(T,\cdot)>0\}$ is countable, see \cite{Siu}. So we can apply the first assertion to such a current $T$.
Note that in this case, Demailly-Paun theorem implies that $\{T\}$ is not big if and only if $\{T\}^2=0$. The last property also implies that $T$ has no positive Lelong number.
\endproof

The following result gives us  the  vanishing  of the  tangent  currents in the setting of foliations and laminations. Its proof will be given in  Section \ref{S:Vanishing}.
   
\begin{theorem} \label{T:main_3}
Let $\Fc$    be  either a holomorphic foliation by Riemann surfaces  with   only hyperbolic singularities,
or a bi-Lipschitz lamination  by Riemann surfaces, in a compact K\"ahler surface $X$. 
Then for every  positive $\ddc$-closed  current $T$  directed  by $\Fc$ which  does not give mass to any  invariant closed analytic curve,  zero is  the unique  tangent  current  to $T\otimes T$  along  the diagonal $\Delta.$
\end{theorem}

Recall that if
a closed subset $Y$ of a complex manifold $X$ is laminated by Riemann
surfaces, then it admits an open covering ${\U_j}$ and on each $\U_j$ there is a homeomorphism $\varphi_j = (h_j , \lambda_ j ) :\ \U_j \cap Y \to  \D \times \T_j$,
where $\T_j$ is a locally compact  metric space and the maps $\varphi_j^{ -1}(z,t)$, with $(z,t)\in\D\times \T_j$, are holomorphic in $z$.
Moreover, on their domains of definition, the transition maps have the form
$$\varphi_k \circ \varphi^{-1}_j (z, t) = \big(h_{jk} (z, t), \lambda_{jk} (t)\big),$$
where $ h_{jk} (z, t)$ is holomorphic with respect to $z$ and $\lambda_{jk}(t)$ do not depend on $z$.
We can choose $\T_j$ as the intersection of a holomorphic disc with $Y$ and $\varphi_j$ such that its restriction to $\T_j$ is the canonical map from $\T_j$ to $\{0\}\times\T_j$. With this choice,
when all $\varphi_j(z,t)$ are bi-Lipschitz maps, we say that the lamination is {\it bi-Lipschitz}.

The last theorem expresses that the current $T\otimes T$ is not too singular along the diagonal of $X\times X$ as its density along the diagonal is  zero.


\subsection{Sketch of the proofs of the main theorems} \label{SS:Main_Theorem}

The following result holds in a more general setting but we only state it in the case we use, see also \cite{DinhSibony18, FornaessSibony05}. Here, we don't need to assume that the singularities of the foliation are hyperbolic.

 \begin{theorem} \label{T:auxiliary}
 Let $T$ be a positive $\ddc$-closed $(1,1)$-current, on a  compact K\"ahler surface $X$, which is directed by 
a  holomorphic foliation or by a bi-Lipschitz lamination by Riemann surfaces. 
\begin{enumerate}
\item[(a)] If $T$ has a positive mass on a leaf $L$, then $\overline L$ is a closed analytic curve and $\overline L\setminus L$ is contained in the set of singularities of the foliation. Moreover, we can write 
$T=T'+T_{\rm an}$, where $T'$ is a directed positive $\ddc$-closed $(1,1)$-current which is diffuse, i.e. having no mass on each single leaf, and 
$T_{\rm an}$ is a finite or countable combination,  with non-negative coefficients, of currents of integration on invariant closed analytic curves. 
\item[(b)] Assume that $T$ gives no mass to any invariant closed analytic  curve. Then $T$ is diffuse and its cohomology class $\{T\}$ is nef. Moreover, $\{T\}$ is also big when $T$ is not closed.
\end{enumerate}
 \end{theorem}
 \proof
(a) Let $T''$ be the restriction of $T$ to $L$. Then, on a flow box outside the singularities, $T''$ is defined by positive harmonic functions on plaques. So $T''$ is positive $\ddc$-closed outside the singularities of the foliation (in the case of a lamination, this set is empty). Since $T''\leq T$, the mass of $T''$ is finite. Hence, as in Corollary \ref{c:nef}, 
 one can extend it by zero to a positive $\ddc$-closed $(1,1)$-current on $X$ that we still denote by $T''$. As in \cite[Prop.\,2.6]{DinhSibony18}, we obtain that 
 $\overline L$ is a compact analytic curve, $\overline L\setminus L$ is contained in the set of singularities of the foliation and $T''=c[L]$ for some constant $c>0$. 

We define $T_{\rm an}$ as the restriction of $T$ to the union of leaves of positive mass. We have seen that these leaves are contained in invariant closed curves and we deduce from the above discussion that $T_{\rm an}$ is positive and closed.
Since the mass of $T$ is finite, this family of leaves is at most countable.  It is now enough to define $T':=T-T_{\rm an}$. Clearly, this is a directed positive $\ddc$-closed $(1,1)$-current which is diffuse.
 
\medskip

(b)  Assume now that $T$ gives no mass to any invariant closed analytic  curve. Clearly, $T$ is diffuse. It follows that $T$ has zero Lelong number at any point outside the singularities of the foliation, see also \eqref{e:Lelong}. By Corollary \ref{c:nef}, the cohomology class $\{T\}$ is nef and it is also big when $T$ is not closed. This ends the proof of the theorem.
 \endproof

 The first step of our proof  consists in proving the following lemma.
 
\begin{lemma}\label{L:closed}
Let $\Fc$ be either a holomorphic foliation  by  Riemann surfaces with only  hyperbolic  singularities, or a bi-Lipschitz
lamination  by Riemann surfaces   in  a compact K\"ahler surface $(X,\omega).$ 
Let $T_1$ and $T_2$ be two positive $\ddc$-closed currents of mass $1$ directed by $\Fc$  such that
neither of them gives mass to any  invariant closed analytic curve.
Then $T_1-T_2$ is   a closed current. If both $T_1$ and $T_2$ are closed, then we have $\{T_1\}^2=\{T_2\}^2=\{T_1\}\smallsmile\{T_2\}=0$.
\end{lemma}
\proof
Since   both $T_1 $ and  $ T_2$ do not give mass to any  invariant closed analytic curve, it follows from Theorem \ref{T:auxiliary} that  
$\nu(T_1,x)=\nu(T_2,x)=0$ for all $x$ outside the singularities of $\Fc.$ 
Since $T_1$ and $T_2$ do not  give mass to this finite set, we  see that $T_1$ and $T_2$  satisfy the assumption of Theorem \ref{T:main_2}.

By  \eqref{e:decompo_posi_har_T12} and Stokes' theorem, we have (the second integral is the mass of $T_j$ which is assumed to be 1)
\begin{equation}\label{e:Omega_omega}
 \int_X \Omega_j\wedge \omega= \int_X T_j\wedge \omega=1\quad\text{for}\quad j=1,2.
\end{equation}
Applying  Theorems \ref{T:main_2} and \ref{T:main_3} to  each one of the  three directed  positive $\ddc$-closed currents $T_1,$ $T_2$ and $T_1+T_2$,  we obtain that all
 $T_1\otimes T_1,$  $T_2\otimes T_2$ and   $(T_1+T_2)\otimes( T_1+T_2)$     admit  zero as the unique  tangent current along the  diagonal $\Delta.$
This, combined with \eqref{e:decompo_posi_har_T12} and \eqref{e:mass_bis}, implies that
\begin{equation}\label{e:Omega_square}
\begin{split}
 \int_X  \Omega_1^2=2\int_X\dbar S_1\wedge \partial\overline S_1, \quad&    \int_X  \Omega_2^2=2\int_X\dbar S_2\wedge \partial\overline S_2  \\
\text{and} \quad  \int_X  (\Omega_1+\Omega_2)^2&=2\int_X\dbar (S_1+S_2)\wedge \partial (\overline S_1+\overline S_2).
 \end{split}
\end{equation}
If both $T_1$ and $T_2$ are closed, we deduce from the discussion after \eqref{e:energy} 
that $\dbar S_1=\dbar S_2=0$ and hence all integrals in \eqref{e:Omega_square}
vanish. This implies $\{T_1\}^2=\{T_2\}^2=\{T_1\}\smallsmile\{T_2\}=0$ as stated in the second assertion of the lemma.

Let $T:=T_1-T_2,$ $\Omega:=\Omega_1-\Omega_2,$   $S:=S_1-S_2$ and $u:=u_1-u_2.$ We infer  from \eqref{e:decompo_posi_har_T12} and \eqref{e:Omega_omega} that 
\begin{equation} \label{e:decompo_posi_har_T}
T=\Omega+\partial S+\overline{\partial S}+i\ddbar u\quad\text{and}\quad\int_X \Omega\wedge \omega=0.
\end{equation}
Moreover, it follows from  \eqref{e:Omega_square} that
\begin{equation} \label{e:integral}
\int_X  \Omega^2 =    \int_X  (\Omega_1 -\Omega_2)^2 =   
2\int_X  \Omega_1^2+2\int_X  \Omega_2^2- \int_X  (\Omega_1+\Omega_2)^2\\
 = 2 \int_X\dbar S\wedge \partial\overline S.
 \end{equation}

On one hand, since $\dbar S$ is an $L^2$  $(0,2)$-form,  the  current $\dbar S\wedge \partial\overline S= \dbar S\wedge \overline{\dbar S}  $ is a positive measure. So the last integral in \eqref{e:integral} is non-negative and  it vanishes if only if $\dbar S=0$  almost everywhere.
On the other hand, since  we know by \eqref{e:decompo_posi_har_T}  that $\int_X \Omega\wedge \omega=0,$ the cohomology class  of $\Omega$ is  a primitive class of $H^{1,1}(X,\R).$ Therefore, it follows
from  the  classical Hodge--Riemann  theorem that  the first integral in \eqref{e:integral} is non-positive, see e.g. \cite{Voisin}.
We conclude that  $\dbar S=0$  almost everywhere.
This and \eqref{e:decompo_posi_har_T} imply that $dT=0.$  The proof of the lemma is thereby completed.
\endproof

\proof[End of the proof of Theorem \ref{T:main_1bis} (see also  \cite{FornaessSibony05})]
We only consider the case of a foliation because the case of a lamination can be obtained in the same way.
It is clear that not more than one property in the theorem holds.
By \cite[Th.\,1.4]{BerndtssonSibony}, there exists a positive $\ddc$-closed current $T_1$ directed by $\Fc$, see also \cite[Th.\,23]{FornaessSibony08}. We can assume that Property (a) in the theorem does not hold. So we can find a current $T_1$ of mass 1 which has no mass on each single leaf of $\Fc$, see Theorem \ref{T:auxiliary}. We show that either Property (b) or (c) holds.

\medskip\noindent
{\bf Case 1.} Assume that there is such a current $T_1$ which is not closed. We show that the foliation satisfies Property (c) in the theorem.
By Theorem \ref{T:auxiliary}, the class $\{T_1\}$ is nef and big. It remains to prove the 
uniqueness of $T_1$. Assume by contradiction that there is another positive $\ddc$-closed current $T_2$ of mass 1 directed by $\Fc$. 
If there is such a current which is closed, then we assume that $T_2$ is closed.
So we have 
$$\int_X T_1\wedge \omega=\int_XT_2\wedge \omega=1.$$
We need to find a contradiction.  

Consider a flow box away from the  set of singularities $\Sing(\Fc)$ of $\Fc$ that we identify with $\D\times \Sigma$. As in the Introduction, we have
$$T_j = \int_\Sigma h^\alpha_j [V_\alpha ]d\mu_j (\alpha),\quad j = 1, 2.$$
Let $\mu = \mu_1 + \mu_2$ and write $\mu_j = r_j \mu$ with a non-negative bounded function $r_j\in L^\infty(\mu)$.  Then we have
$$T_1 - T_2 = \int_\Sigma  \big(h^\alpha_1 r_1 (\alpha) - h^\alpha_2 r_2 (\alpha) \big)[V_\alpha]d\mu(\alpha).$$
Since we know by Lemma  \ref{L:closed} that $T_1 - T_2$ is a  closed  current, 
$h^\alpha_1 r_1 (\alpha) - h^\alpha_2 r_2 (\alpha)$ is constant, for $\mu$-almost every $\alpha$,  that we will denote by $c(\alpha)$. 

We decompose $c(\alpha)\mu(\alpha)$ on the space of plaques $\Sigma$ and obtain that
$ c(\alpha)\mu(\alpha) = \nu_1  '-\nu_2$ for
 mutually singular positive measures $\nu_1$ and $\nu_2$. Then
$$T_1 - T_2 = [V_\alpha ]\nu_1 (\alpha) -[V_\alpha ]\nu_2 (\alpha) = T^+ - T^-$$
for positive closed currents $T^\pm.$ These currents fit together to a global positive closed
currents on $X \setminus \Sing(\Fc).$ Observe  that the mass of $T^\pm$ is  bounded  by   the mass of $T_1+T_2$.  So  the mass of $T^\pm$ is bounded near $\Sing(\Fc)$. Since $\Sing(\Fc)$ is a finite set,
$T^\pm$ extend as positive closed currents through $\Sing(\Fc),$ see e.g. \cite{Sibony, Skoda}. 
Recall that positive $\ddc$-closed $(1,1)$-currents have no mass on finite sets. Therefore,
since we assumed above that $T_1\not=T_2$, we have either $T^+\not=0$ or $T^-\not=0$. It follows from our choice of $T_2$ that $T_2$ is closed and hence $T_1$ is closed as well. This is a contradiction which shows that such a current $T_2$ as above doesn't exist.

\medskip\noindent
{\bf Case 2.} Assume now that all directed positive $\ddc$-closed $(1,1)$-currents are closed. Consider arbitrary directed positive closed $(1,1)$-currents $T_1$ and  $T_2$ of mass 1 which are diffuse. 
So by Theorem \ref{T:auxiliary} applied to $T_1,T_2$, the classes $\{T_1\}$ and $\{T_2\}$ are nef. By Lemma \ref{L:closed}, we have
$\{T_1\}^2=\{T_2\}^2=\{T_1\}\smallsmile\{T_2\}=0$. 
We show that Property (b) in the theorem holds. It is enough to show that $\{T_1\}=\{T_2\}$. 

Since $T_1$ and $T_2$ are of mass 1, we have $(\{T_1\}-\{T_2\})\smallsmile \{\omega\}=0$. So $\{T_1\}-\{T_2\}$ is a primitive class in the Hodge cohomology group $H^{1,1}(X,\R)$ of $X$. By the classical Hodge-Riemann theorem, we have $(\{T_1\}-\{T_2\})^2<0$ unless 
$\{T_1\}-\{T_2\}=0$, see e.g. \cite{Voisin}. Using that $\{T_1\}^2=\{T_2\}^2=\{T_1\}\smallsmile\{T_2\}=0$, we deduce that $\{T_1\}=\{T_2\}$. This ends the proof of the theorem.
\endproof

\proof[End of the proof of Theorem \ref{T:main_1}]
We only consider the case of a foliation because the case of a lamination can be obtained in the same way.
By hypothesis, the foliation has no invariant closed analytic curve. Moreover, by
Theorem \ref{T:main_1bis}, Property (c) in that theorem holds. It follows that 
the foliation admits a unique directed positive $\ddc$-closed current $T$ of mass 1.
This current is not closed and $\{T\}$ is nef and big.
Since  every cluster point of $\tau^x_r$  as $r$ tends to 1 is  a positive  $\ddc$-closed current of mass $1,$
 $\tau^x_r$ converges necessarily to $T$  as $r$ tends to 1.   
\endproof

\section{Existence and  properties  of tangent  currents}\label{S:tangent_currents}

In this   section, we prove Theorem \ref{T:main_2}. In the  first subsection, 
we obtain some estimates which are important in our study. In the second subsection,
we prove the existence of tangent currents and explain how to compute tangent currents using local coordinates, see Proposition \ref{P:local_compu}. The proof of this proposition is given in the same subsection. Part (1) of  Theorem \ref{T:main_2} is a consequence of Proposition \ref{P:local_compu} and Lemma \ref{l:nu}. Part (2) of that theorem will be obtained in the last subsection.


 \subsection{Some test forms and mass estimates}\label{SS:Mass}
 
In this subsection, we will construct some special test forms and also give some estimates for positive $\ddc$-closed currents and  their tensor products. We have the following elementary lemma.

\begin{lemma}\label{l:identity}
 Let  $T_1$ and $T_2$ be as in \eqref{e:decompo_posi_har_T12}.
 Then for every closed smooth form $\Phi$ of bi-degree $(2,2)$ on $X\times X,$ we have
\begin{equation*}
 \langle T_1\otimes T_2, \Phi \rangle =\langle\Omega_1\otimes \Omega_2, \Phi \rangle- \langle\dbar S_1\otimes \partial \overline{S}_2, \Phi \rangle -\langle\partial \overline{S}_1\otimes \dbar S_2, \Phi \rangle.
 \end{equation*}
 If, moreover,  $\Phi$ is  $d$-exact, then
 \begin{equation*}
 \langle T_1\otimes T_2, \Phi \rangle = - \langle\dbar S_1\otimes \partial \overline{S}_2, \Phi \rangle -\langle\partial \overline{S}_1\otimes \dbar S_2, \Phi \rangle.
 \end{equation*}
\end{lemma}
\proof
Observe that when  $\Phi$ is  $d$-exact, since $\Omega_j$ are closed, by Stokes' theorem,  we get
$\langle\Omega_1\otimes \Omega_2, \Phi \rangle=0.$ Hence,  the last identity of the lemma  follows from the first one.
We prove now the first identity. 

Observe that $\Phi$ is $\partial$-closed and $\dbar$-closed. It follows that if $R$ is $\partial$-closed or $\dbar$-closed, by Stokes' theorem,  we have
$$\langle i\ddbar u_1\otimes R, \Phi\rangle =0 \quad \text{and}\quad  \langle R\otimes i\ddbar u_2, \Phi\rangle =0.$$
Therefore, from  \eqref{e:decompo_posi_har_T12}, we get
 \begin{equation*}
 \langle T_1\otimes T_2, \Phi \rangle =\langle\Omega_1\otimes \Omega_2, \Phi \rangle + \langle\partial S_1\otimes \dbar \overline{S}_2, \Phi \rangle + \langle\dbar \overline{S}_1\otimes \partial S_2, \Phi \rangle.
 \end{equation*}
On the other hand, by Stokes' formula again, we have 
$$\langle\partial  S_1\otimes \overline{\partial S}_2, \Phi \rangle  =  \langle \dbar\partial S_1\otimes \overline S_2,\Phi\rangle =-\langle\dbar S_1\otimes \partial \overline{S}_2, \Phi \rangle$$
and
$$\langle\overline{\partial S}_1\otimes \partial S_2, \Phi\rangle  =  \langle \partial \overline{\partial S}_1\otimes S_2, \Phi\rangle  = -\langle\partial \overline{S}_1\otimes \dbar S_2, \Phi \rangle.$$
Hence, the first identity in the lemma follows easily.
\endproof

By Blanchard's theorem  \cite{Blanchard},  $\widehat{X\times X}$ can be endowed with a K\"ahler form
$\widehat \omega.$  The current $\Pi_*(\widehat \omega)$ is positive closed and has positive Lelong numbers along $\Delta$ and is smooth  outside $\Delta.$ Multiplying $\widehat \omega$ by a positive constant allows us to  assume that  the Lelong number of $\Pi_*(\widehat \omega)$ along  $\Delta$ is equal to 1. 
So we have
\begin{equation} \label{e:omega-hat}
\Pi^*(\Pi_*(\widehat \omega))=\widehat \omega+[\widehat\Delta].
\end{equation}

Choose a quasi-psh function $\phi\leq -1$ on $X\times X$ such that $\ddc\phi-\Pi_*(\widehat\omega)$ is a smooth form. This function is smooth outside $\Delta$. Define $\widehat \phi:=\phi\circ\Pi$. We deduce from \eqref{e:omega-hat} that $\ddc\widehat\phi-[\widehat\Delta]$ is a smooth form.

Recall that we only work with a fixed finite atlas of $X$ as mentioned at the end of the Introduction.
Consider a chart $2\B\times 2\B$ in coordinates $(z,w)$ and cover $\Pi^{-1}(2\B\times 2\B)$ with two 
charts denoted by $\widehat \U_1$ and $\widehat \U_2$. The first one $\widehat\U_1$ is given  with local  coordinates 
$$(u,w)=(u_1,u_2,w_1,w_2) \quad \mbox{with} \quad  \|w\|<2 \quad \mbox{and} \quad |u_1|<2,|u_2|<2$$
such that
$$\Pi(u,w)=(u_1,u_1u_2,w_1,w_2)=(z_1,z_2,w_1,w_2).$$

\noindent
{\bf Note. } The second chart $\widehat \U_2$ is defined exactly in the same way, except that the map $\Pi$ is given there by 
$$\Pi(u,w)=(u_1u_2,u_2,w_1,w_2)=(z_1,z_2,w_1,w_2).$$
When we work with local coordinates near $\widehat\Delta$, we will only consider the chart $\widehat \U_1$. The case of $\widehat\U_2$ can be treated in the same way.

\medskip

The function $\phi$ and the forms $\Pi_*(\widehat\omega)$, $\Pi_*(\widehat\omega^2)$ are defined globally on $X\times X$. Their
singularities along $\Delta$ will play an important role in our study. Using local coordinates, we have the following lemma.

\begin{lemma} \label{l:Pi-omega}
There is a constant $c_1>0$ such that for $(z,w)\in 2\B\times 2\B$ we have
$$c_1^{-1}\widetilde\omega\leq \Pi_*(\widehat\omega) \leq c_1 (\ddc\log\|z\| +\widetilde\omega) \quad \text{and} \quad
\Pi_*(\widehat\omega^2) \leq c_1 (\ddc\log\|z\| \wedge \widetilde\omega +\widetilde\omega^2).$$
We also have the following estimates on $2\B\times 2\B$ and $(2\B\times 2\B) \setminus\Delta$ respectively
$$\Pi_*(\widehat\omega)-c_1\widetilde\omega \leq \ddc\phi\leq \Pi_*(\widehat\omega)+c_1\widetilde\omega \quad \text{and} \quad i\partial\phi\wedge \dbar\phi\leq c_1(\|z\|^{-2}\ddc\|z\|^2+\widetilde\omega).$$
\end{lemma}
\proof
Since $\Pi^*(\widetilde\omega)$ is a smooth form, it is bounded by a constant times $\widehat\omega$. This and \eqref{e:omega-hat} imply $c_1^{-1}\widetilde\omega\leq \Pi_*(\widehat\omega)$ for some constant $c_1>0$.
We use the coordinates $(u,w)$ on $\widehat \U_1$ as above. It is not difficult to see that 
$$\widehat\omega \lesssim \ddc\log(1+|u_2|^2)+\ddc |u_1|^2 +\ddc \|w\|^2.$$
This implies the first (double) inequality in the lemma by using the action of $\Pi_*$. 

We obtain the second inequality on $(X\times X)\setminus\Delta$ from the first one
by observing that $(\ddc\log\|z\|)^2=0$ outside $\Delta$. The inequality holds on $X\times X$ because $\Pi_*(\widehat\omega^2)$ has no mass on $\Delta$. To see the last point, one can observe that over each point of $\Delta$ the fiber is  a $\P^1$ and $\widehat\omega^2$ gives it zero mass.

The third (double) inequality is a direct consequence of the definition of $\phi$. It remains to prove the last inequality.
We will only check it on $\Pi(\widehat\U_1)$ because the same proof also works for $\Pi(\widehat\U_2)$.

Recall that $\widehat \phi:=\phi\circ\Pi$ and define $\widehat\psi:=\widehat \phi-\log |u_1|$. 
Since $\ddc\widehat\phi-[\widehat\Delta]$ is smooth and $\widehat\Delta$ is given by the equation $u_1=0$, we deduce that $\ddc\widehat\psi$ is smooth  on $\widehat\U_1$. It follows that
$\widehat \psi$ is  a smooth function on $\widehat\U_1$. Therefore, there are bounded functions $\widehat h,\widehat g_1$ and $\widehat g_2$ on $\widehat\U_1$ such that 
$$\partial \widehat\phi = {1\over 2  u_1} du_1 + \widehat h du_2 + \widehat g_1 dw_1 + \widehat g_2 dw_2.$$
Hence, if we define $h:=\widehat h\circ\Pi^{-1}$,  $g_1:=\widehat g_1\circ\Pi^{-1}$ and $g_2:=\widehat g_2\circ\Pi^{-1}$, we get
$$\partial \phi = {1\over 2 z_1} dz_1 + h d(z_2/z_1) + g_1 dw_1 + g_2 dw_2.$$

Now, using that  $|z_2|\leq 2|z_1|$ on $\Pi(\widehat \U_1)$, we get $\|z\|\lesssim |z_1|$ and we can find bounded functions 
$h_1$ and $h_2$ such that 
$$\partial \phi = \|z\|^{-1} (h_1 dz_1 + h_2 dz_2) + (g_1 dw_1 + g_2 dw_2).$$
Finally, by Cauchy-Schwarz inequality, we can bound $i\partial\phi\wedge\dbar\phi$ by
$$2 \|z\|^{-2} (h_1 dz_1 + h_2 dz_2) \wedge \overline{(h_1 dz_1 + h_2 dz_2)} + 2 (g_1 dw_1 + g_2 dw_2) \wedge \overline{(g_1 dw_1 + g_2 dw_2)}$$
and the desired inequality follows easily.
\endproof

In the following lemma, we only need to consider the integral of the term containing $dy_1\wedge d\overline y_1 \wedge dy_2\wedge  d\overline y_2$ because the other terms vanish on $\{x\}\times X$.

\begin{lemma}\label{L:Lelong_number}
 Let $T$ be a positive $\ddc$-closed  current of mass $1$ on $X.$ Then there exists a constant $c_2>0$, independent of $T$, such that for all $x\in X$,  we have 
$$ \int_{y\in X\setminus\{x\}}  T(y)\wedge \Pi_*(\widehat\omega)(x,y)\leq  c_2.$$
\end{lemma}
\proof 
Observe that  the intersection $\Pi_*(\widehat\omega) \wedge [\{x\}\times X]$ is a current and we can identify it with a positive closed $(1,1)$-form $S_x$ on $\{x\}\times X$ which is smooth outside $x$. Since the cohomology class of $[\{x\}\times X]$ is independent of $x$, the cohomology class of $S_x$ is also independent of $x$. The integral considered in the lemma is equal to
$$\int_{y\in X\setminus\{x\}}  T(y)\wedge S_x.$$
So it is enough to check that the last integral is bounded from above. 

Using a regularization of $\ddc$-closed currents with mass control \cite{DinhSibony04}, it is enough to consider the case where $T$ is smooth. The last integral is then equal to $\langle T, S_x\rangle$ and depends only on the cohomology classes of $T$ and of $S_x$. Since all these cohomology classes are bounded, the result follows easily.
\endproof

\begin{lemma}\label{L:T_1T_2log}
Let $T_1$ and  $T_2$ be two positive  $\ddc$-closed $(1,1)$-currents of mass $1$ on $X.$
Then 
 $$\big\langle T_1\otimes T_2,  \Pi_*(\widehat\omega) \wedge\widetilde\omega \big\rangle_0 \leq 2c_2.$$
\end{lemma}
\proof
We refer to the end of the Introduction for the definition of $\langle \cdot,\cdot\rangle_0$. Since $\widetilde\omega=\pi_1^*(\omega)+\pi_2^*(\omega)$, 
a bi-degree consideration  shows that the considered pairing is equal to
\begin{equation*}
\Big\langle T_2(y)\wedge\omega(y), \int_{x\in X\setminus\{y\}}  T_1(x)\wedge \Pi_*(\widehat\omega)(x,y)  \Big\rangle  + \Big\langle T_1(x)\wedge\omega(x), \int_{y\in X\setminus\{x\}}  T_2(y)\wedge \Pi_*(\widehat\omega)(x,y)  \Big\rangle.
\end{equation*}
On the other hand, by  Lemma \ref{L:Lelong_number},  the integrals in the last line are bounded by $c_2$ because, by hypothesis, the measures $T_1\wedge\omega$ and $T_2\wedge\omega$ have mass 1.
The lemma follows easily. 
\endproof

We will now construct a family of test forms $R_m$ and prove some estimates.
In the chart $\widehat \U_1$ as in the last subsection,  the hypersurface $\widehat\Delta$ is equal to $\{u_1=0\}$ and we have 
$\ddc \log |u_1|=[\widehat\Delta]$.
Moreover, since $\ddc(\phi\circ \Pi)-[\widehat\Delta]$ is a  smooth form, the  function
$\phi\circ\Pi-\log|u_1|$ is also smooth. We deduce that
$\phi-\log\|z\|$ is  bounded in $2\B\times 2\B.$ 
Choose a constant $M\gg 1$ large enough   such that $|\phi-\log{\|z\|}|\leq M$ on each chart $2\B\times 2\B$ of $X\times X$. 
 
 Let $\chi:\ \R\to\R$ be an increasing convex  smooth function  such that
 $\chi(t)=0$ for $t\leq -3M,$ $\chi(t)=t$ for $t\geq 3M,$ ${1\over 10M}\leq  \chi'(t)\leq 1,$ and $\chi''(t)\in \big[ {1\over 8M},{1\over 4M}\big]$  for  $t\in [-2M,2M].$ 
Fix also a constant $A\gg 1$ large enough.  Define for  $m\in\N$  
$$R_m:= \ddc[\chi (\phi+m)]+A\widetilde\omega.$$
 This is  clearly a  smooth closed  $(1,1)$-form on $X\times X.$ We first  show that it is  positive and has bounded  mass.
 A direct   computation gives
 \begin{equation}\label{e:R_m}
 R_m=\chi'(\phi+m)\ddc\phi+{1\over \pi} \chi''(\phi+m) i\partial \phi\wedge \dbar\phi+A\widetilde \omega.
 \end{equation}
 The second term is positive. 
 The first term is bounded below by 
 $-c_1 \widetilde\omega$, see Lemma \ref{l:Pi-omega}. We then deduce that $R_m$ is positive since $A$ is chosen large enough. Furthermore, since $R_m$ is cohomologous to $A \widetilde\omega,$ its mass is equal to the mass
of $A \widetilde\omega$  and hence is bounded independently of $m$.

\smallskip

We have the following lemmas. The goal is to understand the mass repartition of $T_1\otimes T_2$ near $\Delta$ and to prove the basic estimates given in Lemma \ref{L:T_1T_2-testforms}.

\begin{lemma}\label{L:R_m_annulus}
There is a constant $c_3>0$ such that the following properties hold.
\begin{enumerate}
\item For every integer $m\geq 0$, we have
 $$e^{2m}(idz_1\wedge d\overline{z}_1+idz_2\wedge d\overline{z}_2)\leq  c_3R_m
 \quad \text{on}\quad  \big\{e^{-m-1}\leq  \|z\|\leq  e^{-m},\ \|w\|< 2 \big\}. $$
 \item For each $0<r\leq 1,$ if  $m$ is the integer such that $e^{-m-1}<r\leq e^{-m}$, then 
 $$ir^{-2}(dz_1\wedge d\overline{z}_1+dz_2\wedge d\overline{z}_2)\leq  c_3\sum_{n=0}^\infty e^{-2n} R_{m+n}\quad\text{on}\quad \big\{0<\|z\|<r,\ \|w\|< 2\big\}.$$
\end{enumerate}
 \end{lemma}
\proof
(1) In the  considered domain,  we have $|\phi+m|\leq 2M$. Therefore, $\chi'(\phi+m)\geq {1\over 10M}$ and   
$\chi''(\phi+m)\in  \big[{1\over 8M}, {1\over 4M}\big].$
Define $\widehat \phi:=\phi\circ\Pi$ and $\widehat\psi:=\widehat \phi-\log |u_1|$. So $\widehat \psi$ is  a smooth function on $\widehat\U_1$ because $\ddc\widehat\psi$ is smooth. 
Observe that $|u_1|\leq \|z\|$ and hence $|u_1|^{-1}\geq e^m$  
on  the region $\Pi^{-1}\big\{e^{-m-1}\leq \|z\|\leq e^{-m},\ \|w\|\leq 2 \big\}$. We then obtain on the same region that the form $i\partial \widehat \phi\wedge \dbar\widehat \phi$ is  equal to
\begin{eqnarray*}
\lefteqn{i\partial (\widehat \psi+ \log |u_1|)\wedge \overline\partial (\widehat \psi+\log |u_1|)}\\
 &=& i\partial\Big[{M+1\over M}\widehat \psi+{M\over M+1}\log{|u_1|} \Big]\wedge\dbar \Big[ {M+1\over M}\widehat \psi+{M\over M+1}\log{|u_1|} \Big]  \\
 && -{2M+1\over M^2}i \partial\widehat \psi\wedge\dbar\widehat \psi+{2M+1\over (M+1)^2}i\partial \log{|u_1|}\wedge\dbar \log{|u_1|}\\
 &\geq& -{3\over M}i \partial\widehat \psi\wedge\dbar\widehat \psi+{1\over 4M}e^{2m} idu_1\wedge d\overline u_1 \mbox{ since the first term in the last sum is  positive}.
\end{eqnarray*}
Observe that the first term in the last line is bigger than $-\epsilon\widehat\omega$ for some small constant $\epsilon>0$ because $M$ is big. By Lemma \ref{l:Pi-omega}, we also have 
 $\Pi^*(\ddc\phi)\geq \widehat \omega-c_1\Pi^*(\widetilde\omega)$. Therefore,
for $A\gg 1$, using \eqref{e:R_m}, we have
$$\Pi^*(R_m)\geq  {1\over 200M^2} \big(e^{2m} idu_1\wedge d\overline u_1+\widehat \omega\big).$$

Recall that $e^m|u_1|\leq 1$ on  $\big\{e^{-m-1}\leq \|z\|\leq e^{-m},\ \|w\|\leq 2 \big\}$. So using that $z_1=u_1$ and $z_2=u_1u_2$, 
we can find a bounded function $\theta_0$ and bounded forms $\theta_j$  on the region $\Pi^{-1}\big\{e^{-m-1}\leq \|z\|\leq e^{-m},\ \|w\|\leq 2 \big\}$ such that
$$
\Pi^*\big(ie^{2m}(dz_1\wedge d\overline z_1+dz_2\wedge d\overline z_2)  \big)=e^{2m}\theta_0idu_1\wedge d\overline u_1 +e^mdu_1\wedge\theta_1+ e^m d\overline u_1\wedge\theta_2 +\theta_3 .
$$
By Cauchy-Schwarz inequality, the last sum is bounded above by $e^{2m}\theta_0'idu_1\wedge d\overline u_1+\theta_3'$ for some  bounded function $\theta_0'$ and bounded form $\theta_3'$. This, combined with
the previous estimate for $\Pi^* (R_m)$,
implies the inequality in (1) for a suitable constant $c_3$.

\smallskip

(2) Observe   that $r^{-2}\leq e^{2m+2}.$ Applying the first assertion for $m+n$ instead of $m$ yields the desired estimate for a suitable constant $c_3$.
 \endproof

 \begin{lemma}\label{L:R_m}
 Let $T_1$ and $T_2$ be two positive $\ddc$-closed   $(1,1)$-currents   of mass $1$ on $X.$ Then  there is a constant $c_4>0$, independent of $T_1,T_2$, such that 
  $$  \big\langle T_1\otimes T_2  , R_m\wedge  \widetilde\omega \big\rangle\leq  c_4\quad\text{for all}\quad  m\geq 1. $$
 \end{lemma}
 \proof
  Since $\chi''$ is  supported on  $[-3M,3M],$ we  see that  the factor in front of $i\partial \phi\wedge \dbar\phi$ in \eqref{e:R_m} is non-zero only if
 $|\phi+m|\leq 3M.$ Moreover, we know that $|\phi-\log{\|z\|}|\leq M.$ So the above factor is non-zero only if
 $|m-\log{\|z\|}|\leq 4M,$ that is,  $z$ belongs to the ring $ \{ e^{-m-4M}\leq\|z\|\leq e^{-m+4M} \}.$
Therefore,  it is enough to prove an estimate, similar to the one in the lemma, for an integral on a chart $\B\times\B$ as above because these charts cover a neighbourhood of $\Delta$ and hence the support of $R_m$ for $m$ large enough.

By Lemma \ref{l:Pi-omega}, outside the diagonal $\Delta$, we have
 $$i\partial \phi\wedge \dbar\phi\lesssim \|z\|^{-2} \ddc \|z\|^2 +\widetilde\omega.$$ 
 This, coupled  with  the expression of $R_m$ in \eqref{e:R_m},  implies that
 \begin{eqnarray*}
\lefteqn{ \big\langle T_1\otimes T_2  , R_m\wedge  \widetilde\omega \big\rangle_{\B\times\B} \ \lesssim \  \big\langle T_1\otimes T_2,  \widetilde\omega^2 \big\rangle + \big\langle T_1\otimes T_2 , \ddc \phi\wedge \widetilde\omega\big\rangle } \\
 &&+\int_{   e^{-m-4M}\leq\|x-y\|\leq e^{-m+4M} } \Big(T_1(x)\otimes T_2 (y)\Big) \wedge \|x-y\|^{-2}\ddc {\|x-y\|^2}\wedge \widetilde\omega,
 \end{eqnarray*}
 where we recall that $(z,w)=(x-y,y)$.
 
 It is clear that the first term in the last sum is equal to $\langle T_1,\omega\rangle\langle T_2,\omega\rangle=1$.
 By Lemma \ref{L:T_1T_2log}, the second term is also bounded. So it remains to check that the last term is  bounded by a constant independent of $T_1,T_2$ and $m.$

Setting $r:=e^{-m+4M},$ since  $\|x-y\| \approx e^{-m}$ and $\widetilde\omega=\pi_1^*(\omega)+\pi_2^*(\omega)$, the considered term is bounded above by a constant times
$$\int_{\|x\|<1} \Big( r^{-2}\int_{y\in \B(x,r)} T_2 (y) \wedge \ddcy {\|x-y\|^2} \Big) T_1(x)\wedge \omega(x)$$
$$ \quad \quad + \int_{\|y\|< 1}\Big( r^{-2}\int_{x\in \B(y,r)} T_1 (x) \wedge \ddcx {\|x-y\|^2} \Big) T_2(y)\wedge \omega(y)$$
 which is equal to 
$$ \int_{\|x\|<1}\nu(T_2,x,r) T_1(x)\wedge \omega(x)  + \int_{\|y\|<1}\nu(T_1,y,r)T_2(y)\wedge \omega(y).$$
Thus, the lemma follows from 
Lemma \ref{l:lelong} and the fact that both $T_1$ and $T_2$ have mass one.
 \endproof
 
 \begin{lemma}\label{L:R_m-R_n}
 Let $T_1$ and $T_2$ be two positive $\ddc$-closed  $(1,1)$-currents  of mass $1$ on $X.$ Then  there is a constant $c_5>0$, independent of $T_1$ and $T_2$, such that 
  $$ \big\langle T_1\otimes T_2  , R_m\wedge   R_n \big \rangle\leq  c_5\quad\text{for all}\quad m, n\geq 1.$$
 \end{lemma}
 \proof
 Since  $R_m\wedge R_n$ is a closed smooth form of bi-degree $(2,2)$ on $X\times X,$ it follows from Lemma \ref {l:identity} that
$\big\langle T_1\otimes T_2, R_m\wedge R_n \big\rangle$ is equal to 
$$\big\langle\Omega_1\otimes \Omega_2, R_m\wedge R_n \big\rangle- \big\langle\dbar S_1\otimes \partial \overline{S}_2, R_m\wedge R_n \big\rangle - \big \langle\partial \overline{S}_1\otimes \dbar S_2, R_m\wedge R_n \big\rangle.$$
Denote the three terms in the last sum by $I_1,I_2$ and $I_3$ respectively. We will show that they are bounded independently of $T_1,T_2,m$ and $n$. 

Since $\Omega_j$ is cohomologous to $T_j$ which is of mass 1, the cohomology class of $\Omega_j$ is bounded. The forms $R_m$ and $R_n$ are both cohomologous to $A^2\widetilde\omega$. Therefore, the integral $I_1$, which depends only on the cohomology classes of $\Omega_j,R_n$ and $R_m$, is clearly bounded.
 
In order to show  that  the sequences  $I_2$ and $I_3$  are bounded, we only need to prove that for every $L^2$ functions $f_1,f_2$ on $X$ and 
a bounded smooth $(2,2)$-form $\alpha$ on $X\times X$ : 
\begin{equation}\label{e:f_1f_2-L2}
\big|\langle (f_1\otimes  f_2)\alpha, R_m\wedge R_n\rangle \big|\leq  c\|f_1\|_{L^2}\|f_2\|_{L^2}\quad\mbox{for a constant $c$ independent of $m,n.$}
\end{equation}
We only need to consider the case where either $n$ or $m$ is big. 
Assume for simplicity that  $m$ is larger than a fixed constant large enough.
So $R_m\wedge R_n$ has support near the diagonal $\Delta$. Therefore, 
using a partition of unity, we can assume that both $f_1$ and $f_2$ have support in the same chart $\B$ as above. Since we can write $f_1,f_2$ as linear combinations of non-negative functions with bounded $L^2$ norm, we can assume that both $f_1$ and $f_2$ are non-negative. Moreover, since $\alpha$ can be written as a combination of bounded smooth positive $(2,2)$-forms, we can also assume that $\alpha$ is positive.

Observe that the factor before $i\partial\phi\wedge\dbar\phi$ in \eqref{e:R_m} vanishes outside  the region $W_m:=\{e^{-m+4M}\leq \|z\|\leq e^{-m+4M}\}$.
Using  \eqref{e:R_m}  and Lemma \ref{l:Pi-omega}, we obtain
$$R_m\lesssim \Pi_*(\widehat\omega) + \ind_{W_m} i\partial\phi\wedge \dbar\phi \quad \text{and similarly} \quad
R_n\lesssim \Pi_*(\widehat\omega) + \ind_{W_n}  i\partial\phi\wedge \dbar\phi .$$
Using these inequalities, Lemma \ref{l:Pi-omega} and the identity
 $\partial\phi\wedge\partial\phi=0$,  we obtain
\begin{eqnarray*}
R_m\wedge R_n &\lesssim& \Pi_*(\widehat\omega^2) +  \ind_{W_m} (i\partial\phi\wedge \dbar\phi)\wedge \Pi_*(\widehat\omega) +  \ind_{W_n} (i\partial\phi\wedge \dbar\phi)\wedge \Pi_*(\widehat\omega)  \\
& \lesssim& (\|z\|^{-2} + \ind_{W_m}\|z\|^{-4} + \ind_{W_n}\|z\|^{-4}) \widetilde\omega^2.
\end{eqnarray*}

Consider the  integral operator $P$ acting on forms on $\B\times\B$  with a suitable kernel $K(x,y)$ obtained from  the coefficients of the product of $\alpha$ with the  last sum.
Here,  we invoke Examples \ref{Ex:kernel_1} and \ref{Ex:kernel_2} from Appendix \ref{a:Young}  by taking into account  that $\|z\|=\|x-y\|$ and setting $r:=e^{-m+4M}$ or $r:=e^{-n+4M}$.
Applying Lemma \ref{L:Young} to $K$ for $\delta=0,$  we get $\|P(f_2)\|_{L^2}\lesssim \|f_2\|_{L^2}.$
Hence,
\begin{equation*}
\langle (f_1\otimes  f_2)\alpha, R_m\wedge R_n\rangle \lesssim \langle f_1, P(f_2)\rangle \lesssim \|f_1\|_{L^2}\|f_2\|_{L^2} .
\end{equation*}
 This completes the proof of \eqref{e:f_1f_2-L2}.
\endproof

\begin{lemma}\label{L:T_1T_2-testforms} 
Let $T_1$ and $T_2$ be two positive $\ddc$-closed  $(1,1)$-currents of mass $1$ on $X.$ Assume that $T_1$ has no mass on the set 
$\{\nu(T_2,\cdot)>0\}$ and $T_2$ has no mass on the set $\{\nu(T_1,\cdot)>0\}$. Then there is a constant $c_6 > 0$, independent of $T_1,T_2$, and for each $0 < r \leq 1$, there is a constant
$\epsilon_r > 0$ depending on $T_1,T_2$ such that $\epsilon_r \to 0 $ as $r \to 0$ and  the following
estimate holds. For any continuous function $f(z,w)$ with compact support in
$(r\B)\times \B$,  we have
$$\big|\big\langle T_1\otimes T_2, f \gamma\big\rangle\big| \leq  \|f\|_\infty \max(\epsilon_r r^k,c_6r^4).$$ 
Here, $\gamma$ is the wedge-product of four $1$-forms among $dz_1, dz_2, dw_1,dw_2$ or their  complex conjugates, and
$k$  is the total degree of $dz_1,dz_2, d\overline{z}_1, d\overline z_2$ in  $ \gamma$. 
\end{lemma}
\proof
Note that for a bi-degree reason, the pairing in the lemma vanishes unless $\gamma$ is of bi-degree $(2,2)$. 
Since the real and imaginary parts of $f$ can be written as differences of bounded non-negative functions, we can assume that $f$ is a non-negative real-valued function. For simplicity, we can also assume that $\|f\|_\infty=1$.
We distinguishes 5 cases according to the value of $k$.

\medskip\noindent
{\bf Case 1.} Assume that $k=0$ and hence $\gamma=\pm dw_1 \wedge d\overline w_1 \wedge dw_2 \wedge d\overline w_2$.
Observe that positive $\ddc$-closed $(1,1)$-currents on $X$ have no mass on finite sets. Then, 
by applying Fubini's theorem, we obtain that  $T_1\otimes T_2$ has no mass on $\Delta$. 
Therefore, the positive measure $(T_1\otimes T_2)\wedge idw_1 \wedge d\overline w_1 \wedge idw_2 \wedge d\overline w_2$ has no mass
on $\Delta$. It follows  that its mass on $\{\|w\| \leq 2, \|z\| \leq r \}$ tends to $0$ as $r \to 0.$
Hence,
$$|\langle T_1\otimes T_2, f dw_1 \wedge d\overline w_1 \wedge dw_2 \wedge d\overline w_2\rangle| \leq  \epsilon_r$$
for a suitable choice of $\epsilon_r$ satisfying the properties in the lemma.

\medskip\noindent
{\bf Case 2.}  Assume that $k=4$ and hence $\gamma=\pm dz_1 \wedge d\overline z_1 \wedge dz_2 \wedge d\overline z_2$.
Let  $m$ be the integer such that $e^{-m-1}<r\leq e^{-m}.$
So $f idz_1 \wedge d\overline z_1 \wedge idz_2 \wedge d\overline z_2$  is a positive form bounded by
 $e^2 r^4 (ir^{-2}(dz_1\wedge d\overline{z}_1+dz_2\wedge d\overline{z}_2))^2$. Since $T_1\otimes T_2$ has no mass
on $\Delta ,$ it follows from
Lemma \ref{L:R_m_annulus}  that
  $$\big| \big\langle T_1\otimes T_2, fdz_1 \wedge d\overline z_1 \wedge dz_2 \wedge d\overline z_2 \big\rangle\big|  \lesssim e^2 r^4\sum_{n,n'=0}^\infty e^{-2n-2n'}  \big\langle T_1\otimes T_2, R_{m+n}\wedge R_{m+n'} \big\rangle.$$ 
The last sum is bounded according to Lemma \ref{L:R_m-R_n}. This proves
the lemma for Case 2.

\medskip\noindent
{\bf Case 3a.}  Assume that $k=2$ and the bi-degree of $\gamma$ in $dz_1,dz_2,d\overline z_1, d\overline z_2$ is $(1,1)$. It follows that 
  the bi-degree of $\gamma$ in $dw_1,dw_2,d\overline w_1, d\overline w_2$ is also $(1,1)$. Observe that $d z_j\wedge d\overline z_k$ is a linear combination of the positive forms
$$idz_j\wedge d\overline z_j, \quad id(z_j\pm z_k) \wedge d\overline {(z_j\pm z_k)} \quad \text{and} \quad 
id(z_j\pm iz_k) \wedge d\overline {(z_j\pm iz_k)}.$$
Moreover, the last forms are bounded by a constant times $idz_1\wedge d\overline z_1+idz_2\wedge d\overline z_2$ because this form is strictly positive. A similar property holds for the variables $w_1$ and $w_2$. Therefore, it is enough to consider the case where 
$\gamma=dz_j\wedge d\overline z_j \wedge dw_k\wedge d\overline w_k$. 

Recall that $(z,w)=(x-y,y)$. So we have 
\begin{eqnarray*}
\big|\big\langle T_1\otimes T_2, f \gamma \big\rangle\big|
&\lesssim & r^2 \int_{\|y\|<1}  \Big(r^{-2}\int_{x\in \B(y,r)}  T_1 (x) \wedge \ddcx {\|x-y\|^2} \Big) T_2(y)\wedge \omega(y) \\
 &\simeq &   r^2 \int_{\|y\|<1}\nu(T_1,y,r)T_2(y)\wedge \omega(y).
\end{eqnarray*}
Applying 
Lemma \ref{l:lelong}   and  Lebesgue's dominated convergence theorem to the expression in the  last line, we see that
it converges to the limit 
$$\int_{\|y\| <1}\nu(T_1,y)T_2(y)\wedge \omega(y)$$ 
when $r$ tends $0.$
By hypothesis, the last integral is  equal to $0$.  This ends the proof of  Case 3a for a suitable choice of $\epsilon_r$.

\medskip\noindent
{\bf Case 3b.}  Assume that $k=2$ and the bi-degree of $\gamma$ in $dz_1,dz_2,d\overline z_1, d\overline z_2$ is $(2,0)$.
It follows that $\gamma=\pm dz_1\wedge dz_2\wedge d\overline w_1\wedge d\overline w_2$.
Let $\chi$ be a smooth function with compact
support in $\{\|w\| < 2, \|z\| < r\}$ such that $0\leq \chi\leq 1$ and $\chi = 1$ in a
neighbourhood of the support of $f$.  By Cauchy-Schwarz inequality, we can bound $|\langle T_1\otimes T_2, f \gamma \rangle |$ from above by
$$\big| \big\langle T_1\otimes T_2, \chi^2dz_1 \wedge d \overline z_1 \wedge dz_2 \wedge d\overline z_2 \big\rangle \big|^{1/2} \big| \big\langle T_1\otimes T_2, f^2dw_1 \wedge d \overline w_1 \wedge dw_2 \wedge d\overline w_2 \big\rangle \big|^{1/2}.$$
According to Cases 1 and 2, the last product is bounded by $\epsilon_r r^2$ for a suitable choice of $\epsilon_r$. This ends the proof of Case 3b.

\medskip\noindent
{\bf Case 3c.}  Assume that $k=2$ and the bi-degree of $\gamma$ in $dz_1,dz_2,d\overline z_1, d\overline z_2$ is $(0,2)$.
This case can be treated in the same way as Case 3b.

\medskip\noindent
{\bf Case 4a.}  Assume that $k=1$ and the bi-degree of $\gamma$ in $dz_1,dz_2,d\overline z_1, d\overline z_2$ is $(1,0)$.
So $\gamma$ has the form $\gamma=\pm dz_j\wedge d\overline w_k \wedge dw_l\wedge d\overline w_l$. 
With $\chi$ as  before, by Cauchy-Schwarz inequality, 
$|\langle T_1\otimes T_2, f \gamma \rangle |$ is bounded from above by
$$ \big| \big\langle T_1\otimes T_2, \chi^2dz_j \wedge d \overline z_j \wedge dw_l \wedge d\overline w_l \big\rangle \big|^{1/2}
\big| \big\langle T_1\otimes T_2, f^2 dw_k \wedge d \overline w_k \wedge dw_l \wedge d\overline w_l \big\rangle \big|^{1/2}.$$
So Case 4a is a consequence of Cases 1 and 3a.

\medskip\noindent
{\bf Case 4b.}  Assume that $k=1$ and the bi-degree of $\gamma$ in $dz_1,dz_2,d\overline z_1, d\overline z_2$ is $(0,1)$.
This case can be treated in the same way as Case 4a.

\medskip\noindent
{\bf Case 5.} Assume that $k=3$. This case can be treated as in Cases 4a and 4b using Cauchy-Schwarz inequality and the previous cases.
\endproof

\subsection{Tangent currents in the local setting}\label{SS:tangent_currents}

We use the notations introduced earlier. In particular,
over $\Delta\cap (5\B\times 5\B)$, with the coordinates $(z,w)$, $\E$ is identified with $\C^2\times 5\B$, $\pi$ is the projection $(z,w)\mapsto w$ and $A_\lambda$ 
is equal to the map $a_\lambda(z,w):=(\lambda z,w)$. 
Tangent currents can be computed locally according to the following result.

\begin{proposition}\label{P:local_compu} 
The mass of $(T_1\otimes T_2)_\lambda$ on any given compact subset of $\E$ is bounded uniformly on $\lambda$ with $|\lambda|\geq 1$.  Moreover, if $(\lambda_n )$ is a
sequence tending to infinity such that $(T_1\otimes T_2)_{\lambda_n}$ converges to a current $\T$,
then in the above local coordinates $(z,w)$, we have
$$\T = \lim\limits_{n\to\infty} (a_{\lambda_n })_* (T_1\otimes T_2) \quad\text{on}\quad \C^2  \times \B.$$
In particular, $\T$ does not depend on the choice of $\tau$.
\end{proposition}

Note that the last assertion in
the proposition is a consequence of the second one because the identity in the proposition doesn't involve the map $\tau$. 
For the proof of this proposition, we need some notions and results.

\begin{definition}\label{D:negligible_forms}\rm  
Let $(\alpha_\lambda)$ be a family of differential $p$-forms on $X\times X$ or  $\E$, depending
on $\lambda \in \C$ with $|\lambda|$ larger than a fixed constant. We say that this family is {\it fine} and we write 
$\alpha_\lambda\in\Fin(\lambda)$
(resp. {\it negligible} and we write $\alpha_\lambda\in\Neg(\lambda)$)
if the support $\supp(\alpha_\lambda )$ of $\alpha_\lambda$ tends to $\Delta$ as $\lambda \to \infty$  and if Properties (1)\,(2) (resp. (1)\,(2)\,(3)) below hold for  all local coordinate systems $(z,w)$ we consider.

\begin{enumerate}
\item $\supp(\alpha_\lambda) \cap (\B\times\B)$ is contained in $(A|\lambda|^{-1}\B)\times\B$ for some
constant $A > 0$ independent of $\lambda;$
\item  The sup-norm of the coefficient of $\gamma$ in $\alpha_\lambda$ is bounded  by $O(\lambda^k)$, 
where $\gamma$ is a wedge-product of $1$-forms among $dz_1, dz_2, dw_1,dw_2$ or their  complex conjugates, and
$k$  is the total degree of $dz_1,dz_2, d\overline{z}_1, d\overline z_2$ in  $ \gamma$, see also Lemma \ref{L:T_1T_2-testforms}.
\item (only for negligible families) The  sup-norm of the coefficient of $\gamma$ is $o(\lambda^k)$ when $\gamma$ is of maximal degree in 
$dz_1,dz_2,d\overline z_1, d\overline z_2$, or equivalently, when $k=p$.
\end{enumerate}
\end{definition}

Note that Properties (1) and (2) are  often easy to check. Properties (2) and (3) are easier to obtain 
if we use the coordinates $(\lambda z , w)$ instead of $(z , w )$. To check that a family is negligible, it is often enough to understand the leading coefficients of the terms of maximal degree in $dz_1,dz_2,d\overline z_1, d\overline z_2$, see also the proof of Lemma \ref{L:key-technique} below.

Negligible families will be used in our study of tangent currents. They enter into the picture in order to handle non-holomorphic changes of variables, i.e. the use of the map $\tau$.
The following lemma will be used in order to establish properties of tangent currents.

\begin{lemma}\label{L:negligible_test_forms}
Let $(\alpha_\lambda)$ be a negligible family of smooth $4$-forms in $X\times X .$ Let $T_1$ and $T_2$ be
as in Lemma  \ref{L:T_1T_2-testforms}. Then we have 
$$\langle T_1\otimes T_2, \alpha_\lambda \rangle \to 0 \quad \text{as} \quad \lambda\to \infty.$$
\end{lemma}
\proof We can use a partition of unity in order to work in local coordinates $(z,w)$ as above.
So we can assume that the forms $\alpha_\lambda $ have supports in $({1\over 2} \B)\times ({1\over 2}\B)$.
Lemma \ref{L:T_1T_2-testforms}, applied to $r := A|\lambda|^{-1}$ with $A $ from Definition \ref{D:negligible_forms}, gives the result.
\endproof

To  study tangent  currents, we need a description of $\tau$ in local coordinates $(z ,w)$ in $\B\times\B$. Consider the Taylor
expansion of order $2$ of $\tau$ in $z,\overline z$ with functions in $w $ as coefficients. Since
$\tau$  is smooth admissible, when $z$ tends to 0, this map and its differential can be written as
\begin{equation}\label{e:local_tau}
\tau (z  , w ) =\big(z+O(\|z\|^2),  w + a(w)z + O(\|z\|^2 )\big),
\end{equation}
and
\begin{equation}\label{e:local_dtau}
d\tau (z , w ) = \big(dz+O^*(\|z\|^2),  dw + O(1)dz + O (\|z\| )\big ) ,
\end{equation}
where $a(w )$ is a $2\times 2$ matrix whose entries are smooth functions in $w$ and $O^*(\|z\|^k)$
is any smooth 1-form
that can be written as
$$O^* (\|z\|^k ) = O(\|z\|^{k-1} )dz + O(\|z\|^{k-1} )d\overline{z} + O(\|z\|^k ).$$
We also have
\begin{equation}\label{e:local_dtau-1}
d\tau^{-1}(z, w ) = \big(dz+O^*(\|z\|^2),  dw + O(1)dz+  O(\|z\| )\big )  .
\end{equation}

\begin{lemma} \label{L:negligible_forms} 
If $(\alpha_\lambda )$ is a fine (resp. negligible) family of $4$-forms on $\E ,$ then $(\tau^* (\alpha_\lambda ))$
is also a fine (resp. negligible) family of $4$-forms on $X\times X .$
\end{lemma}
\proof 
This is a direct consequence of the above local description of $ d\tau$.
\endproof

Recall that $\tau$ is not holomorphic in general but it is close to a holomorphic map near the diagonal $\Delta$. The following lemma suggests that the non-holomorphicity of $\tau$ doesn't affect the computation of tangent currents. 

\begin{lemma}\label{L:key-technique}
Let $\varphi$ be a smooth function with compact support in $\B\times\B$. 
Then the family  $\ddc( \varphi\circ a_\lambda)$ is fine and the families 
 $\ddc (\varphi\circ a_\lambda\circ \tau )- \ddc( \varphi\circ a_\lambda)$, 
 $ \tau^*\big(\ddc (\varphi\circ a_\lambda )\big)- \ddc( \varphi\circ a_\lambda)$ 
 and 
 $\ddc (\varphi\circ a_\lambda\circ \tau )-\tau^*\big(\ddc (\varphi\circ a_\lambda )\big)$
 are negligible, see Definition \ref{D:negligible_forms}. 
\end{lemma}
\proof 
Observe that Property (1) in Definition \ref{D:negligible_forms} is satisfied for all these families of forms. 
In particular, on the supports of the above forms we have $\|z\|\lesssim |\lambda|^{-1}$.  
In order to check Properties (2) and (3) of this definition, we use the following rules of computation
$$\Fin(\lambda)\wedge\Fin(\lambda)=\Fin(\lambda), \quad \Fin(\lambda)\wedge\Neg(\lambda)=\Neg(\lambda) \quad \text{and} \quad \lambda^{-1}\Fin(\lambda)=\Neg(\lambda).$$

When expanding the forms in the lemma using the coordinates $(z,w)$, the definition of $a_\lambda$ and \eqref{e:local_tau}, \eqref{e:local_dtau}, \eqref{e:local_dtau-1}, we only have fine families of forms and for the non-leading terms, an extra factor $O(\lambda^{-1})$ or $O(\|z\|)$ gives us negligible forms. We leave the details to the reader and only highlight some points in the computation.

For simplicity, write $\zeta=(\zeta_1,\zeta_2,\zeta_3,\zeta_4):=(z_1,z_2,w_1,w_2)$ and $s=(s_1,s_2,s_3,s_4):=a_\lambda(\tau(z,w))$. Recall that $\ddc={i\over\pi}\ddbar$ and we have 

\begin{equation*} 
\begin{split}
 \ddbar(\varphi\circ a_\lambda \circ \tau) &=\sum_{m,n=1}^4  {\partial^2\varphi\over \partial \zeta_m\partial\zeta_n}(s)\partial s_m\wedge 
 \overline\partial s_n+ \sum_{m,n=1}^4  {\partial^2\varphi\over \partial \overline\zeta_m\partial\overline\zeta_n}(s)
 \partial \overline s_m\wedge \overline \partial \overline s_n  \\
 &+\sum_{m,n=1}^4  {\partial^2\varphi\over \partial\overline \zeta_m\partial\zeta_n}(s)\partial \overline s_m\wedge
 \overline\partial s_n +\sum_{m,n=1}^4  {\partial^2\varphi\over \partial \zeta_m\partial\overline\zeta_n}(s)\partial s_m\wedge
 \overline\partial \overline s_n\\
 &+ \sum_{m=1}^4 {\partial\varphi\over \partial \zeta_m}(s)\ddbar s_m 
 +\sum_{m=1}^4 {\partial\varphi\over \partial\overline \zeta_m}(s)\ddbar\overline s_m .
 \end{split}
 \end{equation*}
In the same way, we can expand $\ddc ( \varphi\circ a_\lambda)$ and $ \tau^*\big(\ddc (\varphi\circ a_\lambda )\big)$. It is easy to compare them with $\ddc (\varphi\circ a_\lambda \circ \tau)$. For example, using \eqref{e:local_dtau}, we easily see that $\partial s_1- \partial (\lambda z_1)$ is negligible where 
$s_1$ and $\lambda z_1$ are seen as the first coordinate of $a_\lambda(\tau(z,w))$ and $a_\lambda(z,w)$ respectively. 
So the role of $\tau$ is negligible here.

Another point involved in the computation is the comparison between the coefficients of the above forms. For example, using \eqref{e:local_tau}, we can observe that 
\begin{equation*}
 \Big | {\partial^2\varphi\over \partial \zeta_m\partial\overline\zeta_n}(a_\lambda(\tau(z,w)))  - {\partial^2\varphi\over \partial \zeta_m\partial\overline\zeta_n}( a_\lambda(z,w))\Big |\lesssim \|a_\lambda(\tau(z,w))- a_\lambda(z,w)\| \lesssim \|z\|\lesssim |\lambda|^{-1}.
\end{equation*}
Here again, we see that the role of $\tau$ is negligible. The lemma is then obtained by a direct computation.
\endproof

The  following proposition establishes some properties of tangent currents.

\begin{proposition} \label{P:tangent_limits} 
Let $\Phi$  be a continuous $4$-form with support in a fixed compact
subset of $\E.$ Define $\Phi_\lambda := A^*_\lambda (\Phi)$ and $\Psi_\lambda:= \tau^* A^*_\lambda (\Phi)$. Then, we have the following
properties.
\begin{enumerate}
\item If $\Phi \wedge \pi^* (\Omega) = 0$ for any smooth $(2, 2)$-form $\Omega$ on $\Delta ,$ then the families
of $(\Phi_\lambda )$ and $(\Psi_\lambda)$ are negligible.
\item If $\|\Phi\|_\infty \leq 1,$ then $\limsup_{\lambda\to\infty} |\langle T_\lambda , \Phi\rangle|$ is bounded above by a constant which
does not depend on $\Phi.$
\item If $\Phi \wedge\pi^* (\Omega) \geq 0$ for any smooth  positive $(2, 2)$-form $\Omega$ on $\Delta ,$ then any
limit value of $\langle T_\lambda , \Phi\rangle,$ when $\lambda \to \infty,$ is non-negative. In particular, this
property holds when $\Phi$ is a  positive $(2, 2)$-form.
\item If $\Phi = \ddc \phi$ for some smooth $(1,1)$-form $\phi$ with compact support in $\E ,$
then we have $\langle T_\lambda ,\Phi\rangle \to 0$ as $\lambda\to\infty.$ 
\end{enumerate}
\end{proposition}
\proof 
We continue to use the local coordinates $(z,w)$ as above. Observe that if $(\chi_ k )$ is a finite partition of unity for $\Delta$, then $(\chi_k \circ \pi )$ is a finite
partition of unity for $\E$. Using such a partition, we can reduce the problem
to the case where $\Phi$ and $\phi$  have supports in $(r_0\B)\times ({1\over 2}\B)$ for some constant $r_0>0$. 

\smallskip

(1) The hypothesis in (1) implies that the coefficient
of  $dz_1 \wedge d\overline z_1 \wedge  dz_2 \wedge d\overline z_2$ in $\Phi$  vanishes. 
Then, a direct computation shows that $(\Phi_\lambda)$ is
negligible. By Lemma \ref{L:negligible_forms}, the family $(\Psi_\lambda )$ is also negligible.

\smallskip

(2) Modulo a negligible family of forms, thanks to the first assertion, we have
$$\Phi_\lambda\simeq  f_\lambda (z,w)|\lambda|^4 (idz_1 \wedge d\overline{z}_1 \wedge idz_2 \wedge d\overline{z}_2),$$ 
where $f_\lambda$ is a smooth function supported by $(r_0|\lambda|^{-1}\B)\times ({1\over 2}\B)$ and $|f_\lambda |$ is bounded by a constant. Then, we deduce from \eqref{e:local_dtau} that 
 $\Psi_\lambda$ satisfies a similar property and has
support in $(2r_0|\lambda|^{-1}\B)\times ({1\over 2}\B)$ when $\lambda$ is large enough. By Lemma \ref{L:negligible_test_forms},
negligible families of forms do not change the limit we are considering. Thus, Lemma \ref{L:T_1T_2-testforms} implies the result.

\smallskip

(3) We can assume that $\Phi_\lambda$ is as in (2). The hypothesis of (3) implies that the
coefficient of  $idz_1 \wedge d\overline{z}_1 \wedge i dz_2 \wedge d\overline{z}_2$ in $\Phi$ is non-negative. It follows that $f_\lambda \geq 0$. Using \eqref{e:local_tau}, we can see that $\Psi_\lambda$ is the product of a positive function $g_\lambda$ with
 $idz_1 \wedge d\overline{z}_1 \wedge i dz_2 \wedge d\overline{z}_2$ plus a form in a negligible family. 
 Since $T_1$ and $T_2$ are  positive, we have
$$\big\langle T_1\otimes T_2, g_\lambda idz_1 \wedge d\overline{z}_1 \wedge i dz_2 \wedge d\overline{z}_2 \big\rangle \geq 0.$$ 
The result follows easily.

\smallskip

(4) Using local  coordinates, we can write $\phi$ as a finite combination of forms of type $u\ddc v,$  where $u$ and $v$ are  smooth functions  supported by $(r_0\B)\times ({1\over 2}\B)$. For simplicity, we can assume that $\phi=u\ddc v$.
Define 
$$\phi_\lambda := a^*_\lambda(\phi) = (u \circ a_\lambda)\ddc (v\circ a_\lambda) \quad \text{and} \quad 
\psi_\lambda := (u\circ a_\lambda \circ\tau) \ddc  (v\circ a_\lambda\circ\tau).$$
Write $\tau = (\tau_1, \tau_2 )$ in the natural way with $\tau_1,\tau_2$ having values in $\C^2$.   
We have 
$$u\circ a_\lambda = u(\lambda z, w) \quad \text{and} \quad u\circ a_\lambda \circ\tau = u(\lambda \tau_1(z,w), \tau_2(z,w)).$$  
Similar identities hold for $v$ instead of $u$.

Now, observe that $\tau^* (\ddc \phi_\lambda ) -\ddc \psi_\lambda$ is equal to 
\begin{eqnarray*}
\lefteqn{\tau^* \ddc (u\circ a_\lambda)\wedge \tau^* \ddc (v\circ a_\lambda)- \ddc (u\circ a_\lambda \circ\tau)\wedge 
\ddc (v\circ a_\lambda \circ\tau) }\\
&=& \big[\tau^*\ddc (u\circ a_\lambda)-\ddc(u\circ a_\lambda\circ \tau)\big] \wedge \big[\tau^*\ddc(v\circ a_\lambda)\big]\\
&& + \big[\ddc(u\circ a_\lambda\circ \tau)\big] \wedge\big[\tau^*\ddc (v\circ a_\lambda)-\ddc(v\circ a_\lambda\circ \tau)\big].
\end{eqnarray*}
Using Lemma \ref{L:key-technique}, Definition \ref{D:negligible_forms} and the rules of computations given in the proof of Lemma \ref{L:key-technique}, we
can check that both terms in the last sum belong to negligible families of $4$-forms. 
 
It follows from  Lemma \ref{L:negligible_test_forms} that
$$\big\langle (T_1\otimes T_2)_\lambda , \ddc \phi \big\rangle = \big\langle T_1\otimes T_2, \tau^* (\ddc \phi_\lambda) \big\rangle = \big\langle T_1\otimes T_2, \ddc \psi_\lambda \big\rangle+ o(1) \quad \text{as} \quad\lambda \to\infty.$$
It remains to show that $ \big\langle T_1\otimes T_2, \ddc \psi_\lambda \big\rangle$ tends to 0. 
Using Lemma \ref{l:identity}, we have 
$$ \big\langle T_1\otimes T_2, \ddc \psi_\lambda \big\rangle = - \langle\dbar S_1\otimes \partial \overline{S}_2, \ddc\psi_\lambda\rangle -
 \langle\partial \overline{S}_1\otimes \dbar S_2, \ddc \psi_\lambda \rangle.$$
By Lemmas \ref{L:negligible_forms} and \ref{L:key-technique}, the family $(\ddc\psi_\lambda)$ is fine. Therefore, by Lemma \ref{l:kernel-Delta}, it is enough to show that 
$\ddc\psi_\lambda$ tends to 0 weakly. 

Since the family $(\ddc\psi_\lambda)$ is fine, the mass of $\ddc\psi_\lambda$ is bounded. So, when $\lambda$ tends to infinity, this sequence accumulates to $4$-currents of finite mass supported by $\Delta$. Moreover, since $\ddc\psi_\lambda$ is $d$-exact, any limit $R$ of 
$\ddc\psi_\lambda$ is a $d$-exact 4-current. In particular, $R$ is a normal 4-current supported by $\Delta$. Thus, 
we can identify it to a 0-current on $\Delta$, according to the classical support theorem, see \cite{Federer}. Finally, since the only $d$-exact 0-current on $\Delta$ is zero, we get $R=0$. 
The result follows.
\endproof

We continue the proof of Proposition \ref{P:local_compu}.
The second assertion
in Proposition \ref{P:tangent_limits} implies that
the mass of $(T_1\otimes T_2)_\lambda$ on any given compact subset of $\E$ is bounded uniformly on $\lambda$ with $\lambda$ large enough.  

Consider any sequence $(\lambda_n )$ of complex numbers tending to infinity. After extracting a subsequence, 
we can assume that $(T_1\otimes T_2)_{\lambda_n}$
converges to a $4$-current $\T$ of locally finite mass in $\E$. The first assertion
in Proposition \ref{P:tangent_limits} shows that in the above local coordinates $(z,w)$, if the coefficient
of $dz_1 \wedge d\overline{z}_1 \wedge dz_2 \wedge d\overline{z}_2$  in $\Phi$ vanishes then $\langle \T,\Phi\rangle =0.$
Consequently,  we have $\T\wedge  dw_j =0$ and   $\T\wedge  d\overline w_j =0.$
Hence, $\T$ is a current of bi-degree $(2, 2)$. 

The third assertion of Proposition \ref{P:tangent_limits} implies that $\T$ is positive. Finally, the fourth assertion in that proposition is
equivalent to saying that $\T$ is $\ddc$-closed.

\begin{lemma} \label{l:nu}
There is a positive measure $\vartheta$ on $\Delta$ such that $\T = \pi^* (\vartheta).$ In
particular, the current $\T$ is closed.
\end{lemma}
\proof
We follow  the argument in the proof of  \cite[Lem.\,3.7]{DinhSibony18}.
Consider the family $\Gc$ of all positive $\ddc$-closed $(2, 2)$-currents $R$ on
$\E$  which are {\it vertical} in the sense that $R \wedge \pi^* (\Omega) = 0$ for any smooth form
$\Omega$ of positive degree  on $\Delta.$

\smallskip

\noindent {\bf Claim.}  If $R$ is any current in $\Gc$ and $v$ is a smooth positive function on $\Delta,$ then
$(v \circ \pi )R$ also belongs to $\Gc .$

\smallskip

Indeed, it is clear that $(v\circ\pi)R$ is a positive and vertical $(2,2)$-current. The only point to check is that $(v \circ \pi)R$ is $\ddc$-closed.
Define $\widetilde v := v \circ  \pi$. We have $\ddc R = 0$ and since $R$ is vertical, we get that $d \widetilde v \wedge R = 0$,
 $\dc \widetilde v \wedge R = 0$ and $\ddc \widetilde v \wedge R = 0.$
Consequently, a straightforward calculation shows that
$$\ddc (\widetilde v R) = d(\dc \widetilde v \wedge R) - \dc(d\widetilde v\wedge R) -\ddc\widetilde v\wedge R'+  \widetilde v \ddc R = 0,$$
which completes the proof of the claim.

We infer  from  the claim that every extremal element in $\Gc$ is supported
by a fiber of $\pi$ which is a complex plane. A positive $\ddc$-closed $(2,2)$-current on a complex plane is defined by
a positive  pluriharmonic  function.
On the other hand,  positive  plurisubharmonic functions on a complex plane are necessarily constant.  
Hence, extremal elements in $\Gc$ are proportional to the currents of integration on fibers of $\pi$.
In order to get the lemma, we only need to show that any $R$ in $\Gc$ is an average of those
extremal currents.

Consider the convex cone of positive $\ddc$-closed vertical currents $R$
as above. Observe that the set of currents
with mass 1  is compact and is a basis of the considered cone.
Therefore, Choquet's representation theorem implies that any current in the
cone is an average on the extremal elements. The lemma follows.
\endproof

\proof[End of the proof of Proposition \ref{P:local_compu}] Consider a smooth test $4$-form $\Omega$ with compact
support in $\C^2\times\B.$ Denote by $f (z,w)$ the coefficient of $dz_1 \wedge d\overline z_1\wedge  dz_2 \wedge d\overline z_2$ in $\Omega.$ By the definition of $\T$, Proposition \ref{P:tangent_limits} and the above discussion on negligible families of forms, we see that  only the component $ f (z,w)dz_1 \wedge d\overline z_1\wedge  dz_2 \wedge d\overline z_2$
of $\Omega$ matters in computing the limit. So we have
\begin{eqnarray*}
\langle \T, \Omega\rangle &=& \lim_{n\to\infty} \big\langle T_1\otimes T_2, \tau^* A^*_{\lambda_n} ( f (z,w)dz_1 \wedge d\overline z_1\wedge  dz_2 \wedge d\overline z_2) \big\rangle\\
&=& \lim_{n\to\infty} \big\langle T_1\otimes T_2, |\lambda_n |^4 f (\lambda_n\tau_1 , \tau_2 )\tau_1^*(dz_1 \wedge d\overline z_1\wedge  dz_2 \wedge d\overline z_2) \big\rangle
\end{eqnarray*}
and
\begin{eqnarray*}
\langle T_1\otimes T_2, a_{\lambda_n}^*\Omega\rangle &=& \lim_{n\to\infty} \big\langle T_1\otimes T_2,  a^*_{\lambda_n} ( f (z,w)dz_1 \wedge d\overline z_1\wedge  dz_2 \wedge d\overline z_2) \big\rangle\\
&=& \lim_{n\to\infty} \big\langle T_1\otimes T_2, |\lambda_n |^4 f (\lambda_n z , w )dz_1 \wedge d\overline z_1\wedge  dz_2 \wedge d\overline z_2 \big\rangle.
\end{eqnarray*}
Thus, it is enough to check that the family of forms
$$ |\lambda_n|^4 f (\lambda_n\tau_1 , \tau_2 )\tau_1^*(dz_1 \wedge d\overline z_1\wedge  dz_2 \wedge d\overline z_2)-|\lambda_n|^4 f (\lambda_n z , w )dz_1 \wedge d\overline z_1\wedge  dz_2 \wedge d\overline z_2$$
is negligible. But this can be easily obtained using the same computation as in the proof of Lemma \ref{L:key-technique}.
\endproof

\proof[Proof of Theorem \ref{T:main_2}\,(1)]
This is a direct consequence of Proposition \ref{P:local_compu} and Lemma \ref{l:nu} above. 
\endproof

\subsection{Proof of the mass formula}\label{S:mass_formula}

In this section, we prove Part (2) of Theorem \ref{T:main_2}.  For this purpose, we will use families of smooth test closed $4$-forms $\widehat\Phi_\lambda$ and $\Phi_\lambda$ on $\E$ and on $X\times X$ that we construct below.

\begin{lemma} \label{l:Phi-hat}
There is a smooth closed $(2,2)$-form $\widehat\Phi$ with compact support in $\E$ which is cohomologous to $\Delta$. In particular, we have 
$$\int_{\pi^{-1}(x,x)} \widehat\Phi=1 \quad \text{for every point} \quad (x,x)\in \Delta.$$
\end{lemma}
\proof
Observe that we can compactify $\E$ in order to get a compact K\"ahler manifold $\overline\E$. According to \cite{DinhSibony04}, we can regularize the current $[\Delta]$ on $\overline\E$. More precisely, there is a sequence of smooth closed $(2,2)$-forms $T_n$ on $\overline\E$ converging to $[\Delta]$ in the sense of currents. Each form $T_n$ can be written as the difference of two positive closed $(2,2)$-forms. Moreover, the support of $T_n$ tends to $\Delta$ as $n$ tends to infinity, see \cite[Rk 4.5]{DinhSibony04}. So, for the first assertion in the lemma, it is enough to choose $\Phi=T_n$ with $n$ large enough.

For the second assertion, observe that 
the measure $\Phi\wedge [\pi^{-1}(x,x)]$ is cohomologous to the Dirac mass  $[\Delta]\wedge  [\pi^{-1}(x,x)]$.
Hence, the integral of $\Phi\wedge [\pi^{-1}(x,x)]$ is equal to 1. This ends the proof of the lemma.
\endproof

Define for $|\lambda|\geq 1$
$$\widehat\Phi_\lambda:=(A_\lambda)^*(\widehat\Phi) \quad \text{and} \quad \Phi_\lambda:=\tau^*(\widehat\Phi_\lambda).$$
Clearly, $\widehat\Phi_\lambda$ is a smooth closed $(2,2)$-form and $\Phi_\lambda$ is a smooth closed 4-form. 

\begin{lemma} \label{l:Phi-limit}
The form $\Phi_\lambda$ converges to $[\Delta]$
in the sense of currents when $\lambda$ goes to infinity.  Moreover, 
the three families $(\Phi_\lambda)$, $(\lambda\partial\Phi_\lambda)$ and $(\lambda\dbar\Phi_\lambda)$ are fine. 
\end{lemma}
\proof
It is not difficult to see that $\widehat\Phi_\lambda$ converges to $[\Delta]$ in the sense of currents. The first assertion in the lemma follows easily. 
The second assertion is obtained as in Subsection \ref{SS:tangent_currents} by using that $\partial\widehat\Phi_\lambda=0$, $\dbar\widehat\Phi_\lambda=0$, and also $\dbar \tau=O(|\lambda|^{-1})$, $\partial \overline\tau=O(|\lambda|^{-1})$ on the support of $\Phi_\lambda$.
\endproof

\proof[End of the proof of Theorem \ref{T:main_2}(2).]
Recall that in \eqref{e:decompo_posi_har}, the form $\dbar S_j$ and $\partial\overline S_j$ are uniquely determined by $T_j$. Therefore, we will use here the following decomposition given by Proposition \ref{p:T-rep}
\begin{equation} \label{e:Tj-best}
T_j=\Omega_j+\partial S_j+\dbar \overline S_j + i\ddbar u_j,
\end{equation}
where $\Omega_j$ is a smooth closed $(1,1)$-form, $S_j, \overline S_j, \partial S_j, \partial\overline S_j, \dbar S_j, \dbar\overline S_j$ are forms of class $L^2$ and $u_j, \partial u_j, 
\dbar u_j$ are forms of class $L^p$ for every $1\leq p<2.$  

By Lemma \ref{l:Phi-hat}, we have 
$$\|\vartheta\|=\big\langle \pi^*(\vartheta),\widehat\Phi \big\rangle =\lim_{n\to\infty} \big\langle (T_1\otimes T_2)_{\lambda_n} , \widehat\Phi\big\rangle = \lim_{n\to\infty} \big\langle T_1\otimes T_2,\Phi_{\lambda_n} \big\rangle.$$
We use now \eqref{e:Tj-best}, Stokes' theorem and the fact that $\Phi_{\lambda_n}$ is closed (but not necessarily $\partial$-closed or $\dbar$-closed) in order to expand the last integral $\big\langle T_1\otimes T_2,\Phi_{\lambda_n} \big\rangle$ as in Lemma \ref{l:identity}. Since 
$\Phi_{\lambda_n}$ is closed, the terms like $\big\langle \Omega_1 \otimes i\ddbar u_2,\Phi_{\lambda_n} \big\rangle$ or 
$\big\langle i\ddbar u_1 \otimes i\ddbar u_2,\Phi_{\lambda_n} \big\rangle$ vanish. We have
$$\big\langle T_1\otimes T_2,\Phi_{\lambda_n} \big\rangle = \langle \Omega_1\otimes \Omega_2,\Phi_{\lambda_n}\rangle - \langle  \dbar S_1\otimes \partial\overline S_2, \Phi_{\lambda_n} \rangle 
- \langle  \partial \overline S_1\otimes \dbar  S_2, \Phi_{\lambda_n} \rangle$$
$$  - \langle  \dbar S_1\otimes \overline S_2, \partial \Phi_{\lambda_n} \rangle  -i\langle \dbar u_1\otimes \partial S_2, \partial \Phi_{\lambda_n}\rangle + \text{similar terms involving 
$\partial u_j,$ $\dbar u_j,$ $\partial\Phi_{\lambda_n}$ or $\dbar\Phi_{\lambda_n}$.}$$
We can now apply Lemma \ref{l:Phi-limit}, and then Lemma \ref{l:kernel-Delta} with $c=1$ for  the first three terms in the last sum.
Their sum converges to 
$$\int_X\Omega_1\wedge\Omega_2 - \int_X \dbar S_1 \wedge \partial \overline S_2 - \int_X \partial\overline S_1\wedge \dbar S_2.$$ 
Then, Lemma \ref{l:kernel-Delta-bis} shows that the other terms in the above expression of $\big\langle T_1\otimes T_2,\Phi_{\lambda_n} \big\rangle$ tend to 0. 
Recall that we use \eqref{e:Tj-best} given by Proposition \ref{p:T-rep}.
This completes the proof of the theorem.
\endproof


\section{Vanishing of the tangent currents in the foliation setting} \label{S:Vanishing}

  The purpose of this  section is  to prove  Theorem  \ref{T:main_3}.
  The proof is  given in the  first subsection modulo two auxiliary  propositions  which will   be  proved in
  the last two subsections.

\subsection{Main steps of the proof of the vanishing  theorem} \label{SS:Vanishing_thm}
  
Let $\Fc$ be as in  Theorem  \ref{T:main_3}.  Consider 
a positive $\ddc$-closed $(1, 1)$-current $T$ directed by $\Fc$. 
Recall that if $T$ has positive mass on a leaf,  then this leaf is an invariant closed analytic curve of $\Fc$, see Theorem \ref{T:auxiliary}.
So for Theorem  \ref{T:main_3}, we can assume that  $T$ has no
mass on  each  single leaf of $\Fc$. It follows from  \eqref{e:T-box} and \eqref{e:Lelong} that 
$\nu(T,x)=0$ for all $x$ outside the singularities of $\Fc$.
Since positive $\ddc$-closed $(1,1)$-currents have no mass on finite sets, we can apply 
Theorem \ref{T:main_2} to the tensor product $T\otimes T$. 

Consider a tangent current $\T$ to  $T\otimes T$ along $\Delta$. With the notation as in the above sections, there is a sequence $\lambda_n$ converging to infinity and a positive measure $\vartheta$ on $\Delta\simeq X$ such that 
$$\T=\lim_{n\to\infty} (T\otimes T)_{\lambda_n} = \pi^*(\vartheta).$$ 
We can identify $\vartheta$ with a positive measure on $X$. Recall that by Theorem \ref{T:main_2} the mass $m$ of $\vartheta$ does not depend on the choice of $\T$. The following propositions will be proved in the next subsections.

\begin{proposition}\label{P:Kaufmann} 
For every choice of the tangent current $\T$, the measure  $\vartheta$  is  supported on the singularities of $\Fc$. 
\end{proposition}

Throughout this section, we consider $\lambda$ real such that $\lambda>1$ and $s:=\log\lambda>0$. We refer to \eqref{e:Expec} for the notion of expectation $\Ebf(\cdot)$. 
 
\begin{proposition}\label{P:singularities}
We have 
$$\lim_{s\to\infty} \Ebf((T\otimes T)_\lambda) =0$$
in a neighbourhood of each point $(p,p)\in\Delta$, where $p$ is any singular point of $\Fc$.
\end{proposition}

\proof[End of the proof of Theorem \ref{T:main_3}]
Let $\T'$ be a limit current of $\Ebf((T\otimes T)_\lambda)$ when $s=\log \lambda$ tends to infinity. 
This current belongs to the convex hull of all the above tangent currents $\T$. So we have $\T'=\pi^*(\vartheta')$ for some positive measure $\vartheta'$ of mass $m$ on $\Delta\simeq X$.
By  Propositions \ref{P:Kaufmann} and \ref{P:singularities},  we have $\vartheta'=0$. Therefore, we get $m=0$ and hence, by the mass formula in Theorem \ref{T:main_2}, we have $\T=0$ for any choice of $\T$. This proves the theorem.
\endproof

\subsection{Vanishing of the tangent  currents outside the singularities} \label{SS:outside}

We follow the same lines as in  Kaufmann's work \cite{Kaufmann}. Consider any flow box $\U$ of $\Fc$ outside the singularities, see Subsection \ref{SS:positive_ddc_closed_currents}. So we can choose  holomorphic coordinates
$x=(x_1,x_2)$ in which the plaques of $\Fc$ in
$\U$ are given by
$$L_\alpha = \{x_2  =\phi_\alpha (x_1)\},$$
where $\phi_\alpha :\ 3\D  \to 3\D$ is a holomorphic function  such that $\phi_\alpha (0)=\alpha$ and 
\begin{equation} \label{e:Lip-gr}
\kappa_0^{-1} |\alpha-\beta|\leq |\phi_\alpha(x_1)-\phi_\beta(x_1)|\leq \kappa_0|\alpha-\beta|
\end{equation}
for all $x_1,\alpha,\beta$ in $3\D$ and for some constant $\kappa_0\geq 1$.

Since $T$ is a diffuse positive $\ddc$-closed current  directed by $\Fc$, as in \eqref{e:T-box},  we have the following  decomposition  in the flow
box $\U,$ 
$$T = \int h_\alpha[L_\alpha ] d\mu(\alpha),$$
where $[L_\alpha ] $ denotes the current of integration along the plaque $L_\alpha,$ $h_\alpha$ is  a positive harmonic  function on $L_\alpha$ for $\mu$-almost every $\alpha\in 3\D,$
and $\mu$ is a diffuse positive
measure of finite mass on $3\D.$  We  multiply $\mu$ by  the  positive  function $ h_\alpha(0,\alpha)$  and  divide $h_\alpha$  by $h_\alpha(0,\alpha)$  in order to assume that $h_\alpha(0,\alpha)=1$ for  $\mu$-almost every $\alpha\in 3\D.$
By Harnack's inequality, there is a constant $\kappa \geq 1$  such  that  (we reduce slightly the flow box if necessary)
\begin{equation}\label{e:kappa}
\kappa^{-1}\leq h_\alpha(x)\leq \kappa\quad\mbox{for  $\mu$-almost every $\alpha\in 3\D$ and for $x\in L_\alpha.$}
 \end{equation}

Consider the product foliation
$\Fc \times \Fc$ on $X \times X.$ The above coordinates on the
flow box $\U$ induce natural holomorphic coordinates $(x, y)=(x_1,x_2,y_1,y_2)$ on $\U \times \U$ in which the
plaques of $\Fc \times \Fc$ are given by
$$L_{\alpha,\beta} := L_\alpha \times L_\beta = \big\{x_2 = \phi_\alpha (x_1 ),\ y_2 = \phi_\beta (y_1 ) \big\}.$$
The tensor product $ T \otimes T$ is a positive  current  of bi-dimension $(2,2)$ on $X \times X$ directed
by $\Fc \times \Fc$ which is given on $\U \times \U$ by
$$T\otimes T =\int (h_\alpha \otimes h_\beta) [L_{\alpha,\beta} ] d(\mu\otimes \mu)(\alpha, \beta).$$
Since    $\mu$ has no atoms, by
Fubini's theorem, $ \mu \otimes  \mu$ gives no mass to the set $\{\alpha=\beta\}$ in $3\D \times 3\D$, or equivalently, $T\otimes T$ gives no mass to the diagonal $\Delta$ of $X\times X$. 

To investigate  the tangent currents of $T\times T$ along  $\Delta = \{x = y\},$ 
it is convenient to work in the holomorphic coordinates  $(z,w) := (x-y, y)$ and to use new parameters $\zeta= (\zeta_1 , \zeta_2 )$  with $\zeta_1 := \alpha-\beta$ and $\zeta_2 := \beta$.
Write $z = (z_1 , z_2)$ and  $w=(w_1,w_2).$
In the coordinate system $(z,w)$, the diagonal $\Delta$ is given by the equation $z = 0$. Since $x=z+w$, $y=w$, $\alpha=\zeta_1+\zeta_2$ and $\beta=\zeta_2$, the plaque $L_{\alpha,\beta}$
transforms to (here, we are only interested in parts of $L_{\alpha,\beta}$ near the origin) 
$$\Gamma_\zeta := \big\{(z_1 , f_\zeta (z_1 , w_1 ), w_1 , \phi_{\zeta_2}(w_1 )) \quad \text{with} \quad z_1,w_1\in \D \big\},$$
where 
$$f_\zeta (z_1 , w_1 ) := \phi_{\zeta_1+\zeta_2} (z_1+w_1)-\phi_{\zeta_2}(w_1).$$ 

Always in the coordinates $(z,w)$, the decomposition of $T\otimes T$ becomes
$$T\otimes T =\int h_{\zeta_1+\zeta_2}(z+w) h_{\zeta_2}(w) [\Gamma_\zeta ] d(\mu\otimes \mu) (\zeta).$$
The dilation $A_\lambda$ in the direction normal to $\Delta$ is equal to the map $a_\lambda(z,w):=(\lambda z,w)$.
Note that \eqref{e:Lip-gr} implies that 
the distance between $\Gamma_\zeta$ and $\Delta$ is bounded below by a positive constant times $|\zeta_1|$. 
Such properties allow us to obtain as in  \cite[Lem.\,4.4, 4.5]{Kaufmann} the following lemma.

\begin{lemma}\label{L:Kaufmann}
\begin{enumerate}
 \item The mass of $(a_\lambda )_* [\Gamma_\zeta]$ on any given compact set is
bounded uniformly in $(\lambda, \zeta)$ with $|\zeta_1| \leq |\lambda|^{-1}.$
\item There exists a ball $W$ centered at the origin  such that
$(a_\lambda )_* [\Gamma_\zeta]$ has no mass on $W$ for every pair $(\lambda, \zeta)$ such that $|\zeta_1|> |\lambda|^{-1}.$
\end{enumerate}
\end{lemma}

 \proof[Proof of  Proposition \ref{P:Kaufmann}]
 We only need to show that any  limit of the current
 $(a_\lambda)_*(T\otimes T)$ is zero in $W$ when $\lambda$ tends to infinity.
 Using the estimate  \eqref{e:kappa}, we see that
 $$(a_\lambda)_*(T\otimes T)\leq \kappa^2 S_\lambda,
\quad \text{where} \quad  S_\lambda:=\int (a_\lambda)_*[\Gamma_\zeta]d(\mu\otimes \mu)(\zeta).$$

Write $S_\lambda= S'_\lambda+S''_\lambda$ with
$$S'_\lambda := \int_{|\zeta_1|\leq |\lambda|^{-1}} (a_\lambda)_*[\Gamma_\zeta]d(\mu\otimes \mu)(\zeta) \quad\text{and}\quad
S''_\lambda : = \int_{|\zeta_1|> |\lambda|^{-1}} (S_\lambda)_*[\Gamma_\zeta]d(\mu\otimes \mu)(\zeta).$$
By Lemma \ref{L:Kaufmann}(2), we have $S''_\lambda=0$ on $W$. By Lemma \ref{L:Kaufmann}(1), the mass of $S'_\lambda$ over $W$ is bounded by a constant times
$(\mu\otimes\mu)(\{|\zeta_1| < |\lambda|^{-1} \}).$ The last quantity tends to 0 as $\lambda$ tends to infinity because $\mu\otimes\mu$  gives no mass to the set $\{\zeta_1 = 0\}=\{\alpha=\beta\}$. Therefore, $S'_\lambda$ tends to 0 on $W$ when $\lambda$ tends to infinity. This ends the proof of the proposition.
 \endproof

\subsection{Vanishing of the tangent  currents near the singularities} \label{SS:near_singularities}

In this subsection, we will give the proof of Proposition \ref{P:singularities}. From now on, we only consider real positive parameters $\lambda=e^s$ with $s>0$ and 
place ourselves in the setting of Appendix \ref{a:Harnack}. In particular, the properties of some segments and half-lines in the sector $\S$, described after Lemma \ref{l:Poisson},
are important in our study. We continue to use the notations introduced at the end of the Introduction.

As in \eqref{e:x-zeta}, we will use the following parametrization of $\L_\alpha$
\begin{equation} \label{e:x-zeta-bis} 
x_1 = \alpha e^{i\eta (\zeta +\log {|\alpha|/b})}\quad \text{and}\quad x_2 = e^{i(\zeta +\log{ |\alpha|/b})}\quad \text{with}\quad \zeta = u+iv \in \C
\end{equation}
and similarly, we will use the following parametrization of $\L_\beta$
\begin{equation} \label{e:y-zeta} 
y_1 = \beta e^{i\eta (\check\zeta +\log {|\beta|/b})}\quad \text{and}\quad y_2 = e^{i(\check\zeta +\log{ |\beta|/b})}\quad \text{with}\quad \check\zeta = \check u + i \check v \in \C.
\end{equation}
For $\theta=(\theta_1,\theta_2)\in \D^2$, $\lambda=e^s>1$, and $\alpha,\beta\in \A$, consider the intersection 
$$Z_{\alpha,\beta,\theta}^\lambda := \big\{ (\zeta,\check\zeta)\in \S\times\S, \quad (x,y)\in (L_\alpha\times L_\beta)\cap \{z=\theta/\lambda\}\big\},$$
where the points are counted with multiplicity. 

In $\D^2\times\D^2$, the intersection of the current $T_\alpha \otimes T_\beta$ with the current of integration on the $2$-dimensional complex plane $\{z=\theta/\lambda\}$  is equal to the positive measure
\begin{equation} \label{e:vartheta}
\vartheta_{\alpha,\beta,\theta}^\lambda:=(T_\alpha \otimes T_\beta) \wedge [z=\theta/\lambda] = \sum_{(\zeta,\check\zeta)\in Z_{\alpha,\beta,\theta}^\lambda} H_\alpha(x_1,x_2)H_\beta(y_1,y_2)\delta_{(x,y)},
\end{equation}
where $\delta_{(x,y)}$ is the Dirac mass at the point $(x,y)$. 
Consider also the open set 
\begin{equation} \label{e:Theta}
\Theta:=\big\{\theta=(\theta_1,\theta_2)\in \D^2, \quad  |\theta_1-1|<\epsilon_0, |\theta_2-1| < \epsilon_0  \big\},
\end{equation}
where $\epsilon_0>0$ is a fixed small enough constant depending only on $\eta$. We will show that the masses of the measures $\vartheta_{\alpha,\beta,\theta}^\lambda$ satisfy the following property. 

\begin{proposition} \label{p:theta-slice}
The following property holds for all singularities of $\Fc$. There is a constant $c>0$ such that for $\mu$-almost every $\alpha,\beta\in\A$ and all $s=\log\lambda >0$ we have
$$\sup_{\theta\in\Theta} \int_{\alpha,\beta\in\A} \Ebf\big(\|\vartheta_{\alpha,\beta,\theta}^\lambda\|\big) d\mu(\alpha)d\mu(\beta) \leq c \quad \text{and} \quad \lim_{s\to\infty} \sup_{\theta\in\Theta} \int_{\alpha,\beta\in\A} \Ebf\big(\|\vartheta_{\alpha,\beta,\theta}^\lambda\|\big) d\mu(\alpha)d\mu(\beta) =0.$$
\end{proposition}

Taking into account this result, we first complete the proof of Proposition \ref{P:singularities}.

\proof[Proof of Proposition \ref{P:singularities}]
Consider a cluster value $\T'$ of the family $\Ebf((T\otimes T)_\lambda)$ when $s=\log \lambda$ tends to infinity. We need to show that $\T'=0$. By  Proposition \ref{P:Kaufmann}, we only need to check this property near a singularity of the foliation. Moreover, 
by Theorem \ref{T:main_2}, in the local setting we consider, $\T'$ is a cluster value of $\Ebf((a_\lambda)_*(T\otimes T))$ when $s=\log\lambda$ tends to infinity. We also have $\T'=c [\pi^{-1}(0)]$ for some constant $c\geq 0$. Our goal is to show that $c=0$. 

For this purpose, on the open set 
$$\Theta\times\D^2=\big\{(z,w)\in \D^2\times \D^2, \quad z\in\Theta \big\},$$
we consider the following measures  
$$\vartheta_\lambda:= (a_\lambda)_*(T\otimes T) \wedge (idz_1\wedge d\overline z_1)\wedge (idz_2\wedge d\overline z_2) 
\quad \text{and} \quad \vartheta_\lambda':=\Ebf(\vartheta_\lambda).$$
It is enough to show that the mass of $\vartheta_\lambda'$ 
tends to 0 as $s=\log\lambda$ tends to infinity. Indeed, since $\T'=c [\pi^{-1}(0)]$, this property implies that $\T'$ vanishes on $\Theta\times\D^2$ and hence $c=0$. 

Observe that the mass of $\vartheta_\lambda$ is equal to the mass of $(a_{\lambda})^*(\vartheta_\lambda)$. The last mass is equal to a constant times the average of 
$$m(\theta):=\int_{\alpha,\beta\in\A} \|\vartheta_{\alpha,\beta,\theta}^\lambda\| d\mu(\alpha)d\mu(\beta).$$ 
with respect to the Lebesgue measure on $\theta\in \Theta$. 
The involved constant is the integral of $(idz_1\wedge d\overline z_1)\wedge (idz_2\wedge d\overline z_2)$ on $\Theta$.
We deduce that the mass of  $\vartheta_\lambda'$ is equal to a constant times the average of 
$$m'(\theta):=\int_{\alpha,\beta\in\A} \Ebf(\|\vartheta_{\alpha,\beta,\theta}^\lambda\|) d\mu(\alpha)d\mu(\beta).$$
The estimate in Proposition \ref{p:theta-slice} and Lebesgue's dominated convergence theorem imply that $m'(\theta)$ tends to 0 uniformly on $\theta$. Thus, the mass of $\vartheta_\lambda'$ tends to 0. The result follows.
\endproof

The rest of this section is devoted to prove Proposition \ref{p:theta-slice}. We need to understand the set $Z_{\alpha,\beta,\theta}^\lambda$ which 
is the set of all solutions of the following system of equations with unknown $(\zeta,\check\zeta)$ in $\S\times\S$, see also \eqref{e:x-zeta-bis} and \eqref{e:y-zeta})
\begin{equation} \label{e:system-Z}
\left\lbrace 
\begin{array}{l}
 x_1-y_1=\theta_1/\lambda\\
 x_2-y_2=\theta_2/\lambda
\end{array} \right.
\ \ 
\Longleftrightarrow
\ \ 
\left\lbrace 
\begin{array}{l}
 y_1(\rho_1-1)=\theta_1/\lambda\\
 y_2(\rho_2-1)=\theta_2/\lambda
\end{array} \right.
\ \ 
\Longleftrightarrow
\ \ 
\left\lbrace 
\begin{array}{l}
 x_1(1-1/\rho_1)=\theta_1/\lambda\\
 x_2(1-1/\rho_2)=\theta_2/\lambda \, ,
\end{array} \right.
\end{equation} 
where the ratios $\rho_1$ and $\rho_2$ are defined by 
\begin{equation} \label{e:ratio}
\rho_1:={x_1 \over y_1}= {\alpha \over\beta} e^{i\eta(\zeta-\check\zeta+\log|\alpha|/b-\log|\beta|/b)} \quad \text{and} \quad  \rho_2:={x_2 \over y_2}=e^{i(\zeta-\check\zeta+\log|\alpha|/b-\log|\beta|/b)}. 
\end{equation}
Observe that these ratios are not equal to 0 and nor to 1 because $\theta\in\Theta$. 

\begin{lemma} \label{l:Z-separated}
There is a constant $N>0$ such that  $Z_{\alpha,\beta,\theta}^\lambda$ is $N$-sparse for all 
$\alpha,\beta\in\A$, $\theta\in\Theta$ and $\lambda=e^s>1$, see also Definition \ref{d:dominated}.
\end{lemma} 
\proof
Let $(\zeta^*,\check\zeta^*)\in \S\times\S$ be any point in $Z_{\alpha,\beta,\theta}^\lambda$. Denote by $x^*,y^*,\rho_1^*,\rho_2^*$ the corresponding values of $x,y,\rho_1,\rho_2$.  We only consider $|\rho_2^*|\leq 1$ because the opposite case can be treated in the same way.
Write $\zeta=\zeta^*+\xi$ and $\check\zeta=\check\zeta^*+\check\xi$. Using \eqref{e:x-zeta-bis} and \eqref{e:y-zeta}, we see that
the system \eqref{e:system-Z} is equivalent to
\begin{equation} \label{e:system-Z-local}
\left\lbrace 
\begin{array}{l}
 x_1^*(e^{i\eta \xi}-1)-y_1^*(e^{i\eta\check\xi}-1)=0\\
 x_2^*(e^{i\xi}-1)-y_2^*(e^{i\check\xi}-1)=0 
\end{array} \right.
\quad \Longleftrightarrow \quad
\left\lbrace 
\begin{array}{l}
 \rho_1^*(e^{i\eta \xi}-1)-(e^{i\eta\check\xi}-1)=0\\
 \rho_2^*(e^{i\xi}-1)-(e^{i\check\xi}-1)=0 \, .
\end{array} \right.
\end{equation} 
We are interested in the case where both $\xi$ and $\check\xi$ are small. 

Since $|\rho_2^*|\leq 1$, from the second equation of  
\eqref{e:system-Z-local} we get 
$$\check\xi=f_{\rho_2^*}(\xi) \quad \text{with} \quad f_{\rho_2^*}(\xi):= -i \log \big[1+ \rho_2^* (e^{i\xi}-1)\big],$$
where we use the principal branch of the function $\log$. Substituting this value of $\check\xi$ to the first equation of \eqref{e:system-Z-local} gives
$$\rho_1^*(e^{i\eta \xi}-1) - (e^{i\eta f_{\rho_2^*}(\xi)}-1) =0.$$
The solutions of this equation are the zeros of the function
\begin{eqnarray*}
g_{\rho_1^*,\rho_2^*}(\xi) &:=& {1\over \max(|\rho^*_1-\rho_2^*|,|\rho_2^*-\rho_2^{*2}|)} \big[\rho_1^*(e^{i\eta \xi}-1) - (e^{i\eta f_{\rho_2^*}(\xi)}-1)\big] \\
& = & {1\over \max(|\rho^*_1-\rho_2^*|,|\rho_2^*-\rho_2^{*2}|)} \big[(\rho_1^*-\rho_2^*)(e^{i\eta \xi}-1) +\rho_2^* (e^{i\eta \xi}-1)- (e^{i\eta f_{\rho_2^*}(\xi)}-1)\big] \\
&=&\sum_{n=1}^\infty a_n(\rho_1^*,\rho_2^*) \xi^n,
\end{eqnarray*}
where $a_n(\rho_1^*,\rho_2^*)$ are the Taylor coefficients of $g_{\rho_1^*,\rho_2^*}$ at 0. 

Fix a constant $r>0$ small enough. 
The sum of the second and third terms in the last brackets can be seen as a holomorphic function in $\rho_2^*$ and $\xi$ with 
$|\rho_2^*|<3$ and $\xi\in 2r\D$. Moreover, this function vanishes when $\rho_2^*=0$ or  $\rho_2^*=1$. We easily deduce that
$$\Pc:=\big\{g_{\rho_1^*,\rho_2^*}(\xi), \quad \rho_1^*,\rho_2^*\in \C, |\rho_2^*|\leq 2, \rho_2^*\not=0,1\big\}$$
is a normal family of holomorphic functions in $\xi\in 2r\D$. Note that this family does not depend on $\alpha,\beta,\theta$ and $\lambda$.

\smallskip\noindent
{\bf Claim.} No sequence in $\Pc$ converges to the zero function. 

\smallskip

Taking into the account the claim, we first complete the proof of the lemma. We show that there is $M>0$ such that all $g_{\rho_1^*,\rho_2^*}$ in $\Pc$ admit at most $M$ zeros in $r\D$, counting multiplicity. Assume by contradiction that there is a sequence of functions $g_n$ in $\Pc$ such that $g_n$ has at least $n$ zeros in $r\D$. By the claim, taking a subsequence allows us to assume that $g_n$ converges locally uniformly on $2r\D$ to a non-zero function $g$. By the classical Hurwitz's theorem $g$ has infinitely many of zeros in $r\overline\D$ which is not possible.

We have shown that if $(\zeta^*,\check\zeta^*)$ is a point in $Z_{\alpha,\beta,\theta}^\lambda$ then $Z_{\alpha,\beta,\theta}^\lambda$ admits not more than $M$ points in the ball of center  $(\zeta^*,\check\zeta^*)$ and of radius $r$, counting multiplicity. It is not difficult to deduce that $Z_{\alpha,\beta,\theta}^\lambda$ is $N$-sparse for some constant $N>0$ depending only on $M$ and $r$. The lemma is then proved.

It remains to verify the claim. 
Assume by contradiction that there is a sequence in $\Pc$ converging to 0. We only consider $(\rho_1^*,\rho_2^*)$ in that sequence.
In particular, for each $n$, the coefficient $a_n(\rho_1^*,\rho_2^*)$ tends to 0. 
Using a direct computation, we obtain the following Taylor approximations of order 2 of $f_{\rho_2^*}(\xi)$ $$f_{\rho_2^*}(\xi) \approx -i\rho_2^*(e^{i\xi}-1) +{i\over 2} [\rho_2^*(e^{i\xi}-1)]^2  \approx \rho_2^*\xi + {i\over 2}(\rho_2^*-\rho_2^{*2})\xi^2.$$
We then deduce that
$$a_1(\rho_1^*,\rho_2^*) = i\eta{\rho_1^*-\rho_2^*\over \max(|\rho^*_1-\rho_2^*|,|\rho_2^*-\rho_2^{*2}|)}$$
and
$$2a_2(\rho_1^*,\rho_2^*) = \eta{(-\eta\rho_1^*+\rho_2^*) +(\eta-1)\rho^{*2}_2\over \max(|\rho^*_1-\rho_2^*|,|\rho_2^*-\rho_2^{*2}|)} 
=i\eta a_1(\rho_1^*,\rho_2^*) +\eta(1-\eta) {\rho_2^* -\rho_2^{*2}\over \max(|\rho^*_1-\rho_2^*|,|\rho_2^*-\rho_2^{*2}|)}  \cdot$$

Clearly, $a_1(\rho_1^*,\rho_2^*)$ and $a_2(\rho_1^*,\rho_2^*)$ cannot tend to 0 together. This is a contradiction which ends the proof of the lemma.
\endproof

Let $0<\epsilon<1$ be a constant small enough whose value will be specified later. 
We divide the set $Z_{\alpha,\beta,\theta}^\lambda$ of solutions of 
\eqref{e:system-Z} into three disjoint subsets 
$$Z^{A,\epsilon,\lambda}_{\alpha,\beta,\theta}, \quad Z^{B,\epsilon,\lambda}_{\alpha,\beta,\theta} \quad \text{and} \quad Z^{C,\epsilon,\lambda}_{\alpha,\beta,\theta}$$ 
corresponding to the following conditions (see also \eqref{e:ratio}):

\begin{equation*} \label{e:system-Z-3}
\text{(A)} \ \  \left\lbrace 
\begin{array}{l}
 |\rho_1-1| \leq \epsilon \\
|\rho_2-1|\leq \epsilon
\end{array} \right.
\quad 
\text{(B)} \ \ \left\lbrace 
\begin{array}{l}
\text{either } |\rho_1-1| \leq  \epsilon  \text{ and } |\rho_2-1| > \epsilon \\
\text{or } \ \  \quad |\rho_1-1| >  \epsilon  \text{ and } |\rho_2-1| \leq  \epsilon
\end{array} \right.
\quad 
\text{(C)} \ \  \left\lbrace 
\begin{array}{l}
 |\rho_1-1| > \epsilon \\
|\rho_2-1| > \epsilon.
\end{array} \right.
\end{equation*}
In  the following lemmas, we  will use the notations such as the sector $\S$ and the half-lines  $Q, \Lambda_{1,s}, \Lambda_{2,s}$ introduced in Appendix \ref{a:Harnack}, see the discussion after Lemma \ref{l:Poisson}
and Definition \ref{d:dominated}. For every set $X$ denote by $\Delta_X$ the diagonal of $X\times X$.

Fix a constant $N>1$ large enough which depends only on $\eta$. We will only require it to satisfy Lemma \ref{l:Z-A} below and we only consider $\epsilon$ such that $0<\epsilon\ll N^{-1}e^{-N}$. Consider also the following condition
\begin{gather*}
\text{(AA)} \qquad  |\delta_{\alpha,\beta}|\leq N\epsilon \quad\text{with}\quad \delta_{\alpha,\beta}:=(\arg\alpha -\arg \beta) -i(\log|\alpha|-\log|\beta|)-2n\pi+2m\eta\pi  \\
 \text{for some} \quad n,m\in\Z.
 \end{gather*}
For simplicity, we choose the values of $\arg\alpha$ and $\arg\beta$ in $[0,2\pi)$. Note that since $\alpha, \beta$ belong to $\A$, by considering the imaginary and real parts of $\delta_{\alpha,\beta}$, we see that if $n,m$ exist, they are unique and both $|n|, |m|$ are bounded by a constant independent of $\alpha,\beta,\epsilon$. 
When the condition (AA) fails, we define
$$Q':=\varnothing, \quad \Lambda_{1,s}':=\varnothing \quad \text{and} \quad \Lambda_{2,s}':=\varnothing.$$

Define also
$$s':=s +\log |\delta_{\alpha,\beta}|+N\leq s.$$
When the condition (AA) holds, we set
$$Q':=Q \quad \text{and} \quad \text{for} \quad s'<0 \quad \Lambda_{1,s}'=\Lambda_{2,s}':=\varnothing, \quad \text{for} \quad s'\geq 0 \quad \Lambda_{1,s}':=\Lambda_{1,s'}, \quad  \Lambda_{2,s}':=\Lambda_{2,s'}.$$

\begin{lemma} \label{l:Z-A}
There are constants $N>0$ and $\kappa>0$, independent of $\epsilon$, such that when $\epsilon$ is small enough, the set 
 $Z^{A,\epsilon,\lambda}_{\alpha,\beta,\theta}$ is $\kappa$-dominated by $\Delta_Q\cup \Delta_{\Lambda_{1,s}'}\cup \Delta_{\Lambda_{2,s}'}$ for all $\alpha,\beta\in\A$, $\theta\in\Theta$ and $\lambda=e^s>1$. 
\end{lemma}
\proof
We only consider points $(\zeta,\check\zeta)$ in the set $Z^{A,\epsilon,\lambda}_{\alpha,\beta,\theta}$ which satisfy \eqref{e:system-Z}
and the above condition (A). We first study the dependence of $Z^{A,\epsilon,\lambda}_{\alpha,\beta,\theta}$ on $\epsilon$ 
when $\epsilon$ goes to 0.
So consider $\epsilon$ small enough and tending to 0. All constants we use are independent of $\epsilon$.

Observe that $\rho_1-1$ and $\rho_2-1$ can be expressed in terms of $\zeta-\check\zeta$ using \eqref{e:ratio}. Then, 
using  the fact that $e^t-1\approx t$ for $|t|$ small, the condition (A) gives us the following estimates for some
integers $n'$ and $m'$ 
      \begin{equation} \label{e:zeta-check-zeta}
       \begin{cases}
   \eta\big( \zeta-\check\zeta+{\log|\alpha|-\log|\beta|\over b} \big) +(\arg\alpha -\arg \beta) -i(\log|\alpha|-\log|\beta|)-2n'\pi=O(\epsilon)\\
      \zeta-\check\zeta+{\log|\alpha|-\log|\beta|\over b}-2m'\pi=O(\epsilon).
       \end{cases}
      \end{equation}
Taking a suitable linear combination of these equations gives
\begin{equation} \label{e:alpha-beta}
(\arg\alpha -\arg \beta) -i(\log|\alpha|-\log|\beta|)-2n'\pi+2m'\eta\pi=O(\epsilon).
\end{equation}

Observe that \eqref{e:alpha-beta} cannot be true if $\alpha,\beta$ do not satisfy the condition (AA). We used here that $N$ is large and $\epsilon$ is small. 
In other words, the set $Z^{A,\epsilon,\lambda}_{\alpha,\beta,\theta}$ is empty when the condition (AA) fails. Clearly, the lemma is true in that case. From now on, assume that the condition  (AA) is satisfied.
By considering the real and imaginary parts of the left hand side of \eqref{e:alpha-beta} and of $\delta_{\alpha,\beta}$ in (AA), we obtain that
 $n'=n$ and $m'=m$. Since $|m|$, $|n|$ are bounded,
it follows from \eqref{e:zeta-check-zeta} and \eqref{e:alpha-beta} that 
$|\zeta-\check\zeta|$ is bounded by a constant.
So, in order to complete the proof of the lemma, we only need to show that the distance between $\check\zeta$ and $Q'\cup\Lambda_{1,s}'\cup\Lambda_{2,s}'$ is bounded by a constant. 

As in \eqref{e:zeta-check-zeta}, we obtain
   \begin{equation} \label{e:rho1-rho2}
       \begin{cases}
      \rho_1-1\approx \zeta-\check\zeta+{\log|\alpha|-\log|\beta|\over b}-2m\pi \\
   \rho_2-1 \approx \eta\big( \zeta-\check\zeta+{\log|\alpha|-\log|\beta|\over b} -2m\pi \big) + \delta_{\alpha,\beta}.
       \end{cases}
      \end{equation}
It follows that one of the following three quantities is bounded below and above by positive constants
$${|\rho_1-1|\over |\rho_2-1|} \CommaBin \quad {|\rho_1-1| \over |\delta_{\alpha,\beta}|} \quad \text{and} \quad {|\rho_2-1| \over |\delta_{\alpha,\beta}|}\cdot$$

We consider separately the three cases corresponding to the last three quantities. 
In the first case, using \eqref{e:system-Z}, we get that $|y_1/y_2|$ is bounded from below and above by positive constants. Therefore, 
$|\Re(i(\eta-1)\check\zeta)|$ is bounded by a constant, or equivalently, the distance between $\check\zeta$ and $Q$ is bounded by a constant. So the lemma is true in this case. 

In the second case, the first equation in  \eqref{e:system-Z} implies that $\lambda |y_1| |\delta_{\alpha,\beta}|$ is bounded from below and above by positive constants. Therefore,  we obtain
\begin{equation} \label{e:s-delta}
|\Re(i\eta \check\zeta) +s +\log |\delta_{\alpha,\beta}||\leq c
\end{equation}
for some constant $c$. Observe that $\Re(i\eta \check\zeta)\leq 0$ for $\check\zeta\in\S$. Therefore,
when $s'<0$, $s +\log |\delta_{\alpha,\beta}|$ is negative with a large absolute value. So the inequality \eqref{e:s-delta} is not satisfied for any $\check\zeta\in\S$ and this second case does not occur. Assume now that $s'\geq 0$.
Recall that the equation of $\Lambda_{2,s}'$ is $\Re(i\eta \check\zeta)+s'=0$, see Appendix \ref{a:Harnack}. So, by \eqref{e:s-delta}, the distance between $\check\zeta$ and $\Lambda_{2,s}'$ is bounded by a constant. The lemma is then true as well.

The last case can be treated in the same way as for the second case. We obtain from the second equation in 
\eqref{e:system-Z} that $\lambda |y_2| |\delta_{\alpha,\beta}|$ is bounded from below and above by positive constants. It follows that
\begin{equation*} 
|\Re(i \check\zeta) +s +\log |\delta_{\alpha,\beta}||\leq c
\end{equation*}
for some constant $c$. We conclude as above that the distance between $\check\zeta$ and $\Lambda_{1,s}'$ is bounded by a constant. This ends the proof of the lemma.
\endproof

Recall that we only consider $\epsilon$ small enough. Define
$$\A^{2,\epsilon} :=\big\{ (\alpha,\beta)\in \A^2 \text{ satisfying the condition (AA)} \big\}.$$

\begin{lemma} \label{l:Z-A-sum}
There is a constant $c>0$ independent of $\epsilon$ such that for $\mu$-almost every $\alpha,\beta\in\A$ and for $s=\log\lambda \geq 1$, we have 
$$\sup_{\theta\in\Theta}\Ebf \Big(\sum_{(\zeta,\check\zeta)\in Z^{A,\epsilon,\lambda}_{\alpha,\beta,\theta}} H_\alpha(\zeta)H_\beta(\check\zeta) \Big)\leq c.$$
Moreover, the last expectation vanishes when $(\alpha,\beta)$ is outside $\A^{2,\epsilon}$ and we have
$$\lim_{\epsilon\to 0} (\mu\otimes\mu)(\A^{2,\epsilon})=0.$$
\end{lemma}
\proof
By Lemma \ref{l:Z-A}, it is clear that the considered expectation vanishes when $(\alpha,\beta)$ is outside $\A^{2,\epsilon}$. 
So we will only consider the case where $(\alpha,\beta)$ is inside $\A^{2,\epsilon}$. 
By Lemma \ref{l:Z-A} again, we can divide $Z^{A,\epsilon,\lambda}_{\alpha,\beta,\theta}$ into three sets $Z,Z_1,Z_2$ which are $\kappa$-dominated by  $\Delta_{Q'}$, $\Delta_{\Lambda_{1,s}'}$ and $\Delta_{\Lambda_{2,s}'}$, respectively. We prove now the first assertion. 
It is enough to prove similar estimates for $Z,Z_1$ and $Z_2$ instead of $Z^{A,\epsilon,\lambda}_{\alpha,\beta,\theta}$. 
The following estimates are uniform on $\theta$.

By Lemmas \ref{l:Z-separated} and \ref{l:sum-line}, we have for some constant $c>0$
$$\sum_{(\zeta,\check\zeta)\in Z} H_\alpha(\zeta)H_\beta(\check\zeta) \leq c \quad \text{and hence} \quad 
\Ebf\Big(\sum_{(\zeta,\check\zeta)\in Z} H_\alpha(\zeta)H_\beta(\check\zeta)\Big) \leq c.$$
This is the desired estimate for $Z$. Now, we only consider the case of $Z_1$ because the case of $Z_2$ can be obtained in the same way.

By Lemmas \ref{l:Z-separated}, \ref{l:sum-line} and \eqref{e:G}, we have  
$$\sum_{(\zeta,\check\zeta)\in Z_1} H_\alpha(\zeta)H_\beta(\check\zeta) \lesssim G_{1,\alpha}(s').$$
Recall that the last sum vanishes when $s'<0$. 
Since $s'$ is equal to $s$ plus a constant (depending on $\alpha,\beta$) and $s'\leq s$, we deduce from the last inequality that
$$\Ebf\Big(\sum_{(\zeta,\check\zeta)\in Z_1} H_\alpha(\zeta)H_\beta(\check\zeta)\Big) \lesssim \Ebf(G_{1,\alpha}(s)).$$
By Lemma \ref{l:G-alpha-int}, the last expectation is bounded. This ends the proof of the first assertion in the lemma.

It remains to prove the last assertion in the lemma. Consider $(\alpha,\beta)$ in $\A^{2,\epsilon}$. 
By using the imaginary part of $\delta_{\alpha,\beta}$, the above condition (AA) implies that 
$$|(\log|\alpha|-\log|\beta|)-2mb\pi| \leq N\epsilon.$$
Since $\alpha$ and $\beta$ are in $\A$, we deduce from the last inequality and the definition of $\A$ in Appendix \ref{a:Harnack} that one of the following inequalities holds (these inequalities correspond to $m=1$, $m=-1$ and $m=0$) 
\begin{equation} \label{e:A-2epsilon}
|\alpha|\leq e^{N\epsilon}e^{-2\pi b}, \quad  |\beta|\leq e^{N\epsilon}e^{-2\pi b}  \quad  \text{and} \quad |\log|\alpha|-\log|\beta||\leq N\epsilon.
\end{equation}

Consider the first two inequalities. Observe that when $\epsilon$ goes to 0, the two sets 
$$\big\{(\alpha,\beta)\in \A^2,\ |\alpha|\leq e^{N\epsilon}e^{-2\pi b}\big\} \quad \text{and} \quad \big\{(\alpha,\beta)\in \A^2,\ |\beta|\leq e^{N\epsilon}e^{-2\pi b}\big\}$$
tend to the empty set. So their $\mu\otimes\mu$ measures tend to 0. Therefore, we only need to consider now the set of $(\alpha,\beta)$ in $\A^{2,\epsilon}$ satisfying the last inequality in \eqref{e:A-2epsilon}. Note that in this case the integer $m$ is necessarily equal to 0.

By using the real part of $\delta_{\alpha,\beta}$ and the condition (AA), we obtain
$$|(\arg\alpha-\arg\beta) -2n\pi|\leq N\epsilon.$$
The set $\widetilde \A^{2,\epsilon}$ of all $(\alpha,\beta)\in \A^2$ satisfying this inequality and the last inequality in  \eqref{e:A-2epsilon} with $m=0$ tends to the diagonal of $\A^2$ when $\epsilon$ goes to 0. As $\mu$ contains no atom, the measure $\mu\otimes \mu$ has no mass on  the diagonal of $\A^2$. Thus, the measure $(\mu\otimes\mu)(\widetilde \A^{2,\epsilon})$ tends to 0 as $\epsilon$ tends to 0. This completes the proof of the lemma.
\endproof

We refer to Appendix \ref{a:Harnack} for the definition of $\Lambda^{\! 0}_{1,s}$ and $\Lambda^{\! 0}_{2,s}$.
Consider the following subsets of $\S\times \S$
$$K_s^1:=\big\{(\zeta,\check\zeta) ,\ \zeta\in \Lambda^{\! 0}_{1,s} \text{ and } \check\zeta\in \zeta -\overline\eta \, \R_{\geq 0} \big\}, \quad 
\widecheck K_s^1:=\big\{(\zeta,\check\zeta),\ \check\zeta\in \Lambda^{\! 0}_{1,s} \text{ and } \zeta \in \check\zeta - \overline \eta \, \R_{\geq 0}\big\}$$
and
$$K_s^2:=\big\{(\zeta,\check\zeta),\ \zeta\in \Lambda^{\! 0}_{2,s} \text{ and } \check\zeta \in \zeta +\R_{\geq 0}\big\}, \quad 
\widecheck K_s^2:=\big\{(\zeta,\check\zeta),\ \check\zeta\in \Lambda^{\! 0}_{2,s} \text{ and } \zeta \in \check\zeta +\R_{\geq 0} \big\}.$$
The constants used in the following lemmas may depend on $\epsilon$.

\begin{lemma} \label{l:Z-B}
For every $0<\epsilon<1$, there is a constant $\kappa_\epsilon>0$ such that for $s=\log\lambda$ with $\lambda$ large enough, the set 
 $Z^{B,\epsilon,\lambda}_{\alpha,\beta,\theta}$ is $\kappa_\epsilon$-dominated by $K^1_s\cup \widecheck K^1_s\cup K_2^s\cup \widecheck K^2_s$ for all $\alpha,\beta\in\A$ and $\theta\in\Theta$. 
 \end{lemma}
\proof
Consider $(\zeta,\check\zeta)$ in $Z^{B,\epsilon,\lambda}_{\alpha,\beta,\theta}$. So the above condition (B) is satisfied. 
For simplicity, we assume that the second line in (B) holds. The case where the first line holds  can be treated in the same way.
We deduce that $\rho_2$ is bounded from above and below by positive constants. 
Then, by  \eqref{e:ratio}, we have that $|\Im(\zeta-\check\zeta)|$ is bounded by a constant.
Furthermore, by the second line of \eqref{e:system-Z}, both $\lambda |x_2|$ and $\lambda |y_2|$ are bounded from below by a positive constant. Thus, $\Im(\zeta)$ and $\Im(\check\zeta)$ are bounded from above by $s$ plus a constant.

Now, we have either $|x_1|\geq |y_1|$ or $|x_1|\leq |y_1|$. We only consider the first case as the second one can be obtained in the same way.
We deduce from the inequality $|x_1|\geq |y_1|$ and \eqref{e:ratio} that $\Re(i\eta(\zeta-\check\zeta))$ is bounded from below by a constant. 
Recall that $i\eta=ia-b$ with $b>0$. 
Since $|\Im(\zeta-\check\zeta)|$ is bounded by a constant, we easily deduce that $\Re(\check\zeta)$ is larger than $\Re(\zeta)$ minus a constant. It follows that the distance from $\check\zeta$ to $\zeta+\R_{\geq 0}$ is bounded by a constant.

Since $|\rho_1-1|>\epsilon$ and  $|x_1|\geq |y_1|$, from the first line of  \eqref{e:system-Z}, we obtain that $|x_1|$ is bounded from below and above by positive constants times $\lambda^{-1}$. 
Since $s=\log\lambda$, we deduce that $|\Re(i\eta\zeta)+s|$ is bounded by a constant, see also \eqref{e:x-zeta-bis}. 
Recall that $\Re(i\eta\zeta)+s=0$ is the equation of the real line containing $\Lambda_{2,s}$. 
So, the distance from $\zeta$ to this line is bounded. 
Since $\zeta$ is in $\S$ and $\Im(\zeta)$ is smaller than $s$ plus a constant, we deduce that the distance from $\zeta$ to $\Lambda^{\! 0}_{2,s}$ is bounded by a constant, see the discussion after Lemma \ref{l:Poisson}. 
Now, it is not difficult to conclude that the distance from $(\zeta,\check\zeta)$ to $K_s^2$ is bounded by a constant. This ends the proof of the lemma.
\endproof

\begin{lemma} \label{l:Z-B-sum}
For every fixed $0<\epsilon<1$, there is a constant $c_\epsilon>0$ such that for $\mu$-almost every $\alpha,\beta\in\A$, we have 
$$\sup_{\theta\in\Theta}\Ebf \Big(\sum_{(\zeta,\check\zeta)\in Z^{B,\epsilon,\lambda}_{\alpha,\beta,\theta}} H_\alpha(\zeta)H_\beta(\check\zeta) \Big)\leq c_\epsilon \quad \text{and} \quad \lim_{s\to\infty} \sup_{\theta\in\Theta}\Ebf \Big(\sum_{(\zeta,\check\zeta)\in Z^{B,\epsilon,\lambda}_{\alpha,\beta,\theta}} H_\alpha(\zeta)H_\beta(\check\zeta) \Big) =0.$$
\end{lemma}
\proof
By Lemma \ref{l:Z-B}, we can divide $Z^{B,\epsilon,\lambda}_{\alpha,\beta,\theta}$ into 4 disjoint subsets $Z^1, \widecheck Z^1,
Z^2, \widecheck Z^2$ which are respectively $\kappa_\epsilon$-dominated by $K^1_{s}, \widecheck K^1_{s}, K_2^{s}, \widecheck K^2_{s}$. 
It is enough to show the properties similar to the ones in the lemma but for the sets $Z^1, \widecheck Z^1,
Z^2, \widecheck Z^2$ instead of 
$Z^{B,\epsilon,\lambda}_{\alpha,\beta,\theta}$. For simplicity, we only consider the case of $Z^1$. The other cases can be treated in the same way.

Consider the following lattice of $K^1_s$
$$Z^1_s:=\big\{ (-\overline \eta b^{-1} s + m, -\overline \eta b^{-1}s + m- \overline \eta n), \text{ with } n,m\in \N, m\leq b^{-1}s \big\}.$$
Since $K^1_s$ is $(|\eta|+1)$-dominated by $Z^1_s$, the set $Z^1$ is $(\kappa_\epsilon + |\eta|+1)$-dominated by $Z^1_s$. 
By using Lemmas \ref{l:sum-comparison}, \ref{l:sum-line-0}, and then Lemma \ref{l:axe-sum} (applied to $\hbar:=b^{-1}$ and $s:=mb$) together with \eqref{e:G},  we obtain
\begin{eqnarray*}
\sum_{(\zeta,\check\zeta)\in Z^1} H_\alpha(\zeta)H_\beta(\check\zeta) & \lesssim & \sum_{(\zeta,\check\zeta)\in Z^1_s} H_\alpha(\zeta)H_\beta(\check\zeta) \\
& = & \sum_{0\leq m\leq b^{-1}s} H_\alpha(-\overline \eta b^{-1}s + m)\sum_{n\in\N} H_\beta(-\overline\eta b^{-1}s + m- \overline\eta n) \\
&\lesssim & \sum_{0\leq m\leq b^{-1}s} H_\alpha(-\overline\eta b^{-1}s + m)\int_{l\geq b^{-1}s} H_\beta(m-\overline\eta l)dl  \\ 
&\lesssim & \sum_{0\leq m\leq b^{-1}s} H_\alpha(-\overline\eta b^{-1}s + m)\int_{l\geq m} H_\beta(m-\overline\eta l)dl  \\ 
&\lesssim & \sum_{0\leq m\leq b^{-1}s} H_\alpha(-\overline\eta b^{-1}s + m) \ \lesssim \  G_{1,\alpha}(s).
\end{eqnarray*}
Now, Lemma \ref{l:G-alpha-int} implies the desired properties.
\endproof

We continue to refer to Appendix \ref{a:Harnack} for the notations such as $Q, Q_s^{\infty}$ and $\zeta_s$.

\begin{lemma} \label{l:Z-C}
Let $0<\epsilon<1$ be any fixed constant. Then, there is a constant $\kappa_\epsilon>0$ 
such that for $s=\log\lambda$ with $\lambda$ large enough the following property holds for all $\alpha,\beta\in\A$ and $\theta\in\Theta$. 
There are positive numbers $s_1,s_2,s_3,s_4$ (which may depend on $\alpha,\beta,\theta,\lambda$) such that
 the set $Z^{C,\epsilon,\lambda}_{\alpha,\beta,\theta}$ is $\kappa_\epsilon$-dominated by the union of the following $10$ sets
 $$\Lambda_{1,s}^{\! \infty} \times \Lambda_{2,s}, \quad \Lambda_{1,s} \times \Lambda_{2,s}^{\! \infty}, \quad \Lambda_{2,s}^{\! \infty} \times \Lambda_{1,s} \quad \Lambda_{2,s} \times \Lambda_{1,s}^{\! \infty}$$
 and
 $$\{\zeta_s\}\times Q_s^{\infty}, \quad  \{\zeta_s\}\times \Lambda^{\! \infty}_{1,s_1}, \quad \{\zeta_s\}\times  \Lambda^{\! \infty}_{2,s_2}, \quad Q_s^{\infty}\times \{\zeta_s\}, \quad \Lambda^{\! \infty}_{1,s_3}\times \{\zeta_s\}, \quad \Lambda^{\! \infty}_{2,s_4}\times \{\zeta_s\}.$$
 \end{lemma}
\proof
We only consider $(\zeta,\check\zeta)$ in $Z^{C,\epsilon,\lambda}_{\alpha,\beta,\theta}$. So they satisfy \eqref{e:system-Z} and the condition (C) above.
For simplicity, we assume that $|\rho_1|\leq 1$ because the opposite case can be treated in the same way.
We fix a constant $r>0$ small enough, depending only on $\eta$. Fix also a constant $\kappa_\epsilon$ big enough depending on $\epsilon$ and $r$. 

We deduce from the first equation of \eqref{e:system-Z}, Condition (C) and the inequality $|\rho_1|\leq 1$ that $|\lambda y_1|$ is bounded from below and above by positive constants. 
Therefore, $\Re(i\eta \check\zeta) +\log \lambda$ is bounded from below and above. Since $s=\log\lambda$ and $\check\zeta$ is in $\S$, the distance from $\check\zeta$ to $\Lambda_{2,s}$ is bounded. 

\medskip\noindent
{\bf Case 1.} Assume that $|\rho_2|>r$. We obtain  from the second equation of \eqref{e:system-Z} and Condition (C) that $|\lambda x_2|$ is bounded from below and above 
by positive constants. It follows that $|\Im \zeta -s|$ is bounded by a constant. So the distance between $\zeta$ and $\Lambda_{1,s}$ is bounded by a constant. Thus, $(\zeta,\check\zeta)$ is $\kappa_\epsilon$-dominated by 
$\Lambda_{1,s}\times \Lambda_{2,s}$ for a suitable choice of $\kappa_\epsilon$. 
Moreover, by \eqref{e:ratio}, using $|\rho_2|>r$, we also obtain that $\Im(\check\zeta)$ is larger than $\Im(\zeta)$ minus a constant. Thus, $\Im(\check\zeta)$ is larger than $s$ minus a constant. We conclude that $(\zeta,\check\zeta)$ is $\kappa_\epsilon$-dominated by 
$\Lambda_{1,s}\times \Lambda_{2,s}^{\!\infty}$.

\medskip\noindent
{\bf Case 2.}  Assume that  $|\rho_2|\leq r$. This, Condition (C) and the second equation of \eqref{e:system-Z} imply that $|\lambda y_2|$ is bounded from below and above 
by positive constants. It follows that the distance between $\check\zeta$ and $\Lambda_{1,s}$ is bounded by a constant. We conclude that 
the distance between $\check\zeta$ to $\zeta_s$, which is  the intersection of $\Lambda_{1,s}$ with $\Lambda_{2,s},$ is bounded by a constant.

\medskip\noindent
{\bf Case 2a.} Assume that $|\rho_1|> r$. As above, the first equation of \eqref{e:system-Z} and Condition (C) imply that $|\lambda x_1|$ is bounded from below and above by positive constants and hence
$(\zeta,\check\zeta)$ is $\kappa_\epsilon$-dominated by $\Lambda_{2,s}\times\{\zeta_s\}$ and hence by $\Lambda_{2,s}\times\Lambda_{1,s}^{\!\infty}$.

\medskip\noindent
{\bf Case 2b.} Assume that $|\rho_1|\leq r$.
So, from now on, we only consider $(\zeta,\check\zeta)$ satisfying \eqref{e:system-Z} and the two inequalities $|\rho_1|\leq r$ and $|\rho_2|\leq r$. 
By \eqref{e:system-Z}, both $|\lambda x_1|$ and $|\lambda x_2|$ are bounded from above by a small constant.
Arguing as above, we deduce that  $\zeta$ belongs to $\zeta_s+\S$.

We know that the distance between $\check\zeta$ to $\zeta_s$ is bounded by a constant. 
If $|\rho_1|>r|\rho_2|$ and $|\rho_2|>r|\rho_1|$, by considering $\rho_1/\rho_2$, we deduce from \eqref{e:ratio} that $|\Re(i(\eta-1)(\zeta-\check\zeta))|$ is bounded by some constant. It follows that $|\Re(i(\eta-1)(\zeta-\zeta_s))|$ satisfies the same property. Hence, the distance between 
$\zeta-\zeta_s$ to the real line $\widetilde Q$ containing $Q$ is bounded. Since $\zeta_s$ belongs to $Q$, the distance between $\zeta$ and $\widetilde Q$ is bounded. As $\zeta$ belongs to $\zeta_s+\S$, we see that $\zeta$ is $\kappa_\epsilon$-dominated by $Q_s^\infty$.
It remains to consider the cases where $|\rho_1|\leq r |\rho_2|$ or $|\rho_2|\leq r |\rho_1|$. We only study the first case as the second one can be treated in the same way.

Denote by $Z$ the set of all $(\zeta,\check\zeta)$ satisfying \eqref{e:system-Z} and the inequalities $|\rho_1|\leq r$, $|\rho_2|\leq r$, $|\rho_1|\leq r|\rho_2|$. The inequality  $\rho_1\leq r\rho_2$ and \eqref{e:ratio} imply that $|e^{i(\eta-1)(\zeta-\check\zeta)}|$ is small and hence 
$|e^{i(\eta-1)(\zeta-\zeta_s)}|$ is small as well. The last number is equal to $|e^{i(\eta-1)\zeta}|$ because $\zeta_s$ belongs to $Q$. Hence, $\zeta$ is in the angle limited by $Q$ and $\R_{\geq 0}$, see Appendix \ref{a:Harnack}.

\medskip\noindent
{\bf Claim. } $Z$ is $\kappa_\epsilon$-dominated by $\Lambda\times \{\zeta_s\}$ for some real line $\Lambda$ on the upper half-plane which is parallel to $\R$.

\medskip

Clearly, the claim implies the lemma. Indeed, the intersection of $\Lambda$ with $\S$ is equal to $\Lambda_{1,s_3}$ for some positive number $s_3$. Since $\zeta$ is in the angle limited by $Q$ and $\R_{\geq 0}$, the claim implies that $(\zeta,\check\zeta)$ is 
$\kappa_\epsilon$-dominated by $\Lambda^{\!\infty}_{1,s_3} \times\{\zeta_s\}$ (we increase the value of $\kappa_\epsilon$ if necessary). This is the desired property. So, it remains to prove the claim.

We can assume that $Z$ contains at least two points since otherwise the claim is obvious.
Let $(\zeta^*,\check\zeta^*)$ be a point in $Z$ such that the distance from $\zeta^*$ to the edge $-\overline\eta \R_{\geq 0}$ of $\S$ is smaller than the infimum of all such  distances plus a small positive constant. 
We also denote by $\rho_1^*, \rho_2^*$ the corresponding values of $\rho_1,\rho_2$ for the chosen point of $Z$. As in the proof of Lemma \ref{l:Z-separated}, any $(\zeta,\check\zeta)$ in $Z$ satisfies the equations in \eqref{e:system-Z-local} with $\xi:=\zeta-\zeta^*$ and $\check\xi:=\check\zeta-\check\zeta^*$.

By the choice of $(\zeta^*,\check\zeta^*)$, we have that $|e^{i\eta\xi}|\leq 2$. Recall that the above constant $r$ and the constant 
$\epsilon_0$ in \eqref{e:Theta} are small. Therefore, the inequalities $|\rho_1|\leq r$ and $|\rho_2|\leq r$ imply that $(y_1, y_2)$ is very close to $(-\theta_1/\lambda,-\theta_2/\lambda).$ 
 It follows that
$\check\zeta$ is very close to $\zeta_s$. In particular, this also holds for $\check\zeta^*$. We then deduce that $\check\xi$ is small. This and \eqref{e:system-Z-local} imply that 
$${\rho_1^*(e^{i\eta\xi}-1) \over \rho_2^*(e^{i\xi}-1)} = {(e^{i\eta\check\xi}-1) \over (e^{i\check\xi}-1)}\approx \eta.$$
The first numerator is small in comparison with $\rho_2^*$ because  $\rho_1^*\leq r\rho^*_2$ and  $|e^{i\eta\xi}|\leq 2$. Hence, $|e^{i\xi}-1|$ should be small.
It follows that $\Im(\xi)$ is bounded from above and below by some constants. We conclude that  $\zeta$ has a bounded distance to the real line $\Lambda$  passing through $\zeta^*$ and parallel to $\R$. 
Thus,  $(\zeta,\check\zeta)$ has a bounded distance to $\Lambda\times \{\zeta_s\}$.
This ends the proof of the lemma.
\endproof

\begin{lemma} \label{l:Z-C-sum}
For every fixed $0<\epsilon<1$, there is a constant $c_\epsilon>0$ such that for $\mu$-almost every $\alpha,\beta\in\A$, we have 
$$\sup_{\theta\in\Theta}\Ebf \Big(\sum_{(\zeta,\check\zeta)\in Z^{C,\epsilon,\lambda}_{\alpha,\beta,\theta}} H_\alpha(\zeta)H_\beta(\check\zeta) \Big)\leq c_\epsilon \quad \text{and} \quad \lim_{s\to\infty} \sup_{\theta\in\Theta}\Ebf \Big(\sum_{(\zeta,\check\zeta)\in Z^{C,\epsilon,\lambda}_{\alpha,\beta,\theta}} H_\alpha(\zeta)H_\beta(\check\zeta) \Big) =0.$$
\end{lemma}
\proof
We can divide $Z^{C,\epsilon,\lambda}_{\alpha,\beta,\theta}$ into 10 disjoint subsets $Z_1,\ldots,Z_{10}$ which are respectively $\kappa_\epsilon$-dominated by the 10 sets in Lemma \ref{l:Z-C}. 
It is enough to prove the properties similar to the ones in the lemma  
for each $Z_i$ instead of $Z^{C,\epsilon,\lambda}_{\alpha,\beta,\theta}$. 
We only consider the cases where $i=1,5,6$ because the other cases can be obtained in the same way.
The estimates below are uniform on $\theta\in\Theta$.

For $Z_1$, we will use Lemma \ref{l:Z-separated} and a suitable lattice in $\Lambda_{1,s}^{\!\infty}\times\Lambda_{2,s}$ as in Lemma \ref{l:axe-sum}. After that, 
by using Lemma \ref{l:sum-line-0} and then Lemma \ref{l:axe-sum} (for $\hbar=1$) and \eqref{e:G}, we obtain
$$\sum_{(\zeta,\check\zeta)\in Z_1} H_\alpha(\zeta)H_\beta(\check\zeta)\lesssim \Big(\int_0^\infty H_\alpha(\zeta_s+l)dl\Big)G_{2,\beta}(s)\lesssim G_{2,\beta}(s).$$
Thus, Lemma \ref{l:G-alpha-int} gives us the desired property for $Z_1$.

For $Z_5$, observe that if $(\zeta,\check\zeta)$ is in $Z_5$ then the distance between $\zeta$ and $\zeta_s$ is bounded by $\kappa_\epsilon$. By Harnack's inequality, $H_\alpha(\zeta)$ is bounded by a constant times $H_\alpha(\zeta_s)$. Therefore, 
using the second assertion of Lemma \ref{l:sum-line-0} for $\beta$ instead of $\alpha$, we obtain
$$\sum_{(\zeta,\check\zeta)\in Z_5} H_\alpha(\zeta)H_\beta(\check\zeta)\lesssim H_\alpha(\zeta_s) .$$
So the desired property for $Z_5$ follows from the second assertion of Lemma \ref{l:G-alpha-int}.

Finally, for $Z_6$, arguing as above, using the first inequality in Lemma \ref{l:axe-sum} for $\beta$ instead of $\alpha$, we have
$$\sum_{(\zeta,\check\zeta)\in Z_6} H_\alpha(\zeta)H_\beta(\check\zeta)\lesssim H_\alpha(\zeta_s) \int_{l\geq 0}H_\beta(\zeta_s+l)dl\lesssim H_\alpha(\zeta_s).$$
We then obtain the result by using again  the second assertion of Lemma \ref{l:G-alpha-int}.
\endproof

\proof[End of the proof of Proposition \ref{p:theta-slice}]
By \eqref{e:vartheta}, we have
$$\int_{\alpha,\beta\in\A}\Ebf\big(\|\vartheta_{\alpha,\beta,\theta}^\lambda\|\big) d\mu(\alpha)d\mu(\beta)= \int_{\alpha,\beta\in\A}\Ebf\Big( \sum_{(\zeta,\check\zeta)\in Z_{\alpha,\beta,\theta}^{\lambda}} H_\alpha(x_1,x_2)H_\beta(y_1,y_2)\Big)d\mu(\alpha)d\mu(\beta).$$
We can split the last expression into the sum of the following three terms
$$\int_{\alpha,\beta\in\A} \Ebf\Big(\sum_{(\zeta,\check\zeta)\in Z_{\alpha,\beta,\theta}^{A,\epsilon,\lambda}}\Big)+ \int_{\alpha,\beta\in\A}\Ebf\Big(\sum_{(\zeta,\check\zeta)\in Z_{\alpha,\beta,\theta}^{B,\epsilon,\lambda}}\Big)+ \int_{\alpha,\beta\in\A}\Ebf\Big(\sum_{(\zeta,\check\zeta)\in Z_{\alpha,\beta,\theta}^{C,\epsilon,\lambda}}\Big).$$
By Lemmas \ref{l:Z-A-sum}, \ref{l:Z-B-sum} and \ref{l:Z-C-sum}, when $\epsilon$ is fixed, all these three terms are bounded by a constant independent of $\alpha,\beta, \theta$ and the last two terms tend to 0, uniformly on $\theta$, when $\lambda$ tends to infinity. 
This already gives us the estimate in the proposition.
Moreover, given any $\delta>0$, by taking $\epsilon$ small enough, the last assertion in Lemma \ref{l:Z-A-sum} shows that all limit values of the first term are smaller than $\delta$. The second property in the proposition follows easily.
\endproof


\begin{appendix}

\section{Young's inequality and applications} \label{a:Young}

In this appendix, we recall the classical Young's inequality for integral operators. 
We apply this inequality in the charts of $X\times X$ which cover the diagonal $\Delta$.

Let $k(x, y)$ be a function on  $\B \times \B$,  smooth in $(\B \times \B) \setminus \Delta$. 
Assume that there is a constant $c>0$ and a number $\delta\geq 0$  such that for every $(x,y)\in \B\times \B,$
\begin{equation} \label{e:Young}
 \|k(x,\cdot)\|_{L^{1+\delta}}\leq c\quad\text{and}\quad \|k(\cdot, y)\|_{L^{1+\delta}}\leq c.
\end{equation}
Here, we use  the norm $L^p$ with respect to the normalized Lebesgue measure on $\B$.

Define a linear operator
$P$ on the  space of  measures $\mu$ of bounded mass on $\B$ by
$$(P \mu)(x) := \int_{y\in \B}k(x, y) d\mu(y).$$
We are also interested in the case where $\mu$ is given by an $L^p$ function.

\begin{lemma}[Young's inequality] \label{L:Young}
The operator $P$ maps continuously measures of bounded mass into  $L^{1+\delta}(\B),$  $L^p(\B)$ into $L^q(\B)$, and  $L^\infty$ into $\Cc^0$;  all with norm bounded by $c,$ 
where $q =\infty $ if $p^{-1} + (1 + \delta)^{-1}\leq 1$ and $p^{-1} + (1 + \delta)^{-1} = 1 + q^{ -1}$ otherwise.
\end{lemma}

We list here two examples of kernels needed in our study.

\begin{example}\rm \label{Ex:kernel_1}
  Consider a   kernel $k(x,y)$ of modulus bounded by some constant times
$\|x - y\|^{ -2}$. In this  case,  we can choose any $0\leq \delta<1.$
\end{example}

\begin{example} \rm \label{Ex:kernel_2}
Consider  a family of convolution kernels  
$$k_r(x,y)= r^{-4}g_r(x,y) \textbf{1}_{\{\|x-y\|<r\}},$$
where $\textbf{1}_{\{\|x-y\|<r\}}$ is the  characteristic  function of the set $\{\|x-y\|<r\}\cap (\B\times \B)$
and  $(g_r)$ is a uniformly bounded family of functions. Consider $\delta=0$ and the operator 
$P_r$ with kernel $k_r$. It maps $L^p(\B)$ to itself with norm bounded by a constant independent of $r.$ 
 \end{example}
 
Consider now a family $(K_\lambda)$ of smooth $4$-forms on $X\times X$ depending on a parameter $\lambda\in\C$ with $|\lambda|$ larger than a positive constant. Assume that there is a constant $A>0$ such that $K_\lambda(x,y)$ vanishes when the distance between $x$ and $y$ is larger than $A|\lambda|^{-1}$.

\begin{lemma} \label{l:kernel-Delta}
Assume that $\|K_\lambda\|_\infty =O(|\lambda|^4)$ and that $K_\lambda$ converges weakly to $c[\Delta]$ as $\lambda$ tends to infinity, where $c$ is a constant. Then, for all $2$-forms $f_1$ and $f_2$ of class $L^2$, we have
$$\lim_{\lambda\to\infty} \langle f_1\otimes f_2, K_\lambda\rangle = c\langle f_1,f_2\rangle.$$
\end{lemma}
\proof
Define the integral operator $P_\lambda$ associated to $K_\lambda$ by 
$$P_\lambda(f)(y):=\int_x K_\lambda(x,y) f(x)$$
for all $2$-forms $f$ on $X$. Observe that $P_\lambda(f)$ is also a $2$-form and we have
$$ \langle f_1\otimes f_2, K_\lambda\rangle = \langle f_2, P_\lambda(f_1)\rangle.$$

By hypothesis on the support of $K_\lambda$ and its sup-norm, in local coordinates, the coefficients of $K_\lambda$ satisfy estimates in  \eqref{e:Young} for $\delta=0$.
By Lemma \ref{L:Young} for $\delta=0$, the operator $P_\lambda$ from $L^2$ to $L^2$ has a norm bounded independently of $\lambda$. Therefore, in order to obtain the result, we can assume that $f_1$ is smooth because smooth forms are dense in the space of $L^2$ forms. 
Similarly, we can also assume that $f_2$ is smooth. Now,
by hypothesis, $P_\lambda(f_1)$ converges weakly to $cf_1$ and the result follows easily.
\endproof

\begin{lemma} \label{l:kernel-Delta-bis}
Assume that $\|K_\lambda\|_\infty =O(|\lambda|^3)$. Then we have 
$$\lim_{\lambda\to\infty} \langle f_1\otimes f_2, K_\lambda\rangle = 0$$
if $f_1$ is of class $L^q$ with  $q>4/3$ and $f_2$ is of class $L^2$.
\end{lemma}
\proof
By hypothesis, $K_\lambda$ tends to 0 in $L^1$ when $\lambda$ tends to infinity. Moreover,
in local coordinates, we can check that the coefficients of $K_\lambda$ satisfy estimates in \eqref{e:Young} for all $0\leq \delta<1/3$. We obtain the result exactly as in the last lemma using that $P_\lambda(f_1)$ has a bounded $L^2$ norm, thanks to Lemma \ref{L:Young}.
\endproof

\section{Some properties of $\ddc$-closed currents} \label{a:current}

We recall   some basic notions and properties on positive $\ddc$-closed currents on a complex surface and refer the reader
to \cite{BerndtssonSibony, DinhSibony04, Skoda} for details. 

Let $T$ be a positive $\ddc$-closed $(1,1)$-current on $X$ and let $x$ be a local coordinate system around a point $a$ of $X$.  
It is well-known that such a current gives no mass
to sets of zero Hausdorff 2-dimensional measure, see e.g. \cite[p.\,389]{BerndtssonSibony}.
Define
$$\nu(T,a,r):={1\over\pi r^{2}}\int_{\B(a,r)} T\wedge \ddc \|x\|^2 \quad \text{and} \quad \overline\nu(T,a,r):={1\over\pi r^{2}}\int_{\overline{\B(a,r)}} T\wedge \ddc \|x\|^2.$$
By Skoda \cite{Skoda}, the function
$r\mapsto \nu(T,a,r)$ is  increasing and the limit
\begin{equation}\label{e:Lelong}
\nu(T,a):=\lim_{r\to0+} \nu(T,a,r)=\lim_{r\to0+} \overline\nu(T,a,r)
\end{equation}
is  a non-negative finite  number which is called the {\it Lelong number} of $T$ at $a$. 
Indeed, thanks to Lelong-Jensen identity \cite[Prop.\,1]{Skoda}, we have
\begin{equation} \label{e:Jensen}
\nu(T,a,r)-\nu(T,a) = 2\int_{\B(a,r)\setminus\{a\}} T\wedge \ddc\log\|x\|.
\end{equation}

It is known that the notion of  Lelong number  does not depend on the choice of  local holomorphic coordinates $x$. 
Moreover, it follows from the definition that the functions $a\mapsto \overline\nu(T,a,r)$ and  $a\mapsto \nu(T,a)$ are upper-semi-continuous. 
We have the following result.

\begin{lemma} \label{l:lelong}
Let $T$ be a positive $\ddc$-closed current of mass $1$ on $X$. Then there is a constant $c>0$ such that 
$$\nu(T,x,r) \leq c \quad \text{and} \quad \nu(T,x)\leq c \quad \text{for} \quad \|x\|\leq 5 \quad \text{and} \quad r\leq 4.$$
\end{lemma}
\proof
This is a direct consequence of the inequality $\nu(T,x,r)\leq \overline\nu(T,x,r)$
upper-semi-continuity of  $x\mapsto \overline\nu(T,x,r)$ and the monotone dependence of $\nu(T,x,r)$ on $r$. 
\endproof

The following result was obtained in  \cite{FornaessSibony05}, see also Proposition \ref{p:T-rep} below.

\begin{lemma} \label{l:pot}
Let $T$ be a positive $\ddc$-closed current on $X$. Then
it can be represented as 
\begin{equation} \label{e:decompo_posi_har}
T=\Omega+\partial S+\overline{\partial S} + i\ddbar u
\end{equation}
with $\Omega$ a smooth real closed $(1,1)$-form, $S$ a current of bi-degree $(0,1)$ and $u$ a real function in $L^p$ for $p<2.$ 
Moreover, for every such a representation, the currents $\dbar S$ and $\partial\overline S$ do not depend on the choice of $\Omega,S,u$ and they are
forms of class $L^2$, uniquely determined by $T$. 
\end{lemma}
\proof See  \cite[Prop. 2.6, 2.7 and Thm. 2.9]{FornaessSibony05}.
\endproof

Note that the representation \eqref{e:decompo_posi_har} is not unique but the uniqueness of $\dbar S$ and $\partial \overline S$ allows us to define the {\it energy} $E(T)$ of $T$ as
\begin{equation} \label{e:energy}
E(T):=\int_X\dbar S\wedge \partial\overline S.
\end{equation}
This is a non-negative number which is independent of the choice of $\Omega, S$ and $u$.
It is not difficult to see that $E(T)=0$ if and only if $\dbar S=0$ and if and only if $T$ is closed, see \cite{FornaessSibony05} for details.

Consider a local coordinate system $x=(x_1,x_2)$ in $X$ with $|x_1|<3$ and $|x_2|<3$. Then for almost every $x_2\in 3\D$ the slice 
$T\wedge [(3\D)\times \{x_2\}]$ exists and is a positive measure, see \cite[Th.\,1.18]{Bassanelli}. Denote by $\vartheta_{x_2}$ the restriction of this measure to the disc $(2\D)\times \{x_2\}$. We have the following lemma.

\begin{lemma} \label{l:T-mass-slice}
The mass $m(x_2)$ of  $\vartheta_{x_2}$ is an $L^p$ function in $x_2\in\ 2\D$ for all $1\leq p\leq 2$. 
\end{lemma}
\proof
It is enough to consider the case where $p=2$.
Let $0\leq\chi\leq 1$ be a smooth function on $(3\D)\times (3\D)$ such that $\chi=1$ when $|x_1|\leq 2$ and $\chi=0$ for $|x_1|>5/2$. If $\Phi$ denotes the projection $(x_1,x_2)\mapsto x_2$, then the function $m(x_2)$ satisfies
$$m \leq \Phi_*(\chi T).$$
So, it is enough to prove that the function $\widetilde m:=\Phi_*(\chi T)$ is in $L^2(2\D)$. 

Using the above representation of $T$ and the fact that $i\ddbar T=0$, we have
\begin{eqnarray*}
i\ddbar \widetilde m & = & \Phi_*(i\ddbar\chi\wedge T) + \Phi_*(i\partial\chi\wedge \dbar T) -  \Phi_*(i\dbar\chi\wedge \partial T) \\
& = & \Phi_*(i\ddbar\chi\wedge T) - \Phi_*(i\partial\chi\wedge \ddbar S) +  \Phi_*(i\dbar\chi\wedge \ddbar \overline S).
\end{eqnarray*}
The first term in the last sum is a measure of finite mass. The two other terms belong to the Sobolev space $H^{-1}(3\D)$ because $\dbar S$ and $\partial \overline S$ are in $L^2$. 
So we can write the last sum as $\mu^+-\mu^-+h$, where $\mu^\pm$ are positive measures of finite mass on $3\D$ and $h$ is a distribution in $H^{-1}(3\D)$. Solving the following Laplace's equations
$$i\ddbar \phi^\pm=\mu^\pm \quad \text{and} \quad i\ddbar\phi=h$$
gives us two subharmonic functions $\phi^\pm$ on $3\D$ and a locally $L^2$ function $\phi$ on $3\D$, see e.g. \cite[p.\,355]{Taylor}; indeed, $\phi$ is a locally $H^1$ function. 

Observe now that both functions $\widetilde m$ and $\phi^+-\phi^-+\phi$ satisfy the same Laplace's equation
$$i\ddbar \widetilde m = \mu^+-\mu^-+h \quad \text{and} \quad i\ddbar (\phi^+-\phi^-+\phi) = \mu^+-\mu^-+h.$$
Therefore, their difference $\widetilde m-(\phi^+ - \phi^-+\phi)$ is a harmonic function. Recall that harmonic and subharmonic functions are locally $L^2$ functions. 
So, we easily deduce from the above discussion that $\widetilde m$ is in $L^2(2\D)$. This ends the proof of the lemma.
\endproof

Using the last lemma, we obtain the following result.

\begin{proposition} \label{p:T-rep}
There is a representation as in \eqref{e:decompo_posi_har} such that all currents 
$S,\overline S, \partial S, \partial\overline S, \dbar S$, $\dbar\overline S$ are forms of class $L^2$ and $u, \partial u, \dbar u$ are functions or forms  of class $L^p$ for every $1\leq p<2$. 
\end{proposition}
\proof
It was shown in \cite{FornaessSibony05} that there is such a representation with $S,\overline S, \partial S, \partial\overline S, \dbar S, \dbar\overline S$ in $L^2$ and $u$ in $L^p$ for every $1\leq p<2$. Consider such a representation. With the above notations, it is enough to show that 
$\partial u$ belongs to $L^p(\D\times\D)$ as its complex conjugate $\dbar u$ should satisfy the same property as well. We will only show that
$\partial u/\partial x_1$ is in $L^p(\D\times\D)$ because the same proof works for  $\partial u/\partial x_2$.

We deduce from \eqref{e:decompo_posi_har} that 
$$i\ddbar u = R \quad \text{with} \quad R:=T-\Omega-\partial S-\dbar\overline S.$$
For almost every $x_2\in\D$, the slice $R\wedge [(2\D)\times \{x_2\}]$ exists and is a measure of finite mass. 
Denote by $R_{x_2}$ this measure and by $n(x_2)$ its mass.
Since $\Omega$ is smooth and $\partial S,\dbar\overline S$ are of class $L^2$, we deduce from Lemma \ref{l:T-mass-slice} that $n(x_2)$ is an $L^2$ function (and hence, an $L^p$ function for $1\leq p<2$) on $\D$. 

Consider the following function
$$v(x_1,x_2):= {1\over \pi} \int \log |x_1-\zeta| d R_{x_2}(\zeta) \quad \text{for } x_1\in 2\D \text{ and } x_2\in\D.$$
For each fixed $x_2$, this is the standard logarithmic potential of $R_{x_2}$. It is not difficult to see that there is a constant $c_p>0$ depending only on $p$ such that for each fixed $x_2$
$$\|v\|_{L^p(2\D)} \leq c_p n(x_2) \quad \text{and} \quad \Big\|{\partial v\over \partial x_1}\Big\|_{L^p(\D)}\leq c_p n(x_2).$$
Since $n(x_2)$ is an $L^p$ function, we deduce that $v$ is a function in $L^p((2\D)\times\D)$ and  $\partial v/\partial x_1$ is in $L^p(\D\times\D)$. 
In particular, $u-v$ belongs to $L^p((2\D)\times\D)$ because $u$ is an $L^p$ function. 

Observe now that when $x_2$ is fixed, both $u$ and $v$ satisfy the same Laplace's equation 
$$ (i\ddbar)_{x_1} u = R_{x_2} \quad \text{and} \quad (i\ddbar)_{x_1} v =R_{x_2}.$$ 
We deduce that $u-v$ is harmonic in $x_1$. In particular, there is a constant $c'_p>0$ depending only on $p$ such that for each fixed $x_2\in\D$ we have
$$\Big\|{\partial (u-v)\over \partial x_1} \Big\|_{L^p(\D)} \leq c'_p \|u-v\|_{L^p(2\D)}.$$

Finally, since $u-v$ is in $L^p((2\D)\times\D)$, we deduce from the last estimate that $\partial (u-v)/\partial x_1$ belongs to $L^p(\D\times\D)$. It follows that $\partial u/\partial x_1$ also belongs to $L^p(\D\times \D)$ because we have seen that $\partial v/\partial x_1$ satisfies the same property. This ends the proof of the proposition.
\endproof

\section{Directed $\ddc$-closed currents and Harnack's inequality} \label{a:Harnack}

Let $T$ be a positive $\ddc$-closed $(1,1)$-current directed by $\Fc$ which is a foliation with only hyperbolic singularities or a bi-Lipschitz lamination. 
Assume that $T$ has no mass on every  single leaf of $\Fc$. The local description of $T$ on a regular flow box is given at the beginning of the Introduction. 
It also holds in the case of a bi-Lipschitz lamination. We now discuss the case of a singular flow box, see also \cite{DinhSibony18,FornaessSibony10,Nguyen18a}.

Let $p$ be a hyperbolic singular point of the foliation. 
So, there are local coordinates $x = (x_1 , x_2 )$ centered at $p$
such that in the bidisc $(3\D)^2 := \{|x_1 | < 3, |x_2 | < 3\},$ the foliation $\Fc$ is
defined by the form
$$x_2dx_1-\eta x_1dx_2$$
for some complex number  $\eta=a+ib$ with $a,b\in\R$ and $b\not=0.$ Note that
if we flip $x_1$ and $x_2,$ then $\eta$ is changed to $1/\eta = \overline\eta/|\eta|^2 = a/(a^2 + b^2 ) - ib/(a^2 + b^2 )$. Therefore, we can assume from now on that the axes are chosen so that $b > 0$. 

Observe that the two axes of the
bidisc $(3\D)^2$ are invariant and are the separatrices of the foliation in the
bidisc $(3\D)^2$. Consider the ring $\A$ defined by
$$\A := \left\lbrace  \alpha\in \C,\ e^{-2\pi b} < |\alpha| \leq 1 \right\rbrace.$$
Define also the sectors $\S$ and $\S'$  by
$$\S := \big\lbrace \zeta = u + iv \in \C,\ v > 0\quad\text{and}\quad bu + av > 0\big\rbrace$$
and
$$\S' :=\big\lbrace \zeta = u + iv \in \C,\ v > - \log 3\quad \text{and}\quad bu + av > - \log 3\big\rbrace .$$
Note that the sector $\S$ is spanned by the vectors 
$1, -\overline \eta$, or equivalently, by $1,-\eta^{-1}$ because $\overline\eta =(a^2+b^2)\eta^{-1}$. Moreover, $\S$ is contained in  the upper half-plane $\H := \{u+iv,\ v > 0\}$ and in the sector $\S'$. 
The
angle of $\S$ is $\arctan{(-b/a)}\in (0, \pi )$ and the boundary $b\S$ of $\S$ is formed by two half-lines starting from 0: one is spanned by $-\overline\eta$ (or $-\eta^{-1}$) and the other is $\R_+$ which is spanned by 1.

For $\alpha\in  \C^* ,$ consider the following manifold $\mathcal L_\alpha$ immersed in $\C^2$ and defined by
\begin{equation} \label{e:x-zeta} 
x_1 = \alpha e^{i\eta (\zeta +\log {|\alpha|/b})}\quad \text{and}\quad x_2 = e^{i(\zeta +\log{ |\alpha|/b})}\quad \text{with}\quad \zeta = u + iv \in \C.
\end{equation}
So we have 
\begin{equation} \label{e:|z|}
|x_1| =e^{\Re(i\eta\zeta)}= e^{-bu-av}\quad\text{and}\quad    |x_2| =e^{\Re(i\zeta)}= e^{-v}.
\end{equation}
Observe that $v$ and $bu+av$ are equal to constants times the distances from $u+iv$ to the two edges of $\S$.
The map $\zeta \mapsto (x_1 , x_2 )$ is injective because $\eta \not\in\R.$ The following properties are not difficult to check.
\begin{enumerate}
\item $\mathcal L_\alpha$ is tangent to the above vector field and is a submanifold of $\C^{*2} .$
\item  $\mathcal L_{\alpha_1}$ is equal to $\mathcal L_{\alpha_2}$ if $\alpha_1 /\alpha_2 = e^{2ki\pi\eta}$ for some $k\in\Z$ and they are
disjoint otherwise. In particular, $\mathcal L_{\alpha_1}$ and $\mathcal L_{\alpha_2}$ are disjoint if $\alpha_1 , \alpha_ 2 \in \A$
and $\alpha_1 \not = \alpha_2 .$
\item  The union of $\mathcal L_\alpha$ is equal to $\C^{*2}$ for $\alpha \in \C^* ,$ and then also for $\alpha\in  \A.$
\item  The intersection $L_\alpha :=\mathcal  L_\alpha \cap \D^2$ of $\mathcal L_\alpha$ with the unit bidisc $\D^2$ is given by
the same equations as in the definition of $\mathcal L_\alpha$  but with $\zeta \in \S.$ Moreover,
$L_\alpha$ is a connected submanifold of $\D^{*2} .$ In particular, it is a leaf of $\Fc \cap \D^2 .$
\item Similarly, the intersection $L'_\alpha := \mathcal L_\alpha \cap (3\D)^2$ is given by the same
equations with $\zeta \in \S'.$  Moreover, $L'_\alpha$ is a connected submanifold of
$(3\D^*)^2$ and is the leaf of $\Fc \cap (3\D)^2$ which contains $L_\alpha .$
\end{enumerate}

Recall that $T$ is assumed to have no mass on every single leaf of $\Fc$. So, it gives no mass to the separatrices of the singularities and admits the following decomposition.

\begin{lemma}[see {\cite[Lem.\,4.1]{DinhSibony18}}] \label{L:T_T-alpha}
There is a positive measure $\mu$ of finite mass on $\A$, without atoms, and positive 
harmonic functions $h_\alpha$ on $L'_\alpha$ for $\mu$-almost every $\alpha\in \A$   such that we have in $(3\D)^2$
$$T =\int_\A T_\alpha d\mu(\alpha),\quad \text{where}\quad T_\alpha := h_\alpha[L'_\alpha ].$$
Moreover, the mass of $T_\alpha$ in $(2\D)^2$ is $1$ for $\mu$-almost every $\alpha \in\A.$ 
\end{lemma}

Using \eqref{e:x-zeta}, we define
\begin{equation} \label{e:H_alpha}
H_\alpha(\zeta ) := h_\alpha \big( \alpha e^{i\eta(\zeta +\log{ |\alpha|/b)}}, e^{i(\zeta +\log{ |\alpha|/b})}\big).
\end{equation}
This is a positive harmonic function on the sector $\S'$.  
Consider the map 
\begin{equation*} \label{e:Phi-gamma}
\Phi : \zeta \mapsto \zeta^\gamma  \quad \text{with} \quad \gamma := {\pi\over \arctan(-b/a)} >1.
\end{equation*}
It sends $\S$ bi-holomorphically to the upper half-plane $\H$. 
Define the real variables $u,v, U, V$ and the function $\widetilde H_\alpha$ by
\begin{equation*} \label{e:UV}
u+iv:=\zeta, \quad U+iV:=\zeta^\gamma=(u+iv)^\gamma \quad \text{and} \quad \widetilde H_\alpha:=H_\alpha\circ\Phi^{-1}.
\end{equation*}
The function $\widetilde H_\alpha$ is positive harmonic on $\Phi(\S')$ which contains the closed half-plane $\overline \H$. 

\begin{lemma} \label{l:Poisson}
There is a constant $c>0$ such that  for $\mu$-almost every $\alpha\in \A$,  we have the following Poisson formula
$$\widetilde H_\alpha(U+iV)={1\over \pi}\int_{t\in\R} \widetilde H_\alpha (t)
{V\over V^2 + (t - U )^2} dt \quad \text{for } U+iV \text{ in } \H$$
and the estimates
$$\int_{t\in\R} \widetilde H_\alpha  (t)|t|^{ -1+1/\gamma} dt \leq  c \quad \text{and} \quad \widetilde H_\alpha(t) \leq c \quad \text{for } t\in \overline\H.$$
\end{lemma}
\proof
This result was obtained in 
\cite[Lem.\,4.2,\,4.3,\,4.4]{DinhSibony18} and \cite[Prop.\,1]{FornaessSibony10} except that the inequality $\widetilde H_\alpha(t) \leq c$ was proved for $t\in\R$. 
However, the above Poisson formula implies that this inequality still holds for
 $t\in\overline\H$.
\endproof

We now describe some segments and half-lines in $\S$ which play an important role in our study. 
Several of them are parallel to the edges of $\S$.
We consider a parameter $s\geq 0$. 
Let $\Lambda_{1,s}$ denote the half-line, starting from the point $-\overline\eta b^{-1}s = -\eta^{-1} b^{-1}(a^2+b^2)s$ on the boundary of $\S$, which is parallel to the edge $\R_+$ of $\S$. This is the restriction to $\S$ of the real line $is+\R$ which is also the line of equation $\Re(i\zeta)+s=0$.
Denote by $\Lambda_{2,s}$ the half-line starting from the point $b^{-1} s$ on the boundary of $\S$ and parallel to the other edge $-\overline\eta\R_{\geq 0}$ of $\S$, i.e. the edge 
containing the points $-\overline\eta$ and $-\eta^{-1}$. This is the restriction to $\S$ of the real line $b^{-1}s -\overline \eta\R$ which is of equation $\Re(i\eta\zeta)+s=0$.

Define $\zeta_s:= (1-\overline\eta)b^{-1}s$ which is the only intersection point of $\Lambda_{1,s}$ and $\Lambda_{2,s}$. Denote by $Q$ the half-line starting from 0 and passing through $\zeta_s$. It does not depend on $s$. Denote also by $Q_s^{\infty}$ the half-line starting from the point $\zeta_s$ which is contained in $Q$. Note that the equation of $Q$ is $\Re(i(\eta-1)\zeta)=0$ because $\zeta_s$ satisfies this equation. Moreover, the part of $\S$ limited by $Q$ and $\R_{\geq 0}$ is defined by the inequality $\Re(i(\eta-1)\zeta)\leq 0$ because this inequality is true for $\zeta=1$. 
The quantity $|\Re(i(\eta-1)\zeta)|$ is equal to a constant times the distance from $\zeta$ to the real line $\widetilde Q$ containing $Q$.

Finally, the point $\zeta_s$ divides $\Lambda_{1,s}$ into two intervals: the bounded one is denoted by 
$\Lambda^{\! 0}_{1,s}$ and the unbounded one is denoted by $\Lambda^{\! \infty}_{1,s}$. Similarly, the point $\zeta_s$ divides $\Lambda_{2,s}$ into two intervals: the bounded one is denoted by 
$\Lambda^{\! 0}_{2,s}$ and the unbounded one is denoted by $\Lambda^{\! \infty}_{2,s}$.
So, $\Lambda^{\! 0}_{1,s}$ is the segment joining
$-\overline\eta b^{-1}s$ and $\zeta_s$; $\Lambda^{\! 0}_{2,s}$ is the segment joining $b^{-1}s$ and $\zeta_s$. 
The following lemma gives us estimates on some integrals on $\Lambda_{1,s}$ and $\Lambda_{2,s}$, see also \eqref{e:G}.

\begin{lemma} \label{l:axe-sum}
Let $\hbar>0$ be any fixed constant. Then there is a constant $c_\hbar>0$ such that for every $s>0$ and $\mu$-almost every $\alpha\in\A$ we have 
$$\int_{l\geq \hbar s} H_\alpha(-\overline\eta b^{-1}s+l) dl \leq c_\hbar \quad \text{and} \quad 
\int_{l\geq \hbar s}  H_\alpha(b^{-1}s- \overline\eta l)dl \leq c_\hbar.$$
Moreover, we have 
$$\lim_{s\to\infty} \int_{l\geq \hbar s} H_\alpha(-\overline\eta b^{-1}s+l) dl =0 \quad \text{and} \quad 
\lim_{s\to\infty}\int_{l\geq \hbar s}  H_\alpha(b^{-1}s- \overline\eta l)dl =0.$$
\end{lemma}
\proof
We only prove the lemma for $\Lambda_{1,s}$ because the case of $\Lambda_{2,s}$ can be obtained in the same way. 
We prove now the first inequality in the lemma. The constants we use below may depend on $\hbar$.

Write $\zeta:=u+iv:=-\overline\eta b^{-1}s+l$ and consider $U,V$ as above. By the first assertion of Lemma \ref{l:Poisson}, we have  
$$\int_{l\geq \hbar s}  H_\alpha(-\overline\eta b^{-1}s+l) dl = {1\over \pi}\int_{t\in\R} \widetilde H_\alpha(t) \Big(\int_{l\geq\hbar s} {V\over V^2+(t-U)^2} dl \Big) dt.$$
By the second assertion of Lemma \ref{l:Poisson}, it is enough to show that the integral between the parentheses 
is bounded by a constant times $|t|^{-1+1/\gamma}$. 

Observe that $v=s$. 
Define 
$$r:=s^{-1}l,\quad U':=s^{-\gamma} U \quad \text{and} \quad V':=s^{-\gamma} V.$$
Since $l\geq \hbar s$, we have $r\geq \hbar$. According to \cite[Lem.\,5.6]{DinhSibony18}, we have
$$U'=r^\gamma + O(r^{\gamma-1}) \quad \text{and} \quad V'=\gamma r^{\gamma-1}+O(r^{\gamma-2}) \quad \text{as} \quad r\to\infty.$$

We deduce that the above integral between the parentheses is bounded by a constant times 
(we use the variable $R:=s^\gamma |t|^{-1} r^\gamma$)
$$\int_{\hbar}^\infty {s^{\gamma+1} r^{\gamma-1} \over s^{2\gamma} r^{2\gamma-2} +(t-s^\gamma r^\gamma)^2} dr
=\gamma^{-1} |t|^{-1+1/\gamma} \int_{\hbar^\gamma s^\gamma |t|^{-1}}^\infty {s \over s^2 |t|^{-1/\gamma} R^{2-2/\gamma}+ |t|^{1/\gamma}(\pm 1-R)^2} dR.$$
By the estimate in Lemma \ref{l:Poisson}, we only need to show that the last integral is bounded by a constant. 
For this purpose, it is enough to consider the case where the $\pm 1$ in the last line is $1$. Denote the considered integral by 
$I(s,t)$.
We split it into two parts : $I_1(s,t)$ is the integral for $R$ in $[1/2,+\infty)$ and $I_2(s,t)$ is the integral for $R$ such that  $\hbar^\gamma s^\gamma |t|^{-1} \leq R \leq 1/2$.

In order to bound  $I_1(s,t)$, we define $R':=s^{-1}|t|^{1/\gamma}(1-R)$. Then we have
$$I_1(s,t)\lesssim \int_{-\infty}^\infty {s \over s^2 |t|^{-1/\gamma} + |t|^{1/\gamma}(1-R)^2} dR
\leq \int_{-\infty}^\infty {1\over 1 + R'^2} dR'.$$
So $I_1(s,t)$ is bounded by a constant.
For the integral $I_2(s,t)$, observe that the domain of this integral is non-empty only when $|t|\geq 2\hbar^\gamma s^\gamma$. 
So we have $I_2(s,t)=0$ when $|t| < 2\hbar^\gamma s^\gamma$. Moreover, when $|t|\geq 2\hbar^\gamma s^\gamma$, we obtain
$$I_2(s,t)\lesssim \int_0^{1/2} {s\over |t|^{1/\gamma}} dR\leq \int_0^{1/2} {1\over \hbar} dR.$$
Clearly, $I_2(s,t)$ is bounded by a constant as well. This ends the proof of the first inequality in the lemma.

Note that when $\hbar$, $t$ are fixed and $s\geq 1$, the above estimates on $I_1(s,t)$ show that 
$$I_1(s,t) \lesssim \int_{\hbar^\gamma s^\gamma |t|^{-1}}^\infty {s\over s^2+R^2} dR= \int_{\hbar^\gamma s^{\gamma-1} |t|^{-1}}^\infty {1\over 1+R''^2} dR''$$
with $R'':=s^{-1}R$. As $\gamma>1$, we see that $I_1(s,t)$ tends to 0 when $s$ tends to infinity. 
This is one of the instances where we use $\gamma>1$, i.e. the hyperbolicity of the singularities of the foliation.
Since $I_2(s,t)$ vanishes when $s$ is large enough, we obtain that $I(s,t)$ tends to 0 as $s$ tends to infinity. 

On the other hand, we have seen in the above discussion that 
$$\int_{l\geq \hbar s} H_\alpha(-\overline\eta b^{-1}s+l) dl \lesssim \int_{t\in\R} \widetilde H_\alpha(t) |t|^{-1+1/\gamma} I(s,t)dt.$$
Recall that $I(s,t)$ is bounded. Now, we easily deduce the first limit in the lemma from the estimate in Lemma \ref{l:Poisson} and Lebesgue's dominated convergence theorem. This completes the proof of the lemma.
\endproof

For any function or more generally a current $f(s)$, depending on the parameter $s>0$, we denote {\it the expectation} of $f(s)$ on the interval 
$(0,s]$ by $\Ebf(f(s))$. This is the mean value of $f$ on the interval $(0,s]$ which is given by the formula
\begin{equation} \label{e:Expec}
\Ebf(f(s)):=s^{-1}\int_0^s f(\check s)d\check s.
\end{equation}

For $s\geq 0$, consider also the following  integrals of $H_\alpha$ on the half-lines $\Lambda_{1,s}$ and $\Lambda_{2,s}$
\begin{equation} \label{e:G}
G_{1,\alpha}(s):=\int_{l\geq 0} H_\alpha(-\overline\eta b^{-1}s+l)dl \quad \text{and} \quad G_{2,\alpha}(s):=\int_{l\geq 0} H_\alpha(b^{-1}s- \overline \eta l)dl.
\end{equation}
We have the following result.

\begin{lemma} \label{l:G-alpha-int} 
There is a constant $c>0$ such that for $\mu$-almost every $\alpha\in\A$, all $s>0$ and for $i=1,2$, we have 
$$\Ebf(G_{i,\alpha}(s))\leq c \quad \text{and} \quad \lim_{s\to\infty} \Ebf(G_{i,\alpha}(s))=0.$$
Moreover, we have  for $\mu$-almost every $\alpha\in\A$ and all $s>0$
$$\Ebf(H_\alpha(\zeta_s))\leq c \quad \text{and} \quad \lim_{s\to\infty} \Ebf(H_\alpha(\zeta_s))=0.$$
\end{lemma}
\proof
We only prove the lemma for $i=1$ because the case where $i=2$ can be obtained in the same way. 
Consider the first assertion.
Define 
$$G'_{1,\alpha}(s):=\int_{0\leq l\leq s} H_\alpha(-\overline\eta b^{-1}s+l)dl \quad \text{and} \quad 
G''_{1,\alpha}(s):=\int_{l\geq s} H_\alpha(-\overline\eta b^{-1}s+l)dl.$$
By Lemma \ref{l:axe-sum} applied to $\hbar=1$, we obtain the same properties as in the first assertion of the lemma for $G''_{1,\alpha}$ instead for $G_{1,\alpha}$. So we only need to prove such properties for $G'_{1,\alpha}$.

We use the same notations as in the proof of Lemma \ref{l:axe-sum} but here, since we consider $0\leq l\leq s$, we have 
$0\leq r\leq 1$. Define also $t':=s^{-\gamma} t$.
According to \cite[Lem\,5.5]{DinhSibony18}, for some constants $\rho>0, \beta>0$ and $c>0$ depending only on $\eta$, we have
$$U'=-\rho+O(r),\quad V'= \beta r+O(r^2)\quad \text{and} \quad V'^2+(t'-U')^2 \geq c [r^2+(\rho+t')^2].$$
As in Lemma \ref{l:axe-sum}, we get 
\begin{eqnarray*}
G'_{1,\alpha}(s) & \lesssim & \int_{t\in\R} \widetilde H_\alpha(t) \Big( s^{1-\gamma} \int_{0<r<1} {r\over r^2 +(\rho+t')^2} dr\Big) dt \\
& \lesssim & \int_{t\in\R} \widetilde H_\alpha(t) \Big( s^{1-\gamma} \log{1+(\rho+t')^2 \over (\rho+t')^2} \Big) dt \\
& \lesssim & \int_{t\in\R} \widetilde H_\alpha(t) \Big( s^{1-\gamma} \log\Big[1+{1 \over (\rho-|t'|)^2} \Big] \Big) dt.
\end{eqnarray*}
Recall that $t'=s^{-\gamma}t$. By using $s_*:=|t|^{-1/\gamma}s$, we obtain
$$G'_{1,\alpha}(s) \lesssim \int_{t\in\R} \widetilde H_\gamma(t) |t|^{-1+1/\gamma} g(s_*) dt \quad \text{with} \quad g(s_*):=  s_*^{1-\gamma} \log\Big[1+{s_*^{2\gamma} \over (\rho s_*^\gamma - 1)^2} \Big].$$
Since $s_*$ depends linearly on $s$, it follows that
$$\Ebf(G'_{1,\alpha}(s)) \leq  \int_{t\in\R} \widetilde H_\gamma(t) |t|^{-1+1/\gamma} \Ebf(g(s_*)) dt  \quad \text{with} \quad  
\Ebf(g(s_*)):=s_*^{-1} \int_0^{s_*} g(\check s_*) d\check s_*.$$

Now, observe that $g(s_*)$ tends to 0 when $s_*$ tends to 0 or infinity. Moreover, $g(s_*)$ has a unique singularity at the point $\rho^{-1/\gamma}$ which is a logarithmic singularity. Therefore, $\Ebf(g(s_*))$ is a bounded continuous function tending to 0 when $s_*$ tends to 0 or infinity.
We apply now the estimate in Lemma \ref{l:Poisson} and Lebesgue's dominated convergence theorem. It is not difficult to obtain 
that $\Ebf(G'_{1,\alpha}(s))$ is bounded by a constant and tends to 0 when $s$ tends to infinity.
This completes the proof of the first assertion.

Consider now the second assertion. 
We apply Harnack's inequality to positive harmonic functions on the sector $\S'$ which contains $\S$. So there is a constant $\kappa \geq 1$ such that  $\mu$-almost every $\alpha$, 
we have $H_\alpha (\zeta_s)\leq \kappa H_\alpha(\zeta)$ when $|\zeta-\zeta_s|\leq 1$. It follows from \eqref{e:G}  that 
$$H_\alpha (\zeta_s)\leq \kappa G_{1,\alpha}(s).$$
So the second assertion is a consequence of the first one by replacing $c$ with $\kappa c$.
\endproof

We need the following lemma in order to estimate some integrals on the half-line $Q$.

\begin{lemma} \label{l:diag-sum}
Let $\zeta$ be any fixed point in the interior of the angle $\S$. Then there is a constant $c_\zeta>0$ such that for $\mu$-almost every  $\alpha\in \A$ we have
$$\int_0^\infty H_\alpha (l\zeta) dl \leq c_\zeta.$$
\end{lemma}
\proof 
We will use the above notations with $\zeta=u+iv$. 
Note that the constants we use in this lemma may depend on $u,v$ or equivalently on $U,V$. 
Since $u+iv$ is in the interior of $\S$, we have $V>0$. 
The integral in the lemma is equal to
$${1\over \pi}\int_{t\in\R} \widetilde H_\alpha (t) \Big[\int_0^\infty  {l^\gamma V\over l^{2\gamma} V^2 + (t - l^\gamma U )^2} dl\Big] dt.$$
Observe that $l^{2\gamma} V^2 + (t - l^\gamma U )^2$ is larger than a positive constant times $l^{2\gamma}+t^2$. This is easy to see by considering $|t|>2l^\gamma U$ and $|t|\leq 2l^\gamma U$.
We then deduce that the integral in the above brackets is smaller than a constant times
$$\int_0^\infty  {l^{\gamma} \over l^{2\gamma} + t^2} dl =\gamma^{-1} |t|^{-1+1/\gamma} \int_0^\infty {\widetilde l^{1/\gamma}\over \widetilde l^2+1} d\widetilde l,$$
where we use the new variable $\widetilde l$ with $l^\gamma=|t|\widetilde l$. Since the last integral is finite, we easily deduce the first estimate in the lemma from the integral estimate in Lemma \ref{l:Poisson}.
\endproof

We will describe some applications of Harnack's inequality which allow us to estimate some infinite sums used in our computation.

\begin{definition} \rm \label{d:dominated}
Let $Z$ and $Z'$ be two subsets of $\R^n$, where the points are counted with multiplicity. 
We say that $Z$ is {\it $N$-sparse} for some constant $N>0$ if any open ball of radius 1 in $\R^n$ contains at most $N$ points of $Z$ counted with multiplicity. 
We say that $Z$ is {\it $\kappa$-dominated by $Z'$} for some constant $\kappa>0$  if the distance between $A$ and $Z'$ is less than $\kappa$ for every point $A$ in $Z$.
\end{definition}

Note that $Z$ is $\kappa$-dominated by $Z'$ if and only if each point of $Z$ is $\kappa$-dominated by $Z'$.

\begin{lemma} \label{l:sum-comparison-0} 
Let $Z$ be an $N$-sparse subset of $\S$ which is $\kappa$-dominated by another subset $Z'$ of $\S$. Then 
there is a constant $c_{N,\kappa}>0$ independent of $Z$ and $Z'$ such that for $\mu$-almost every $\alpha \in\A$, we have
$$\sum_{\zeta\in Z} H_\alpha(\zeta) \leq c_{N,\kappa} \sum_{\zeta'\in Z'} H_\alpha(\zeta').$$
\end{lemma}
\proof
This lemma can be proved using the same arguments as in the next lemma which is slightly more complicated. The details are left to the reader.
\endproof

\begin{lemma} \label{l:sum-comparison}
Let $Z$ be an $N$-sparse subset of $\S\times\S$ which is $\kappa$-dominated by another subset $Z'$ of $\S\times \S$. Then 
there is a constant $c_{N,\kappa}>0$ independent of $Z$ and $Z'$ such that for $\mu$-almost every   $\alpha,\beta\in\A$, we have
$$\sum_{(\zeta,\xi)\in Z} H_\alpha(\zeta)H_\beta(\xi) \leq c_{N,\kappa} \sum_{(\zeta',\xi')\in Z'} H_\alpha(\zeta')H_\beta(\xi').$$
\end{lemma}
\proof
By hypotheses, the balls $B_{\zeta',\xi'}$ of center $(\zeta',\xi')\in Z'$ and radius $\kappa$ cover the set $Z$. Moreover, since $Z$ is $N$-sparse, the cardinality of $B_{\zeta',\xi'}\cap Z$ is bounded by some constant $N'$ which only depends on $N$ and $\kappa$. On the other hand, by Harnack's inequality, there is a constant $c>0$ independent of $Z,Z',\alpha,\beta, \zeta',\xi'$ such that 
$$H_\alpha(\zeta)\leq c H_\alpha(\zeta') \quad \text{and} \quad H_\beta(\xi)\leq c H_\beta(\xi')\quad \text{for all } 
(\zeta,\xi)\in B_{\zeta',\xi'}\cap (\S\times\S).$$
We easily deduce the lemma by taking $c_{N,\kappa}:= N'c^2$.
\endproof

In the same way, we obtain the following results.

\begin{lemma} \label{l:sum-line-0}
Let  $\zeta_0$ and $\xi_0$ be two points in $\S$ with $\xi_0 \not= 0$. Let $Z$ be any $N$-sparse subset of $\S$ which is $\kappa$-dominated by the half-line $L:=\zeta_0+\xi_0\R_{\geq 0}$. Then, there is a constant $c_{N,\kappa}>0$ independent of $\zeta_0,\xi_0,Z$ such that for $\mu$-almost every $\alpha\in\A$, we have
$$\sum_{\zeta\in Z} H_\alpha(\zeta) \leq c_{N,\kappa} |\xi_0| \int_{l\geq 0} H_\alpha(\zeta_0+l\xi_0)dl.$$
In particular, if $L$ is a half-line starting from $0$ in the interior of $\S$ (i.e. $\zeta_0=0$ and $\xi_0$ is in the interior of $\S$), then there is a constant $c_{N,\kappa,L}>0$ independent of $Z$ such that
for $\mu$-almost every $\alpha\in\A$, we have
$$\sum_{\zeta\in Z} H_\alpha(\zeta) \leq c_{N,\kappa,L}.$$
\end{lemma}
\proof
Using the change of variable $l=:|\xi_0|^{-1} l'$, we can assume that $|\xi_0|=1$. 
Observe that the second assertion is a consequence of the first one and Lemma \ref{l:diag-sum} applied to $\zeta:=\xi_0$. It remains to prove the first assertion.

By Lemma \ref{l:sum-comparison-0}, we can assume that $Z$ is the subset $\zeta_0+\xi_0\N$ of the half-line $L$. 
By Harnack's inequality, there is a constant $c>0$ such that $H_\alpha(\zeta)\leq c H_\alpha(\xi)$ 
for $\zeta,\xi\in\ \S$ with $|\zeta-\xi|\leq 1$.  It follows that 
$$H_\alpha(\zeta_0+n\xi_0)  \leq c\int_n^{n+1} H_\alpha(\zeta_0+l\xi_0)  dl$$
for every $n\in\N$. Therefore, we have
$$\sum_{\zeta\in Z} H_\alpha(\zeta)  \leq c \int_0^\infty H_\alpha(\zeta_0+l\xi_0) dl.$$
This ends the proof of the lemma.
\endproof

\begin{lemma} \label{l:sum-line}
Let $L$ be a half-line as in Lemma \ref{l:sum-line-0} and let $\Delta_L$ denote  the diagonal of $L\times L$. Let $Z$ be any $N$-sparse subset of $\S\times\S$ which is $\kappa$-dominated by $\Delta_L$. Then there is a constant $c_{N,\kappa}>0$  independent of $\zeta_0,\xi_0,Z$ such that for $\mu$-almost every $\alpha\in\A$, we have
$$\sum_{(\zeta,\xi)\in Z} H_\alpha(\zeta) H_\beta(\xi) \leq c_{N,\kappa} |\xi_0| \int_{l\geq 0} H_\alpha(\zeta_0+l\xi_0)dl.$$
Moreover, if $L$ is  a half-line starting from $0$ in the interior of $\S$, then there is a constant $c_{N,\kappa,L}>0$ independent of $Z$ such that
for $\mu$-almost every $\alpha,\beta\in\A$, we have
$$\sum_{(\zeta,\xi)\in Z} H_\alpha(\zeta) H_\beta(\xi) \leq c_{N,\kappa,L}.$$
\end{lemma}
\proof
We can assume that $|\xi_0|=1$.
By Lemma \ref{l:sum-comparison}, we can replace $Z$ by the diagonal $Z'$ of the set $(\zeta_0+\xi_0\N)\times (\zeta_0+\xi_0\N)$  because $Z$ is 
 $(\kappa+2)$-dominated by $Z'$. By Lemma \ref{l:Poisson} applied for $\beta$ instead of $\alpha$, we have that  $H_\beta(\xi)$ is bounded by 1. Therefore, the lemma is a direct consequence of Lemma \ref{l:sum-line-0}.
\endproof

\end{appendix}


\end{document}